\documentclass[11pt]{amsart}

\usepackage{epigamath}
\hypersetup{linkcolor=black,citecolor=black,filecolor=black,urlcolor=black}

\usepackage[english]{babel}

\numberwithin{equation}{section}

\usepackage{amsmath, amsthm, amssymb} 
\usepackage{epigraph}
\usepackage{mathtools}
\usepackage{longtable} 
\usepackage{pifont} 
\usepackage{scalerel}
\usepackage{array}
\usepackage{lipsum}
\usepackage{scrextend}
\usepackage{mathrsfs} 
\usepackage{verbatim}

\usepackage[shortlabels]{enumitem}
\setlist[enumerate,1]{label={\rm(\arabic*)},ref={\rm\arabic*}}

\usepackage{tikz-cd} 
\newcommand*{\Scale}[2][4]{\scalebox{#1}{$#2$}}%

\usepackage[all]{xy}

\newtheorem{thm}{Theorem}[section]
\newtheorem{proposition}[thm]{Proposition}
\newtheorem{corollary}[thm]{Corollary}
\newtheorem{lem}[thm]{Lemma}
\newtheorem{conjecture}[thm]{Conjecture}

\theoremstyle{definition}
\newtheorem{defin}[thm]{Definition}

\theoremstyle{remark}
\newtheorem{remark}[thm]{Remark}
\newtheorem{setting}[thm]{Setting} 
\newtheorem{question}[thm]{Question}
\newtheorem{example}[thm]{Example}

\newcommand{\OOO}{\mathscr{O}}

\let\emptyset\varnothing

\newcommand{\pr}{\mathrm{pr}}

\newcommand{\lra}{\longrightarrow}
\DeclareMathOperator{\Cr}{Cr}
\DeclareMathOperator{\Cl}{Cl}
\DeclareMathOperator{\Aut}{Aut}
\DeclareMathOperator{\GL}{GL}
\DeclareMathOperator{\PGL}{PGL}
\DeclareMathOperator{\Pic}{Pic}
\DeclareMathOperator{\Bir}{Bir}
\DeclareMathOperator{\Sing}{Sing}

\DeclareMathOperator{\Exc}{Exc}
\DeclareMathOperator{\mult}{mult}
\DeclareMathOperator{\rk}{rk}

\newcounter{foo}
\makeatletter
\newcommand{\otherlabel}[2]{\protected@edef\@currentlabel{#2}\label{#1}}
\makeatother

\newcommand{\supth}[1]{\ensuremath{#1^{\mathrm{th}}}}

\EpigaVolumeYear{10}{2026} \EpigaArticleNr{13} \ReceivedOn{January 10, 2025}
\InFinalFormOn{November 22, 2025}
\AcceptedOn{February 11, 2026}

\title{Finite abelian groups acting on rationally connected threefolds~I: Groups of product type}
\titlemark{Finite abelian groups acting on rationally connected threefolds~I}

\author{Konstantin Loginov}
\address{Steklov Mathematical Institute of Russian Academy of Sciences, Moscow,
  Russia\\
  Centre of Pure Mathematics, MIPT, Moscow, Russia}
\email{loginov@mi-ras.ru}

\authormark{K.~Loginov}

\AbstractInEnglish{We initiate the study of finite abelian groups that faithfully act on $3$-dimensional rationally connected varieties.  We show that these groups can be naturally divided into three types: The groups of product type are finite abelian groups that are products of two groups that belong to the Cremona group of rank~$1$ and $2$, respectively; the groups of K3 type faithfully act on $G\mathbb{Q}$-Fano threefolds $X$ preserving a K3 surface $S\in|-K_X|$ with at worst du Val singularities; the third type consists of groups that act on $G\mathbb{Q}$-Fano threefolds with empty anti-canonical linear system.  The classification of groups of product type follows from a result of J.~Blanc.  For the groups of K3 type, we establish a boundedness result.  We also formulate a conjecture regarding the groups of the third type.}



\acknowledgement{This work was performed at the Steklov International Mathematical Center and supported by the Ministry of Science and Higher Education of the Russian Federation (agreement no. 075-15-2022-265), supported by the Simons Foundation. The work is supported by the state assignment of MIPT (project FSMG-2023-0013). The author is a Young Russian Mathematics award winner and would like to thank its sponsors and jury.} 

\begin{document}



\maketitle

\begin{prelims}

\DisplayAbstractInEnglish

\bigskip

\DisplayKeyWords

\medskip

\DisplayMSCclass

\end{prelims}


\newpage

\setcounter{tocdepth}{1}

\tableofcontents


\section{Introduction}
We work over the field of complex numbers $\mathbb{C}$. 

Finite abelian groups are among the most elementary objects studied in algebra. In algebraic geometry, rational varieties form a rather simple and well-behaved class of varieties. Consider the following question. 

\begin{question}
\label{main-question}
Which finite abelian groups can faithfully and biregularly act on rational (or rationally connected) varieties of a given dimension? 
\end{question}

Since we consider varieties endowed with the action of a finite group, there exists an equivariant resolution of singularities, so we can consider the class of finite abelian groups that act on rational (or rationally connected) smooth projective varieties. Also,  the notions of biregular and birational group actions would give us the same list of groups in answer to Question~\ref{main-question} since we can resolve the indeterminacy of the group action.

In dimension $1$, the answer to Question~\ref{main-question} follows from the classification of finite subgroups of $\Aut(\mathbb{P}^1)=\PGL(2, \mathbb{C})$ or, equivalently, from the classification of Platonic solids. The answer is: either a cyclic group or the Klein four-group. We formulate this elementary result for later use.

\begin{proposition}
\label{prop-abelian-subgroup-of-pgl2}
Let $G$ be a finite abelian subgroup of $\Aut(\mathbb{P}^1)=\PGL(2, \mathbb{C})$. Then $G$ is isomorphic to one of the following groups: 
\begin{enumerate}
\item
$\mathbb{Z}/n$\ for\ $n\geq1$,
\item
$(\mathbb{Z}/2)^2$. 
\end{enumerate}
\end{proposition}

Let us examine this proposition more closely. 
We see that there are two types of finite abelian groups that can faithfully act on $\mathbb{P}^1$: The groups of the first type are subgroups of a torus $\mathbb{C}^\times$ in $\PGL(2, \mathbb{C})$, and the group of the second type is rather special. First of all, there is only one such group, so there is a ``boundedness property.'' Then, the group of the second type induces a non-trivial automorphism of a divisor $D=\{0\}+\{\infty\}$ (in suitable coordinates) on $\mathbb{P}^1$ that is anti-canonical; that is, $D\sim -K_{\mathbb{P}^1}$,  {where $K_{\mathbb{P}^1}$ is the canonical class}. We can consider~$D$ as a ``Calabi--Yau variety'' embedded into $\mathbb{P}^1$ as its anti-canonical section. 
The goal of this work is to generalize this observation to higher dimensions. 

To put Question~\ref{main-question} into context, recall some definitions.  In what follows, we denote by $\Bir(X)$ the group of birational automorphisms of a variety $X$.  The \emph{Cremona group} $\Cr_n(\mathbb{C})=\Bir(\mathbb{P}^n)$ is the group of birational automorphisms of the $n$-dimensional projective space $\mathbb{P}^n$ over the field~$\mathbb{C}$. For $n=1$, this group is isomorphic to $\PGL(2, \mathbb{C})$. In contrast, for $n\geq 2$, the Cremona group becomes a more complicated object to work with. One way to understand its structure is by means of its finite subgroups.  {The study of finite subgroups of $\Cr_2(\mathbb{C})$ is classical and goes back to the 19th century, to the works of Kantor and Wiman \cite{Kan95,Wim96}.  However, the complete classification of finite subgroups of $\Cr_2(\mathbb{C})$ was obtained much later by Dolgachev and Iskovskikh; see \cite{DI09}. } Note that to answer Question~\ref{main-question} {in the case of rational varieties} of dimension $n$ is the same as to classify finite abelian subgroups of $\Cr_n(\mathbb{C})$ up to isomorphism.  The answer to Question~\ref{main-question} in dimension~$2$ is as follows.

\begin{thm}[\textit{cf.} {\cite{Bl07}}]
\label{thm-surfaces}
Let $X$ be a rational surface, and let $G\subset \Bir(X)=\Cr_2(\mathbb{C})$ be a
finite abelian group. Then~$G$ is isomorphic to one of the following groups:
\begin{enumerate}
\item
\label{cremona-plane-abelian-1}
$\mathbb{Z}/n\times \mathbb{Z}/m$\ for\ $n\geq
1, m\geq 1$,

\item
\label{cremona-plane-abelian-2}
$\mathbb{Z}/2n\times (\mathbb{Z}/2)^2$\ for\ $n\geq 1$,

\item
\label{cremona-plane-abelian-3}
$(\mathbb{Z}/4)^2\times \mathbb{Z}/2$,

\item
\label{cremona-plane-abelian-4}
$(\mathbb{Z}/3)^3$,

\item
\label{cremona-plane-abelian-5}
$(\mathbb{Z}/2)^4$.
\end{enumerate}
\end{thm}

Again, we can divide these groups into two classes. 
The groups~\eqref{cremona-plane-abelian-1},~\eqref{cremona-plane-abelian-2}, and~\eqref{cremona-plane-abelian-5} are products of groups that act on $\mathbb{P}^1$. In particular, they are realized as subgroups in
\[
\Cr_1(\mathbb{C})\times\Cr_1(\mathbb{C})\subset
\Cr_2(\mathbb{C}).
\] 
The groups~\eqref{cremona-plane-abelian-3} and~\eqref{cremona-plane-abelian-4} can be realized as subgroups of the automorphism group of a del Pezzo surface of degree~$2$ associated with a Fermat quartic curve, and of a Fermat cubic surface, respectively. In fact, the groups~\eqref{cremona-plane-abelian-3} and~\eqref{cremona-plane-abelian-4} are cyclic extensions of abelian groups acting on elliptic curves that are $G$-invariant anti-canonical sections of the corresponding del Pezzo surfaces. These curves admit complex multiplication (\textit{cf.}~Proposition~\ref{prop-elliptic-curve}), so we can consider them as anti-canonically embedded $1$-dimensional Calabi--Yau varieties with ``extra'' symmetries (see Section~\ref{sect-log-CY} for our definition of Calabi--Yau varieties). In more detail, this is explained in Section~\ref{sec-warm-up}, where Theorem~\ref{thm-surfaces} is re-proved. Note that there are only two groups of type~\eqref{cremona-plane-abelian-3}
and~\eqref{cremona-plane-abelian-4}, so again there is a ``boundedness property.'' It is also interesting to note that groups of type~\eqref{cremona-plane-abelian-3} and~\eqref{cremona-plane-abelian-4} cannot faithfully act on a conic bundle; see Proposition~\ref{prop-invariant-pencil}.
 
We recall some results related to Question~\ref{main-question}.  By a \emph{$p$-group} we mean a finite group of order $p^k$, where~$p$ is a prime number.  The following result can be deduced from Theorem~\ref{thm-surfaces}; however, it was obtained independently.

\begin{thm}[\textit{cf.} {\cite{Beau07}}]
\label{thm-p-surfaces}
Let $X$ be a rational surface, and let $G\subset \Bir(X)=\Cr_2(\mathbb{C})$ be a
finite abelian $p$-group. Then $G$ can be generated by at most $r$
elements, where
\begin{itemize}
\item
$r=4$ when $p=2$,
\item
$r=3$ when $p=3$,
\item
$r=2$ when $p\geq 5$.
\end{itemize}
Moreover, these bounds are sharp.
\end{thm}

In dimension $3$, the geometry becomes much more complicated, and the complete classification of finite subgroups of $\Cr_3(\mathbb{C})$ seems to be out of reach.
However, there exist classification results for some classes of finite groups; see \textit{e.g.}~\cite{Pr09} in the case of simple groups. There are also numerous boundedness results for finite subgroups of the Cremona group; see \cite{PS16} and references therein. Starting from dimension $3$, from the point of view of the minimal model program, it is more natural to consider the groups of birational automorphisms of rationally connected varieties of a given dimension. In a series of works, a generalization of Theorem~\ref{thm-p-surfaces} to dimension $3$ was obtained using explicit methods of the minimal model program. 

\begin{thm}[\textit{cf.} {\cite{Pr11,Pr14,PS17,Kuz20,Xu20,Lo22}}]
\label{thm-threefolds}
Let $X$ be a rationally connected variety of dimension~$3$, and let
$G\subset \Bir(X)$ be a finite $p$-group. Then $G$ can be
generated by at most $r$ elements, where
\begin{itemize}
\item
$r=6$ when $p=2$,
\item
$r=4$ when $p=3$,
\item
$r=3$ when $p\geq 5$.
\end{itemize}
Moreover, these bounds are attained in $\Cr_3(\mathbb{C})$.
\end{thm}

Recently, a generalization to higher dimensions was obtained using a different technique.

\begin{thm}[\textit{cf.} {\cite[Corollary 11]{KZh24}}]
Let $X$ be a rationally connected variety of dimension $n$, and let $G\subset \Bir(X)$ be a finite abelian $p$-group. Then $G$ can be generated by $r$ elements, where 
\[
r \leq \frac{pn}{p-1}\leq 2n.
\]
\end{thm}

Our goal is to generalize Theorem~\ref{thm-surfaces} to dimension $3$, that is, to obtain the classification of finite abelian groups that can faithfully act on $3$-dimensional rationally connected varieties. We introduce the following definition.  A finite abelian group $G$ is called \emph{a group of product type} if $G\simeq G_1\times G_2$, where $G_i\subset \Cr_i(\mathbb{C})$.  In particular, $G$ is isomorphic to a subgroup in
\[
\Cr_1(\mathbb{C})\times\Cr_2(\mathbb{C})\subset
\Cr_3(\mathbb{C}).
\] 
Hence, the list of groups of product type can be easily deduced from Proposition~\ref{prop-abelian-subgroup-of-pgl2} and Theorem~\ref{thm-surfaces}; see Table~\ref{table-1} in Section~\ref{sec-product-type}. 

In this paper, by a K3 surface we mean a normal projective surface $S$ with at worst canonical (that is, du Val) singularities such that $\mathrm{H}^1(S, \OOO_S)=0$ and $K_S\sim 0$.  We say that a finite abelian group~$G$ is of \emph{K3 type} if~$G$ faithfully acts on a $G\mathbb{Q}$-Fano threefold $X$ preserving a K3 surface $S\in|-K_X|$ with at worst du Val singularities. In particular, we have that $G$ is an extension of a finite abelian group $H$ that faithfully acts on a K3 surface by a cyclic group (cf. Corollary~\ref{cor-G-preserves-K3-type}):
\[
0\lra \mathbb{Z}/m\lra G \lra H\lra 0.
\]
Note that if a finite group faithfully acts on a singular K3 surface, then this action lifts to its minimal resolution.  Our main result is the following theorem (see Section~\ref{subsec-mfs} for the definitions of a $G\mathbb{Q}$-Fano variety and a $G\mathbb{Q}$-Mori fiber space).

\begin{thm}
\label{thm-first-main-thm}
Let $X$ be a rationally connected variety of dimension $3$, and let $G\subset \Bir(X)$ be a finite abelian group. Then 
\begin{enumerate}
\item
either $G$ is of product type, 
\item
or $G$ is of K3 type,
\item
or $G$ faithfully acts on a $G\mathbb{Q}$-Fano threefold $X'$ with $|-K_{X'}|=\emptyset$, and, moreover, any $G\mathbb{Q}$-Mori fiber space with rationally connected total space and a faithful action of $G$ is a $G\mathbb{Q}$-Fano threefold with empty anti-canonical system.
\end{enumerate}
\end{thm}

These three types are not mutually exclusive.  Actually, we expect
that there is no need for the groups of the third type in
Theorem~\ref{thm-first-main-thm}.  We formulate this as a conjecture. 

\begin{conjecture}
In the notation of Theorem~\ref{thm-first-main-thm}, any group of the third type is either of product type or of K3 type.
\end{conjecture}

On the other hand, there exist examples of groups of K3 type that are not of product type; see Example~\ref{exam-k3-not-pt}.  It follows that such groups do not admit a faithful action on a Mori fiber space whose total space is rationally connected and whose base is positive-dimensional; see Corollary~\ref{cor-exclude-positive-dim-base}.

Our next result is a boundedness property for groups of K3 type.  By the Weil index of a Fano variety $X$ with at worst canonical singularities, we mean the number $C(X)=\max\{ k\ |\, -K_X \sim k H\}$, where $H$ is a Weil divisor on $X$.

\begin{proposition}
\label{thm-second-main-thm}
Assume that $G$ is a group of K3 type. Then in the corresponding exact sequence\color{black}
\[
0\lra \mathbb{Z}/m\lra G \lra H\lra 0, 
\] 
where the group $H$ faithfully acts on a K3 surface, we have $m\leq C$, where $C$ is  the maximal Weil index of a  Fano threefold with at worst canonical singularities.   
\end{proposition}

The existence of the constant $C$ follows from the boundedness for Fano threefolds with canonical singularities; see \cite{B21}.  {Recently, an effective estimate $C=66$ was obtained in \cite{JL25}.}  On the other hand, a complete classification of finite groups that can faithfully act on a smooth K3 surface is given in \cite{BH23}. Since there exist only finitely many isomorphism classes of such groups, we obtain the following.

\begin{corollary}

There are only finitely many isomorphism classes of groups of K3 type.\color{black}
\end{corollary}

We will work with groups of K3 type in a forthcoming paper.  We also establish the following result on the extensions of finite abelian groups that can act on rational curves and rational surfaces.

\begin{proposition}[{Proposition~\ref{prop-abstact-extension}}]
\label{prop-abstact-extension-intro}
Let $H\subset \Cr_1(\mathbb{C})$ and $K\subset \Cr_2(\mathbb{C})$ be finite abelian groups. Then an abelian extension $G$ of\, $H$ by $K$ $($or $K$ by $H)$ is of product type.
\end{proposition}

Note that we do not claim that such an extension is necessarily split. We prove that under the assumptions of the proposition, $G$ is isomorphic to a product of some finite abelian groups $H'$ and $K'$ that are subgroups of $\Cr_1(\mathbb{C})$ and of $\Cr_2(\mathbb{C})$, respectively, and $H'$, $K'$ may be different from $H$, $K$.  This purely algebraic result on finite abelian groups is very useful for us. From this result it follows that if a finite abelian group $G$ faithfully acts on a rationally connected $3$-dimensional Mori fiber space with a positive-dimensional base, then such a group can faithfully act on a Fano threefold that is a product of a rational curve and a rational surface; see Corollary~\ref{cor-exclude-positive-dim-base}. This allows us to reduce our classification problem to the case of $G\mathbb{Q}$-Fano threefolds.  An analogous result holds in dimension $2$; see Proposition~\ref{prop-abstact-extension-cr2}.  This leads to the following natural question.

\begin{question}
\label{second-question}
Is there a finite group $G$ such that $G$ faithfully acts on a Mori fiber space $f\colon X\to Z$ with rationally connected total space~$X$, $\dim X=n$, and $\dim Z>0$, but $G$ does not admit a faithful action on a Fano variety with terminal singularities of dimension $n$?
\end{question}

Although the classification of finite subgroups of $\Cr_2(\mathbb{C})$ is obtained in \cite{DI09}, the answer to Question~\ref{second-question} is not known even in the case $n=2$; see \cite[Section 5.7]{DI09}.  By Proposition~\ref{prop-abstact-extension-intro}, the answer to Question~\ref{second-question} is negative for finite abelian groups in the case $n=3$.

We expect that Theorem~\ref{thm-first-main-thm} could be generalized to arbitrary dimension in the following way: Finite abelian groups that can faithfully act on a rationally connected variety~$X$ of dimension~$n$ either are products of groups that faithfully act on rationally connected varieties $X_i$ of dimension $n_i < n$ with $\sum n_i=n$, or they are cyclic extensions of groups faithfully acting on Calabi--Yau varieties of dimension $n-1$ anti-canonically embedded into some birational modification of $X$, or they act on a Fano variety with empty anti-canonical linear system.

\subsubsection*{Our approach} We present a sketch of a proof of the main results and explain the structure of the paper.  The idea of the proof of Theorem~\ref{thm-first-main-thm} is as follows.  Let $X$ be a $3$-dimensional projective rationally connected variety, and let~$G$ be a finite abelian group such that $G\subset \Bir(X)$. After passing to a $G$-equivariant resolution, we may assume that~$X$ is smooth and that $G\subset \Aut(X)$. We run a $G$-equivariant minimal model program that terminates on a $G\mathbb{Q}$-Mori fiber space $\pi\colon X'\to Z$.

In Section~\ref{sec-product-type}, we study the extension properties of finite abelian groups. In particular, we prove Proposition~\ref{prop-abstact-extension-intro} (= Proposition~\ref{prop-abstact-extension}). This excludes the case $\dim Z>0$ and allows us to assume that $X$ is a $3$-dimensional $G\mathbb{Q}$-Fano variety; see Corollary~\ref{cor-exclude-positive-dim-base}.

In Sections~\ref{sec-group-action-on-curves} and~\ref{sec-group-action-on-surfaces}, we do some preparations, considering the action of finite abelian groups on some classes of curves and surfaces. After that, it is easy to apply our methods to finite abelian subgroups of $\Cr_2(\mathbb{C})$ and re-prove Theorem~\ref{thm-surfaces}. This is done in Section~\ref{sec-warm-up}. After that, we consider the action of finite abelian groups on a germ of terminal singularity of dimension $3$. In Section~\ref{subsec-terminal-sing}, we prove that if there exists an (either smooth or singular) $G$-fixed point on $X$, then $G$ is of product type. We also prove some more precise results that will be used later on.

The main idea is to consider the action of $G$ on the anti-canonical linear system $|-K_X|$ of a $3$-dimensional $G\mathbb{Q}$-Fano variety.  If it is empty, we are in the last case of Theorem~\ref{thm-first-main-thm}. So we may assume that $|-K_X|$ is non-empty. Since the abelian group $G$ acts on $\mathrm{H}^0(X, \OOO(-K_X))$, there exists a $G$-invariant element $S\in |-K_X|$. We analyze singularities of the pair $(X, S)$; see Section~\ref{subsec-pairs-sing} for the relevant definitions. In Section~\ref{sec-sing-of-pair}, we show that either the pair $(X, S)$ is lc, or $G$ is of product type. On the other hand, we show that if $(X, S)$ is plt, then $G$ is of K3 type, and we are in the second case of Theorem~\ref{thm-first-main-thm}.

So we assume that $(X, S)$ is lc but not plt.  We prove that $G$ is of product type in this case, which is the most technical part of this paper.  We consider three possibilities.  If $S$ is irreducible, in Section~\ref{sec-sing-of-pair} we carefully apply the theory of elliptic surface singularities as well as threefold terminal singularities to show that $G$ is of product type.
In Section~\ref{sec-dual-complex}, we consider the case where $(X, S)$ is dlt and $S$ is reducible. Here we use the theory of the dual complex of Calabi--Yau pairs developed in \cite{KX16}. We prove that the dual complex of $(X, S)$ is homeomorphic either to a line segment or to a sphere $\mathbb{S}^2$; see Corollary~\ref{not-a-circle}. Using this result, we study combinatorics of the group action on the dual complex and conclude that~$G$ is of product type in this case. 
In Section~\ref{sec-reducible-lc}, we consider the case where $(X, S)$ is lc but not dlt, and $S$ is reducible. Here we also need the theory of the dual complex, as well as a boundedness result for the number of components of $S$ passing through a terminal point; see Proposition~\ref{prop-lc-bounded-number-of-components}. Again, this allows us to show that $G$ is of product type in this case. This proves Theorem~\ref{thm-first-main-thm}.

Finally, in Section~\ref{sect-action-on-K3}, we establish a boundedness result for groups of K3 type. We do this in Proposition~\ref{prop-K3-boundedness}, thus proving Proposition~\ref{thm-second-main-thm}. After that, we present some examples.

\subsection*{Acknowledgments}
The author thanks Ivan Cheltsov, Paul Hacking, J\'anos Koll\'ar, Joaqu\'in Moraga, Antoine Pinardin, Dmitry Pirozhkov, Yuri Prokhorov, Constantin Shramov, Andrey Trepalin, and Chenyang Xu for helpful conversations, and Alexander Kuznetsov for reading the draft of the paper and making many useful remarks.
The author is grateful to the anonymous referees for numerous helpful suggestions on the exposition. 

\section{Preliminaries}
 All the varieties are projective and defined over the field of complex numbers $\mathbb{C}$, unless stated otherwise. We use the language of the minimal model program (the MMP for short); see \textit{e.g.}~\cite{KM98}.

\subsection{Contractions} By a \emph{contraction} we mean a projective
morphism $f\colon X \to Y$ of normal varieties such that~$f_*\OOO_X =
\OOO_Y$. In particular, $f$ is surjective and has connected fibers. A
\emph{fibration} is defined as a contraction $f\colon X\to Y$ such
that $\dim Y<\dim X$.

Let $G$ be a finite group. By a \emph{$G$-variety} we mean a variety
$X$ together with an action of a group $G$. If $f$ is a
$G$-equivariant contraction (resp.\ a $G$-equivariant fibration) of
$G$-varieties, we call it a \emph{$G$-contraction} (resp.\ a
\emph{$G$-fibration}).

\subsection{Pairs and singularities} 
\label{subsec-pairs-sing}
A \emph{pair} $(X, B)$ consists of a normal variety $X$ and a boundary Weil $\mathbb{Q}$-divisor $B$ with coefficients in $[0, 1]$ such that $K_X + B$ is $\mathbb{Q}$-Cartier.  Let $\phi\colon W \to X$ be a log resolution of $(X,B)$, and let
\[
K_W +B_W = \phi^*(K_X +B)
\]
be the log pullback of $(X, B)$.  The \emph{log discrepancy} of a prime divisor $D$ on $W$ with respect to $(X, B)$ is $1 - \mathrm{coeff}_D B_W$, and it is denoted by $a(D, X, B)$. We say $(X, B)$ is lc (resp.\ klt, $\epsilon$-lc) if $a(D, X, B)\geq 0$ (resp.\ $> 0$, $\geq \epsilon$) for every $D$. Note that if $(X, B)$ is $\epsilon$-lc, then automatically $\epsilon \leq 1$ because $a(D, X, B) = 1$ for almost all $D$. We say that the pair $(X, B)$ is \emph{strictly lc} if it is lc but not klt.  We say that the pair $(X,B)$ is \emph{plt} if $a(D, X, B)>0$ holds for any $\phi$-exceptional divisor $D$ and for any log resolution $\phi$.  We say that the pair $(X,B)$ is \emph{dlt} if $a(D, X, B)>0$ holds for any $\phi$-exceptional divisor $D$ and for some log resolution~$\phi$.

An \emph{lc-place} of $(X, B)$ is a prime divisor $D$ over $X$, that
is, on a birational model of $X$, such that $a(D, X, B) = 0$. An
\emph{lc-center} is the image on $X$ of an lc-place.

\subsection{Log Fano and log Calabi--Yau pairs}
\label{sect-log-CY}
Let $(X, B)$ be an lc pair. We say $(X,B)$ is a \emph{log Fano pair} if $-(K_X+B)$ is ample.  If $B=0$, then $X$ is called a \emph{Fano variety}.  We say that an lc pair $(X,B)$ is a \emph{log Calabi--Yau pair} if $K_X + B \sim_{\mathbb{Q}} 0$. If $B=0$, then $X$ is called a \emph{Calabi--Yau variety}.

\subsection{Mori fiber space}
\label{subsec-mfs}
Let $G$ be a finite group.  Recall that a $G$-variety $X$ is called $G\mathbb{Q}$-factorial if every $G$-invariant Weil divisor on~$X$ is $\mathbb{Q}$-Cartier.  A $G\mathbb{Q}$-\emph{Mori fiber space} is a $G\mathbb{Q}$-factorial variety $X$ (we assume that the action of $G$ on $X$ is faithful) with at worst terminal singularities, together with a $G$-contraction $f\colon X\to Z$ to a variety $Z$ such that $\rho^G(X/Z)=1$ and $-K_X$ is ample over $Z$. If $X$ is $G$-factorial (\textit{e.g.}~if $X$ is smooth), we call it a $G$-\emph{Mori fiber space}.  If $Z$ is a point, we say that $X$ is a $G\mathbb{Q}$-Fano variety.  If $\dim X=2$ and $\dim Z=1$, then we say that $f\colon X\to Z$ is a \emph{$G$-conic bundle}.  If $\dim X=2$ and $Z$ is a point, then we say that $f\colon X\to Z$ is a \emph{$G$-del Pezzo surface}.

\subsection{Dual complex}
\label{subsec-dual-complex}
Let $D=\sum D_i$ be a Cartier divisor on a smooth projective variety~$X$. Recall that $D$ has \emph{simple normal crossings} (snc for short) if all the components $D_i$ of $D$ are smooth and any point in $D$ has an open neighborhood in the analytic topology that is analytically equivalent to the union of coordinate hyperplanes.

\emph{The dual complex}, denoted by $\mathcal{D}(D)$, of a simple normal crossing divisor $D=\sum_{i=1}^{r} D_i$ on a smooth variety~$X$ is a CW-complex constructed as follows.  The simplices $v_Z$ of $\mathcal{D}(D)$ are in bijection with the  irreducible components $Z$ (called \emph{strata of} $D$) of the intersection $\bigcap_{i\in I} D_i$ for any non-empty subset $I\subset \{ 1, \ldots, r\}$, and the vertices of $v_Z$ correspond to the components $D_i$ with $i\in I$.  In particular, the dimension of $v_Z$ is equal to $\#I-1$.  The gluing maps are constructed as follows.  For any non-empty subset $I\subset \{ 1, \ldots, r\}$, let $Z\subset \bigcap_{i\in I} D_i$ be an irreducible component, and for any $j\in I$, let $W$ be the unique component of $\bigcap_{i\in I\setminus\{j\}} D_i$ containing~$Z$. Then the gluing map is the inclusion of $v_W$ into $v_Z$ as a face of $v_Z$ that does not contain the vertex $v_i$ corresponding to $D_i$. Note that the dimension of $\mathcal{D}(D)$ does not exceed $\dim X-1$. If $\mathcal{D}(D)$ is empty, that is, $D=0$, we say that $\dim \mathcal{D}(D)=-1$.

In what follows, for a divisor $D$, we denote by $D^{=1}$ the sum of the components of $D$ with coefficient $1$. For an lc log CY pair $(X, D)$, we define $\mathcal{D}(X, D)$ as $\mathcal{D}(D_Y^{=1})$, where $f\colon (Y, D_Y)\to (X, D)$ is a log resolution of $(X, D)$, so the formula
\[
K_{Y} + D_Y= f^*(K_X + D)
\]
is satisfied. It is known that the PL-homeomorphism class of $\mathcal{D}(D_Y^{=1})$ does not depend on the choice of a log resolution; see \cite[Proposition 11]{dFKX17}.  For more results on the topology of dual complexes of Calabi--Yau pairs, see \cite{KX16}.

\begin{remark}
Note that if a finite group $G\subset \Aut(X)$ preserves a divisor $D$ on a projective variety $X$ such that the pair $(X, D)$ is lc, then $G$ acts on the topological space $\mathcal{D}(X, D)$ by PL-homeomorphisms.
\end{remark}

\subsection{Group action on spheres}
We recall some facts about the action of finite groups on spheres of lower dimensions. 

\begin{lem}
\label{lem-action-on-spheres-1}
Let $G$ be a finite abelian group that faithfully acts on a circle $\mathbb{S}^1$ by PL-homeomorphisms. Then $G=\mathbb{Z}/n$ for some $n\geq 1$ or $G=\mathbb{Z}/2\times\mathbb{Z}/2$. 
\end{lem}

\begin{proof}
Observe that $G$ can be realized as a subgroup of symmetries of a regular $n$-gon in a plane. 
 
Hence, $G$ is a subgroup of the dihedral group $D_n=\mathbb{Z}/n \rtimes \mathbb{Z}/2$. Thus $G$ fits into an exact sequence $0 \to H \to G \to K \to 0$, where $K$ is either trivial or $\mathbb{Z}/2$ and $H$ is cyclic contained in $\mathbb{Z}/n$. If $K$ is trivial, $G$ is cyclic. If $K = \mathbb{Z}/2$, then $H$ is fixed by the action of $\mathbb{Z}/2$ on $\mathbb{Z}/n$, but the action is $x \mapsto -x$, so $H$ is either trivial or $\mathbb{Z}/2$. This proves that $G$ is either cyclic or isomorphic to $\mathbb{Z}/2 \times \mathbb{Z}/2$.
\end{proof}

\begin{lem}[\textit{cf.} {\cite[Section 2]{Zim12}}]
\label{lem-action-on-spheres}
Let $G$ be a finite group that faithfully acts on a sphere~$\mathbb{S}^2$ by PL-ho\-meo\-mor\-phisms. Then $G$ is conjugate to a subgroup of the orthogonal group~$O(3)$. In particular, if\, $G$ is abelian, then either $G\subset \mathbb{Z}/n\times\mathbb{Z}/2$ for some $n\geq 1$, or $G=(\mathbb{Z}/2)^3$. Consequently, either $G$ has a fixed point on $\mathbb{S}^2$, or $G$ has an orbit of length~$2$ on~$\mathbb{S}^2$.
\end{lem}

\subsection{Group actions on algebraic varieties.}
We collect some generalities on the actions of finite groups on algebraic varieties. 

\begin{defin}
Let $G$ be a finite group. By $\mathfrak{r}(G)$, we denote the \emph{rank} of $G$, that is, the minimal number of generators of $G$.
\end{defin}

\begin{remark}
\label{rem-rank-subgroup}
Note that if $G$ is a finite abelian group and $H\subset G$ is a subgroup, then $\mathfrak{r}(H)\leq \mathfrak{r}(G)$. 
\end{remark}

\begin{lem}[\textit{cf.} {\cite[Lemma 4]{Po14}}]
\label{lem-faithful-action}
Let $X$ be an algebraic variety, and let $G$ be a finite group such that $G\subset \Aut(X)$. Assume $P\in X$ is a fixed point of\, $G$.  Then the induced action of\, $G$ on the tangent space $T_P X$ is faithful.
\end{lem}

\begin{corollary}
\label{cor-faithful-action}
If a finite abelian group $G\subset \Aut(X)$ has a fixed point $P\in X$, where $X$ is an algebraic variety of dimension $n$ and~$P$ is a smooth point, then $\mathfrak{r}(G)\leq n$.
\end{corollary}

In what follows, we will need the following results. 

\begin{lem}[\textit{cf.} {\cite[Lemma 2.6]{Pr11}}]
\label{lem-fixed-curve-surface}
Let $X$ be a $3$-dimensional algebraic variety $X$ with isolated singularities, and let $G$ be a finite abelian group such that $G\subset \Aut(X)$.
\begin{enumerate}
\item\label{lfcs-1} 
If there is a curve $C\subset X$ of\, $G$-fixed points, then $\mathfrak{r}(G) \leq 2$.
\item\label{lfcs-2} 
If there is a $($possibly, reducible$)$ divisor $S\subset X$ of\, $G$-fixed points, then $\mathfrak{r}(G)\leq 1$.  If moreover $S$ is singular along a curve, then $G$ is trivial.
\end{enumerate}
\end{lem}

\begin{proof}
We prove the first claim. Let $P$ be a general point in $C$. Since $X$ has isolated singularities, $P$ is smooth on $X$. In particular, $\dim T_{P}X=3$. By Lemma~\ref{lem-faithful-action}, the induced action of $G$ on $T_{P}X$ is faithful. On the other hand, the action of $G$ on the $1$-dimensional subspace $T_P C \subset T_{P}X$ is trivial, since by assumption the action of $G$ is trivial on $C$. Hence $G$ faithfully acts on a $G$-invariant complement $T'$ to $T_P C$ in $T_P X$. So we have $\dim T'=2$. Since $G$ is a finite abelian group, it follows that $\mathfrak{r}(G)\leq 2$ in this case.  

The first part of the second claim is proven analogously. We prove the second part.  By assumption, the action of $G$ on $S$ it trivial. It follows that the induced action of $G$ on $T_P S$ is trivial for any point $P\in S$. For a general point $P\in \Sing(S)$, we have $T_P S=T_P X$ since $X$ is smooth at $P$ and $P$ is singular on $S$. By Lemma~\ref{lem-faithful-action}, we conclude that $G$ is trivial.
\end{proof}

\begin{corollary}
\label{cor-G-preserves-K3-type}
If a finite abelian group $G$ faithfully acts on a $3$-dimensional algebraic variety~$X$ with isolated singularities preserving a surface $S$, then there is an exact sequence
\[
0\lra \mathbb{Z}/m \lra G \lra H \lra 0
\]
for $m\geq 1$, where $H\subset \Aut(S)$. 
\end{corollary}

\begin{proof}
Denote the image of a natural map $\phi\colon G\to \Aut(S)$ by $H$. By Lemma~\ref{lem-fixed-curve-surface}\eqref{lfcs-2}, the kernel of $\phi$ is cyclic. The claim follows. 
\end{proof}

\begin{lem}
\label{lem-diagonal-in-product}
Let $X$ be an algebraic variety, and let $G$ be a finite abelian group such that $G\subset \Aut(X)$.  Let $S=\sum_{i=1}^N S_i$ be a divisor on $X$ such that the induced action of\, $G$ on $S$ is faithful.  Denote by $G_S$ the image of\, $G$ in the symmetric group on the set $\{S_i\}$.  Assume that the action of\, $G_S$ on the set of components $\{S_i\}$ is transitive.  Then in the exact sequence
\begin{equation}
\label{eq-action-on-a-surface}
0\lra H\lra G\lra G_S\lra 0,
\end{equation}
the group $H$ is isomorphic to a subgroup of $\Aut(S_i)$.
\end{lem}

\begin{proof} 
Fix a component $S_1$ of $S$. 
Since the action of $G_S$ on the set $\{S_i\}$ is transitive, for any component $S_j$, there exists an element $g'_j\in G_S$ such that $g'_j(S_1)=S_j$. Pick preimages $g_j\in G$ of the elements $g'_j\in G_S$. 
For any $j$, we identify $S_1$ and $S_j$ via the isomorphism $S_1\simeq g_j(S_1)=S_j$. 
Since by assumption the action of $G$ on~$S$ is faithful, $G$ can be included in the exact sequence \eqref{eq-action-on-a-surface} with 
$
H\subset \prod_{i=1}^{N} \Aut(S_i).
$
Consider an arbitrary element 
\[
h\in H\subset \prod_{i=1}^{N} \Aut(S_i)\simeq \prod_{i=1}^N \Aut(S_1), 
\]
where we used the isomorphisms chosen above. 
Let $\pr_j\colon \prod_{i=1}^N \Aut(S_1)\to \Aut(S_1)$ be the projection to $\supth{j}$ factor. 
Since $G$ is abelian, for any $j$, we have $g_jhg_j^{-1}=h$. It follows that $\pr_i(h)=\pr_j(h)$ for any $i, j$. 
Let 
\[
\Delta\colon \Aut(S_1)\lra \prod_{i=1}^N \Aut(S_1)
\] 
be the diagonal embedding: $\Delta(g)=(g,\ldots, g)$. Then $H\simeq \Delta(\pr_1(H))$. Thus, $H\simeq \pr_1(H)\subset  \Aut(S_1)$, and the claim follows. 
\end{proof}

\section{Groups of product type}
\label{sec-product-type}
In this section, we fix some notation related to finite abelian groups and introduce the definition of abelian groups of product type. 
We also study extension properties of finite abelian groups. 
We start with the following lemma, which follows from the classification of finite abelian groups; see also \cite[Theorem 5]{Ka69} for a more general statement. 

\begin{lem}
\label{lem-decreasing-rank}
Let $G$ be a finite abelian group. Let $g\in G$ be an element, and let $\langle g \rangle\subset G$ be a subgroup generated by~$g$. Assume that $\langle g \rangle$ is a maximal cyclic subgroup in $G$, which means that if there is a cyclic subgroup $H\subset G$ such that $\langle g\rangle\subset H$, then $H=\langle g\rangle$.  
Then $G\simeq \langle g \rangle\times K$ for some subgroup $K\subset G$. In particular, $\mathfrak{r}(G/\langle g \rangle)< \mathfrak{r}(G)$.
\end{lem}

\begin{corollary}
\label{cor-another-splitting}
Let $G$ be a finite abelian group such that $G$ fits in the exact sequence
\begin{equation}
\label{eq-one-more-exact-sequence}
0 \lra H \lra G \overset{\phi}\lra  K \lra 0.
\end{equation}
Assume that $H$ is cyclic. Then either the sequence \eqref{eq-one-more-exact-sequence} splits, or we have $G\simeq H^+\times K^-$, where $H^+$ is a cyclic group that contains $H$ as a subgroup and $K^-=K/\phi(H^+)$\color{black}. 
\end{corollary}

\begin{proof}
Assume that there is no cyclic subgroup $H^+$ of $G$ such that $H\subset H^+$ and $H\neq H^+$. Then by Lemma~\ref{lem-decreasing-rank}, we have $G\simeq H\times K$, and the exact sequence \eqref{eq-one-more-exact-sequence} splits. So we may assume that there is a cyclic subgroup $H^+$ of $G$ such that $H\subset H^+$ and $H\neq H^+$. We may assume that $H^+$ is maximal with this property. 
Put $K^-=K/\phi(H^+)$. 
Then we have an exact sequence
\[
0 \lra H^+ \lra G \lra K^- \lra 0.
\]
Since the cyclic subgroup $H^+$ is maximal in $G$, by Lemma~\ref{lem-decreasing-rank}, we have $G\simeq H^+\times K^-$. 
\end{proof}

Let $G$ be a finite abelian group. 
In what follows, by $G_p$ we denote the $p$-Sylow subgroup of $G$, where $p$ is a prime number, so we have
\[
G=\mathop{\Scale[1.2]{\prod}}_{p\geq 2} G_p.
\]
We occasionally say that $G_p$ is the \emph{$p$-part} of $G$.  

\begin{remark}
Assume that $\mathfrak{r}(G_p)\leq k$ for any prime number $p$. Then $\mathfrak{r}(G)\leq k$.
\end{remark}

For example, by Theorem~\ref{thm-threefolds}, if $G$ faithfully acts on a rationally connected threefold, then $\mathfrak{r}(G_2)\leq 6$, $\mathfrak{r}(G_3)\leq 4$, and $\mathfrak{r}(G_p)\leq 3$ for $p\geq 5$. 

For a finite abelian group $G$, we put 
\[
G_{\neq p}=\mathop{\Scale[1.2]{\prod}}_{q\neq p} G_q, \quad \text{so }  G = G_p\times G_{\neq p}. 
\]
For an abelian $p$-group $G_p$, we say that $G_p$ has type 
\[
\lambda=[\lambda_1,\ldots,\lambda_k] \quad \text{for }  \lambda_1\geq\cdots\geq \lambda_k
\]
if 
\[
G_p=\mathbb{Z}/p^{\lambda_1}\times\cdots\times\mathbb{Z}/p^{\lambda_k}.
\]
\begin{remark}
\label{rem-pass-to-p-part}
A sequence of finite abelian groups 
\begin{equation}
\label{extension-just-another-extension}
0\lra H\lra G\lra K \lra 0
\end{equation}
is exact if and only if for any prime $p$ the $p$-parts $H_p$, $G_p$, $K_p$ of the groups $H$, $G$, $K$, respectively, form an exact sequence
\begin{equation}
\label{extension-p-groups}
0\lra H_p\lra G_p\lra K_p \lra 0.
\end{equation}
We say that the exact sequence \eqref{extension-p-groups} is the \emph{$p$-part} of the exact sequence \eqref{extension-just-another-extension}.
\end{remark}

To any type $\lambda=[\lambda_1,\ldots,\lambda_k]$ corresponds the Young diagram with $\lambda_i$ squares in the $\supth{i}$ row. For two Young diagrams $\lambda=[\lambda_1,\ldots,\lambda_k]$ and $\mu=[\mu_1,\ldots \mu_k]$, one can define their product $\lambda\cdot \mu$ as a formal linear combination of Young diagrams with non-negative coefficients; see \textit{e.g.}~\cite[Section 2]{Fu00}. Then the \emph{Littlewood--Richardson coefficient} $c^\nu_{\lambda\mu}$ is the coefficient at the Young diagram $\nu=[\nu_1,\ldots \nu_k]$ in the product of Young diagrams~$\lambda\cdot \mu$. 

\begin{thm}[\textit{cf.} {\cite[Section 2]{Fu00}}]
\label{thm-fulton}
In the above notation, assume that $G_p$ has type $\lambda=[\lambda_1,\ldots,\lambda_k]$, $H_p$ has type $\mu=[\mu_1,\ldots,\mu_k]$, and $K_p$ has type $\nu=[\nu_1,\ldots,\nu_k]$. 
Then an extension of the form \eqref{extension-p-groups} exists if and only if for the Littlewood--Richardson coefficient, we have $c^\lambda_{\mu\nu}>0$. 
\end{thm}

\begin{remark}
\label{rem-LR-coef-is-symmetric}
Let $G$, $H$, and $K$ be finite abelian groups. Since $c^\lambda_{\mu\nu}=c^\lambda_{\nu\mu}$ (\textit{cf.} \cite[Corollary~I.5.1.3]{Fu06}), from Remark~\ref{rem-pass-to-p-part} and Theorem~\ref{thm-fulton}. it follows that an abelian extension of $K$ by $H$ isomorphic to $G$ exists if and only if there exists an abelian extension of $H$ by $K$ isomorphic to $G$.
\end{remark}

\begin{defin}
\label{cremona-product-list}
We say that a finite abelian group $G$ is \emph{a group of product type} if~$G\simeq G_1\times G_2$, where $G_i\subset \Cr_i(\mathbb{C})$.  In particular, $G$ is isomorphic to a subgroup in
\[
\Cr_1(\mathbb{C})\times\Cr_2(\mathbb{C})\subset
\Cr_3(\mathbb{C}).
\] 
\end{defin}

\begin{proposition}
\label{prop-two-conditions-for-abelian-groups}
For a finite abelian group $G$, the following conditions are equivalent: 
\begin{enumerate}
\item\label{ptcfag-1}
$G\subset A_1\times A_2$, where the $A_i$ are arbitrary groups;
\item\label{ptcfag-2}
$G\simeq G_1\times G_2$, where $G_i\subset A_i$.
\end{enumerate}
Thus, a finite abelian group $G$ is of product type if and only if\, $G$ is isomorphic to a subgroup of\, $\Cr_1(\mathbb{C})\times\Cr_2(\mathbb{C})$.
\end{proposition}
\begin{proof}

The implication $(2)\Rightarrow(1)$ is clear. Assume that $G\subset A_1\times A_2$. 
Let $P_i=\pi_i(G)\subset A_i$. 
It is enough to argue separately for each $p$-primary component. Thus,
\[
G_p\subset P_{1,p}\times P_{2,p}.
\]

Write
\[
G_p\simeq \bigoplus_{j=1}^{r}\mathbb Z/p^{\lambda_j}\mathbb Z,
\qquad
P_{1,p}\times P_{2,p}
\simeq
\bigoplus_{j=1}^{s}\mathbb Z/p^{\gamma_j}\mathbb Z,
\]
where the exponents are arranged in non-increasing order. Clearly, one has $\lambda_j\leq \gamma_j$ for all $j$. 
Also, each cyclic factor $\mathbb Z/p^{\gamma_j}\mathbb Z$ belongs either to
$P_{1,p}$ or to $P_{2,p}$. Assign the corresponding factor
$\mathbb Z/p^{\lambda_j}\mathbb Z$ to $G_{1,p}$ or $G_{2,p}$ accordingly.
Since $\lambda_j\leq\gamma_j$, we have
$
\mathbb Z/p^{\lambda_j}\mathbb Z
\subset
\mathbb Z/p^{\gamma_j}\mathbb Z.
$
Hence
\[
G_p\simeq G_{1,p}\times G_{2,p},
\qquad
G_{i,p}\subset P_{i,p}\subset A_i.
\]

Taking the product over all primes $p$, define
\[
G_i=\prod_p G_{i,p}.
\]
Then
$G\simeq G_1\times G_2$, $G_i\subset A_i$.
This proves $(1)\Rightarrow(2)$.
\color{black}
\end{proof} 

Note that Proposition~\ref{prop-two-conditions-for-abelian-groups} does not hold for arbitrary (finite) groups. 

Using Proposition~\ref{prop-abelian-subgroup-of-pgl2} and Theorem~\ref{thm-surfaces}, one checks that the groups of product type are exactly the 
following groups:

\begin{center}
  \begin{table}[h!]
    {\renewcommand{\arraystretch}{1.1}
\begin{tabular}{ | m{1.3em} | m{5.0cm} | m{3.5cm} | } 
  \hline
   & $G$ &  \\ 
  \hline
  (1) & $\mathbb{Z}/k\times \mathbb{Z}/l\times
\mathbb{Z}/m$ & $k\geq 1,\ l\geq 1,\ m\geq 1$ \\ 
  \hline
  (2) & $\mathbb{Z}/2k\times(\mathbb{Z}/4)^2\times
\mathbb{Z}/2$ & $k\geq 1$ \\ 
    \hline
  (3) & $\mathbb{Z}/3k\times(\mathbb{Z}/3)^3$ & $k\geq 1$ \\ 
    \hline
  (4) & $\mathbb{Z}/2k\times \mathbb{Z}/2l\times
(\mathbb{Z}/2)^2$ & $k\geq 1,\ l\geq 1$ \\ 
    \hline
  (5) & $\mathbb{Z}/2k\times (\mathbb{Z}/2)^4$ & $k\geq 1$ \\ 
    \hline
  (6) & $(\mathbb{Z}/4)^2\times (\mathbb{Z}/2)^3$ & \\ 
    \hline
  (7) & $(\mathbb{Z}/2)^6$ &  \\ 
    \hline
\end{tabular}}
\caption{Groups of product type}
  \label{table-1}
\end{table}
\end{center}

Note that the class of groups of product type is closed under taking subgroups and quotient groups.  We will use the following observation repeatedly.

\begin{remark}
\label{rem-rank-bounded-by-3-pt}
If $\mathfrak{r}(G)\leq 3$, then $G$ is of product type.  
\end{remark}

Another simple observation is the following proposition, which can be proven by an easy case-by-case analysis. 

\begin{proposition}
\label{prop-abstact-extension-cr2}
Let $H, K\subset \Cr_1(\mathbb{C})$ be finite abelian groups. Let $G$ be an abelian group that is an extension of\, $H$ by $K$ $($or $K$ by $H)$. Then we have 
\[
G\simeq H'\times K'
\] 
for some subgroups $H', K'\subset \Cr_1(\mathbb{C})$.  
In particular, $G$ is isomorphic to a subgroup of 
\[
\Cr_1(\mathbb{C})\times \Cr_1(\mathbb{C})\subset \Cr_2(\mathbb{C}).
\] 
Moreover, $G$ has type~\eqref{cremona-plane-abelian-1},~\eqref{cremona-plane-abelian-2}, or \eqref{cremona-plane-abelian-5} as in Theorem~\ref{thm-surfaces}.
\end{proposition}

It is a kind of surprise that an analogous result holds in dimension $3$.

\begin{proposition}
\label{prop-abstact-extension}
Let $H\subset \Cr_1(\mathbb{C})$ and $K\subset \Cr_2(\mathbb{C})$ be finite abelian groups. Then an abelian extension $G$ of\, $K$ by $H$ $($or $H$ by $K)$ is of product type.
\end{proposition}

\begin{proof}
According to Remark~\ref{rem-LR-coef-is-symmetric}, it is enough to consider the case of an extension of $K$ by~$H$. 
We have an exact sequence
\begin{equation}
\label{exact-sequence-product-type-2}
0\lra H \lra G \lra K \lra 0, 
\end{equation}
where $G$ is a finite abelian group. 
According to Proposition~\ref{prop-abelian-subgroup-of-pgl2}, there are two possibilities for $H$: 
\begin{enumerate}
\item
$\mathbb{Z}/k$ for $k\geq 1$, 
\item
$(\mathbb{Z}/2)^2$. 
\end{enumerate}
According to Theorem~\ref{thm-surfaces}, there are the following possibilities for $K$:
\begin{enumerate}
\item\label{K1}
$\mathbb{Z}/n\times \mathbb{Z}/m$ for $n\geq
1, m\geq 1$,

\item\label{K2}
$\mathbb{Z}/2n\times (\mathbb{Z}/2)^2$ for $n\geq 1$,

\item\label{K3}
$(\mathbb{Z}/4)^2\times \mathbb{Z}/2$,

\item\label{K4}
$(\mathbb{Z}/3)^3$,

\item\label{K5}
$(\mathbb{Z}/2)^4$.
\end{enumerate}
If the exact sequence \eqref{exact-sequence-product-type-2} splits, the group $G$ is of product type. So we may assume that it does not split. 
Assume that $H$ is cyclic.
Then, according to Corollary~\ref{cor-another-splitting}, we have $G\simeq H^+\times K^-$, where $H^+$ is cyclic and $K^-$ is a quotient of $K$. Since finite abelian subgroups of $\Cr_2(\mathbb{C})$ are closed under taking quotient groups, we have that $G$ is of product type. 
So we may assume that $H=(\mathbb{Z}/2)^2$. 
We apply Theorem~\ref{thm-fulton} and consider five cases.

\emph{Case}~\eqref{K1}.~ Assume that $K$ has type~\eqref{K1} as in Theorem~\ref{thm-surfaces}.  Consider the $2$-part of the exact sequence \eqref{exact-sequence-product-type-2}.  Let $(k, l)$ be the type of $K_2$, where $k\geq l$. We may assume that $l\geq 1$; otherwise, $G=G_2\times G_{\neq 2}$ is of product type, because $\mathfrak{r}(G_{\neq 2})\leq 2$.  Then the product of the Young diagrams corresponding to $K_2$ and $H=H_2$ is as follows:
\[
[k,l]\cdot[1,1] = [k+1,l+1] + [k+1, l, 1] + [k, l+1, 1] + [k,l,1,1].
\]
Here $[k,l+1,1]$ occurs only when $k>l$.
Since $G=G_2\times G_{\neq 2}$ and $\mathfrak{r}(G_2)\leq 3$ for the first three summands, we have $\mathfrak{r}(G)\leq 3$. So
we need only consider the type $[k,l,1,1]$. 
Since $\mathfrak{r}(G_{\neq 2})\leq 2$, it follows that $G=G_2\times G_{\neq 2}$ is of product type.

\emph{Case}~\eqref{K2}.~ 
Assume that $K$ has type~\eqref{K2}.
Consider the product of the Young diagrams for the $2$-part of the exact sequence~\eqref{exact-sequence-product-type-2}:
\[
[k,1,1]\cdot[1,1] = [k+1, 2, 1] + [k+1, 1, 1, 1] + [k, 2, 2] + [k, 2, 1, 1] + [k, 1, 1, 1, 1].
\]
Here $[k,2,2]$ and $[k,2,1,1]$ occur only when $k>1$.
Since $\mathfrak{r}(G_{\neq 2})\leq 1$, it follows that $G=G_2\times G_{\neq 2}$ is of product type.

\emph{Case}~\eqref{K3}.~ 
Assume that $K$ has type~\eqref{K3}.
Consider the product of the Young diagrams for the $2$-part of the exact sequence~\eqref{exact-sequence-product-type-2}:
\[
[2, 2, 1]\cdot[1,1] = [3, 3, 1] + [3, 2, 2] + [3, 2, 1, 1] + [2, 2, 2, 1] + [2, 2, 1, 1, 1].
\]
It follows that $G=G_2$ is of product type. 

\emph{Case}~\eqref{K4}.~  
Assume that $K$ has type~\eqref{K4}.
Clearly, in this case we have $\mathfrak{r}(G)\leq 3$, so $G$ is of product type. 

\emph{Case}~\eqref{K5}.~  
Finally, assume that $K$ has type~\eqref{K5}.
Consider the product of the Young diagrams for the exact sequence~\eqref{exact-sequence-product-type-2}:
\[
[1,1,1,1]\cdot[1,1] = [2, 2, 1, 1] + [2, 1, 1, 1, 1] + [1, 1, 1, 1, 1, 1].
\]
It follows that $G=G_2$ is of product type.
\end{proof}

\begin{corollary}
\label{cor-inv-curve-or-surface}
Let $X$ be a $3$-dimensional algebraic variety $X$ with isolated singularities, and let $G$ be a finite abelian group such that $G\subset \Aut(X)$.
If there exists a $G$-invariant rational curve $C$ on $X$, then $G$ is of product type.
If there exists a $G$-invariant rational surface $S$ on $X$, then $G$ is of product type.
\end{corollary}

\begin{proof}
Assume that there exists a $G$-invariant rational curve $C$ on $X$. Then we have an exact sequence
\begin{equation}
0\lra G_N \lra G \lra G_C \lra 0, 
\end{equation}
where $G_C$ faithfully acts on $C$, so $G_C\subset \Cr_1(\mathbb{C})$, and $G_N$ faithfully acts on the normal bundle to $C$, and so $\mathfrak{r}(G_N)\leq 2$ (\textit{cf.}~Lemma~\ref{lem-fixed-curve-surface}). In particular, $G_N\subset \Cr_2(\mathbb{C})$. Now from Proposition~\ref{prop-abstact-extension} it follows that $G$ is of product type. 

Assume that there exists a $G$-invariant rational surface $S$ on $X$. Then in the notation of Corollary~\ref{cor-G-preserves-K3-type}, we have $\mathbb{Z}/m\subset \Cr_1(\mathbb{C})$ and $H\subset \Cr_2(\mathbb{C})$. By Proposition~\ref{prop-abstact-extension}, we conclude that $G$ is of product type. 
\end{proof}

The next useful lemma follows from the Lefschetz principle. 

\begin{lem}
\label{lem-lefschetz}
Let $\mathbb{K}$ be a field of characteristic $0$. Assume that a finite group $G$ is isomorphic to a subgroup of\, $\Cr_n(\mathbb{K})$ for some $n\geq 1$. Then $G$ is isomorphic to a subgroup of\, $\Cr_n(\mathbb{C})$. 
\end{lem}

\begin{corollary}
\label{cor-exclude-positive-dim-base}
Let $G$ be a finite abelian group. 
Let $f\colon X\to Z$ be a $G\mathbb{Q}$-Mori fiber space such that $\dim X=3$, $\dim Z>0$, and $X$ is rationally connected. Then $G$ is of product type. 
\end{corollary}

\begin{proof}
Note that $Z$ is rational. Assume that $\dim Z=2$ (the case $\dim Z=1$ is analogous). Then $\dim X_\eta=1$, where $X_\eta$ is the scheme-theoretic generic fiber. We have an exact sequence
\begin{equation}
\label{eq-mfs-exact-sequence}
0\lra G_f \lra G \lra G_Z \lra 0, 
\end{equation}
where $G_Z$ faithfully acts on $Z$, and so $G_Z\subset \Cr_2(\mathbb{C})$, and $G_f$ faithfully acts on $X_\eta$. Since $X_\eta$ is a smooth Fano curve over $\overline{\mathbb{C}(Z)}$, it is rational over $\overline{\mathbb{C}(Z)}$, so $G_f\subset \Cr_1(\overline{\mathbb{C}(Z)})$. By Lemma~\ref{lem-lefschetz}, it follows that $G_f$ is isomorphic to a subgroup in $\Cr_1({\mathbb{C}})$. 
Then the claim follows from Proposition~\ref{prop-abstact-extension}.
\end{proof}

We formulate the following result for later use.

\begin{lem}
\label{cor-r2-4-rp-2}
Let $H$ be a group of type~\eqref{cremona-plane-abelian-1},~\eqref{cremona-plane-abelian-2}, or~\eqref{cremona-plane-abelian-5} as in Theorem~\ref{thm-surfaces}.  Let $G$ be an abelian cyclic extension of\, $H$. Assume that $\mathfrak{r}(G_2)= 4$. Then $\mathfrak{r}(G_p)\leq 2$ for $p>2$. In particular, we have
\[
G=(\mathbb{Z}/2)^2\times\mathbb{Z}/2n\times\mathbb{Z}/2m \quad  \text{for some } n,m\geq 1.
\] 
\end{lem}

\begin{proof}
Since $\mathfrak{r}(G_2)=4$ and $G$ is a cyclic extension of $H$, it follows that for the $2$-part $H_2$ of~$H$, we have $\mathfrak{r}(H_2)\geq 3$. Consequently, $H$ has type~\eqref{cremona-plane-abelian-2} or~\eqref{cremona-plane-abelian-5} as in Theorem~\ref{thm-surfaces}. In particular, $\mathfrak{r}((H)_p) \leq 1$ for $p>2$, where $H_p$ is the $p$-part $H$. Thus $\mathfrak{r}(G_p)\leq 2$ for $p>2$. Applying Theorem~\ref{thm-fulton} to the $2$-part shows that $G_2$ has type $[a,b,1,1]$; together with the estimate above for odd $p$, this gives the stated decomposition.
\end{proof}

\section{Group action on curves}
\label{sec-group-action-on-curves}
In this section, we prove some results on the action of finite abelian groups on smooth projective curves.  The following proposition is a slightly more precise version of Proposition~\ref{prop-abelian-subgroup-of-pgl2}.

\begin{proposition}
\label{prop-rational-curve}
Let $G$ be a finite abelian subgroup of $\Aut(\mathbb{P}^1)\simeq\PGL(2, \mathbb{C})$.

If\, $G=\mathbb{Z}/n$ for $n\geq 2$, then~$G$ has two fixed points on $\mathbb{P}^1$, and any other point on $\mathbb{P}^1$ has a $G$-orbit of cardinality~$n$. 

If\, $G=(\mathbb{Z}/2)^2$, then there exist three $G$-orbits on~$\mathbb{P}^1$ of cardinality $2$, and any other point on $\mathbb{P}^1$ has a $G$-orbit of cardinality $4$. 
\end{proposition}

\begin{proof}
A cyclic subgroup $\mathbb{Z}/n$ of $\PGL_2(\mathbb{C})$ for $n\geq 2$ is, up to conjugation, generated by $[x_0 : x_1] \mapsto [\xi x_0 : x_1]$, where $\xi$ is a primitive root of unity of degree $n$. Hence there are two fixed points: $[0 : 1]$ and $[1 : 0]$. All other points have orbits of cardinality $n\geq 2$.

A subgroup of $\PGL_2(\mathbb{C})$ isomorphic to $(\mathbb{Z}/2)^2$ is, up to conjugation, generated by $[x_0 : x_1] \mapsto [-x_0 : x_1]$ and $[x_0 : x_1] \mapsto [x_1 : x_0]$. There are three orbits of size 2: $\{[0 : 1],[1 : 0]\}$, $\{[1 : 1],[1 : -1]\}$, $\{[1 : i],[1 : -i]\}$. Any other $G$-orbit on $\mathbb{P}^1$ has cardinality~$4$.
\end{proof}

\begin{proposition}
\label{prop-elliptic-curve}
Let $C$ be a smooth curve of genus $1$, and let $G$ be a finite abelian group such that $G\subset \Aut(C)$.  Then the following cases are possible $($we denote by $\sigma$ the antipodal involution on $C$; by ord.\ of CM we mean the order of the automorphism group of\, $C$ considered as an elliptic curve$)$:
\end{proposition}

\begin{center}
  \begin{table}[h!]
    {\renewcommand{\arraystretch}{1.1}
\begin{tabular}{ | m{1.3em} | m{4.0cm} | m{4cm} | m{1.8cm} | m{1.97cm} | } 
  \hline
   & $G$ & Min. length of an orbit & $\sigma\in G$ & Ord.\ of CM \\ 
  \hline
  (1) & $\mathbb{Z}/n\mathbb{Z}\times\mathbb{Z}/m\mathbb{Z}$ & $mn$ & & \\ 
  \hline
  \refstepcounter{foo}\otherlabel{t22}{2}{(2)} 
  & $\mathbb{Z}/2$ & $1$ & \checkmark & \\ 
    \hline
  \refstepcounter{foo}\otherlabel{t23}{3}{(3)} 
   & $(\mathbb{Z}/2)^2$ & $2$ & \checkmark & \\ 
    \hline
  \refstepcounter{foo}\otherlabel{t24}{4}{(4)} 
   & $(\mathbb{Z}/2)^3$ & $4$ & \checkmark & \\ 
    \hline
  \refstepcounter{foo}\otherlabel{t25}{5}{(5)} 
   & $\mathbb{Z}/4$ & $1$ & \checkmark & $4$ \\ 
    \hline
  \refstepcounter{foo}\otherlabel{t26}{6}{(6)} 
   & $\mathbb{Z}/2\times\mathbb{Z}/4$ & $2$ & \checkmark & $4$  \\ 
    \hline
  \refstepcounter{foo}\otherlabel{t27}{7}{(7)} 
  & $\mathbb{Z}/3$ & $1$ & & $6$ \\ 
    \hline
  \refstepcounter{foo}\otherlabel{t28}{8}{(8)} 
   & $(\mathbb{Z}/3)^2$ & $3$ & & $6$ \\ 
    \hline
  \refstepcounter{foo}\otherlabel{t29}{9}{(9)} 
  & $\mathbb{Z}/6$ & $1$ &  \checkmark & $6$ \\ 
    \hline
\end{tabular}}
\caption{Finite abelian groups acting on smooth curves of genus $1$}
\label{table-2}
\end{table}
\end{center}

\begin{proof}
  It is well known  that 
\[
\Aut(C)\simeq C_0\rtimes \mathbb{Z}/k\mathbb{Z}, 
\]
where $C_0$ is the
subgroup of translations that can be identified (after fixing the base point on $C$) with $C$, and $k\in\{2,4,6\}$. 
If $k>2$, then $C$ is said to have complex multiplication. 
If $G\subset C_0$, it follows that $\mathfrak{r}(G)\leq 2$ and $G$ acts by translations, so we are in the first case of the proposition. Hence we may assume that $G\not\subset C_0$. 
Consider the exact sequence
\begin{equation}
\label{eq-extension-elliptic-curve}
0\lra H\lra G \lra K \lra 0, 
\end{equation}
where $H = G\cap C_0$ and $K\subset \mathbb{Z}/k\mathbb{Z}$. 
We have $K =\mathbb{Z}/r$ for $r\in\{2,3,4,6\}$. 
We consider four cases.

\emph{Case $r=2$.}~ 
Note that $K=\mathbb{Z}/2$ acts on $C$ by the antipodal involution $\sigma$. One checks that the only translations that commute with $\sigma$ are translations by the $2$-torsion points. Hence we obtain $G\subset
(\mathbb{Z}/2)^3$ (note that we cannot get an element of order~$4$ with an extension \eqref{eq-extension-elliptic-curve}). Observe that the action of $G$ on $C$ is not free, since $\sigma$ has four fixed points on $C$. We obtain cases~\eqref{t22},~\eqref{t23},~\eqref{t24} in Table~\ref{table-2}.

\emph{Case $r=3$.}~
By the Lefschetz fixed-point formula, we see that $K$ has three fixed points on $C$. Note that the only translations that commute with $K$ are translations by $3$-torsion points, which form a subgroup $\mathbb{Z}/3\times\mathbb{Z}/3$ of $C_0$.
Observe that only a ``diagonal'' subgroup $\mathbb{Z}/3$ of $\mathbb{Z}/3\times\mathbb{Z}/3$ commutes with $K =\mathbb{Z}/3$. Thus, $G\subset\mathbb{Z}/3\times\mathbb{Z}/3$. 
If $G=\mathbb{Z}/3\times\mathbb{Z}/3$, then $G$ has an orbit of length $3$. If $G=\mathbb{Z}/3$, then $G$ fixes the origin. 
We obtain cases~\eqref{t27} and~\eqref{t28} in Table~\ref{table-2}.

\emph{Case $r=4$.}~ 
As in the case $r=2$, note that the only translations that commute with $K$ belong to the subgroup $\mathbb{Z}/2\times\mathbb{Z}/2$ of $C_0$ that consists of translations by $2$-torsion points. On the other hand, only a ``diagonal'' subgroup $\mathbb{Z}/2$ of $\mathbb{Z}/2\times\mathbb{Z}/2$ commutes with $K=\mathbb{Z}/4$. Thus, $G\subset \mathbb{Z}/2\times\mathbb{Z}/4$.
If $G=\mathbb{Z}/2\times\mathbb{Z}/4$, then $G$ has an orbit of length $2$. If $G=\mathbb{Z}/4$, then~$G$ fixes the origin. 
In fact, in the latter case, $G$ has exactly two fixed points, which follows, for example, from the Lefschetz fixed-point formula. 
This gives cases~\eqref{t25} and~\eqref{t26} in Table~\ref{table-2}.
 
\emph{Case $r=6$.}~ 
From considerations in the cases $r=2$ and $r=3$, it follows that $H$ is trivial in this case.  
Then~$G$ has only one fixed point on $C$.
We obtain case~\eqref{t29} in Table~\ref{table-2}. The proof is completed.
\end{proof}

\begin{corollary}
\label{cor-elliptic-type}
Let $C$ be a smooth curve of genus $1$, and let $G$ be a finite abelian group such that $G\subset \Aut(C)$. 
Then $G\subset \Cr_2(\mathbb{C})$.
\end{corollary}

\begin{proof}
  This follows from Proposition~\ref{prop-elliptic-curve} and Theorem~\ref{thm-surfaces}.
\end{proof}

\begin{example}
The maximal groups in Proposition~\ref{prop-elliptic-curve} that do not consist of translation elements are realized, for example, for the curves given by the following equations:
\begin{enumerate}\setlength{\itemsep}{2pt}
\item
$x_1^2+x_2^2+x_3^2+x_4^2=a_1x_1^2+a_2x_2^2+a_3x_3^2+a_4x_4^2=0\ \text{in}\ \mathbb{P}^3$ and $G=(\mathbb{Z}/2)^3$, 
\item
$x_1^3+x_2^3+x_3^3=0\ \text{in}\ \mathbb{P}^2$ and $G=(\mathbb{Z}/3)^2$, 
\item
$x_1^2+x_2^4+x_3^4=0\ \text{in}\  \mathbb{P}(2,1,1)$ and $G=\mathbb{Z}/2\times\mathbb{Z}/4$, 
\item
$x_1^2+x_2^3+x_3^6=0\ \text{in}\ \mathbb{P}(3,2,1)$ and $G=\mathbb{Z}/6$.
\end{enumerate}
\end{example}

\section{Group action on surfaces}
\label{sec-group-action-on-surfaces}
In this section, we prove some results on the action of finite abelian groups on algebraic surfaces.

\begin{proposition}
\label{prop-invariant-pencil}
Let $S$ be a projective surface, and let $G$ be a finite abelian group such that $G\subset \Aut(S)$. 
Assume that $S$ admits a $G$-invariant pencil $\mathcal{H}$ whose general element is an irreducible reduced rational curve. Then $G$ is of type~\eqref{cremona-plane-abelian-1},~\eqref{cremona-plane-abelian-2}, or~\eqref{cremona-plane-abelian-5} as in Theorem~\ref{thm-surfaces}.
\end{proposition} 

\begin{proof}
Resolving singularities, we may assume that $S$ is a smooth projective surface and $G\subset \Aut(S)$. 
Resolving the base locus of $\mathcal{H}$ in the $G$-equivariant category, we obtain a morphism $f\colon S'\to \mathbb{P}^1$, where $S'$ is a smooth projective surface and $G\subset \Aut(S')$. By the Bertini theorem, a general fiber of~$f$ is a smooth rational curve.
Then we have an exact sequence
\[
0\lra G_\eta \lra G\lra G_Z\lra 0, 
\]
where $G_Z$ faithfully acts on the base $Z=\mathbb{P}^1$, and so $G_Z\subset \Cr_1(\mathbb{C})$ and $G_\eta$ faithfully acts on the generic fiber $S'_\eta$, which is a Fano curve defined over the function field of the base~$\mathbb{C}(t)$. By Lemma~\ref{lem-lefschetz}, we see that $G_\eta\subset \Cr_1(\mathbb{C})$.  
Now the claim follows from Proposition~\ref{prop-abstact-extension-cr2}. We also note that the groups of types~\eqref{cremona-plane-abelian-1},~\eqref{cremona-plane-abelian-2}, and~\eqref{cremona-plane-abelian-5} can faithfully act on $\mathbb{P}^1\times\mathbb{P}^1$. 
\end{proof}

\begin{lem}
\label{lem-2-dim-CY-pairs}
Let $(S, \Delta)$ be an lc log Calabi--Yau pair of dimension $2$ with
the action of a finite group $G$. 
Assume that $\Delta\neq 0$. 
Then $S$ is $G$-birational either to a
$G$-del Pezzo surface, or to a $G$-conic bundle over a curve $Z$ such
that either $Z\simeq \mathbb{P}^1$, or $Z$ is a smooth elliptic curve.
\end{lem}

\begin{proof}
  Consider the minimal resolution $f\colon S'\to S$, which is automatically $G$-equivariant, and let $(S', \Delta')$ be the log pullback of $(S, \Delta)$.
Run a $G$-equivariant minimal model program on $K_{S'}$ to obtain a $G$-Mori fiber space $\pi\colon S''\to Z$, and let $(S'', \Delta'')$ be the pushdown of $(S', \Delta')$. We have the following diagram:
\begin{equation}
\label{eq-diag-minimal-model-ruled-surface}
\begin{tikzcd}
(S', \Delta')  \arrow{r}{g} \arrow[swap]{d}{f} & (S'', \Delta'')  \arrow{d}{\pi} \\
S & Z\rlap{.}
\end{tikzcd}
\end{equation}
Note that $\Delta''\neq 0$.
Assume that $Z$ is a curve. 
By the canonical bundle formula, we have   
\[
0\sim_{\mathbb{Q}} K_{S''} + \Delta'' \sim_{\mathbb{Q}} \pi^*(K_Z + \Delta_Z + M_Z), 
\]
where $\Delta_Z$ is a boundary $\mathbb{Q}$-divisor and $M_Z$ is a pseudo-effective $\mathbb{Q}$-divisor class. Since $M_Z$ is pseudo-effective, we have $\deg M_Z\geq 0$. Taking degrees, we have $\deg K_Z\leq 0$. Thus either $Z=\mathbb{P}^1$ and $\pi\colon S''\to Z$ is a $G$-conic bundle, or $Z$ is an elliptic curve and $\Delta_Z=0$, $M_Z\equiv 0$. 

Assume that $Z$ is a point. Then $S''$ is a $G$-del Pezzo surface; in particular,  $S$ is rational. The claim follows.
\end{proof}

\begin{remark}
\label{rem-albanese-map}
In the proof of Lemma~\ref{lem-2-dim-CY-pairs}, if $Z$ is an elliptic curve, then the ($G$-equivariant) map $\pi\circ g\colon S'\to Z$ is the Albanese map. This map induces a $G$-equivariant rational map $S\dashrightarrow Z$.
\end{remark}

\begin{corollary}
\label{cor-2-dim-cy-pair-pt}
Assume that a finite abelian group $G$ faithfully acts on an lc log Calabi--Yau pair $(S, \Delta)$ with $\dim S=2$ and $\Delta\neq 0$. Then $G$ is of product type.
\end{corollary}

\begin{proof}
By Lemma~\ref{lem-2-dim-CY-pairs}, the surface $S$ is either rational, or birationally ruled over an elliptic curve. In the former case, we have $G\subset\Cr_2(\mathbb{C})$, so $G$ is of product type. In the latter case,~$G$ is an extension of a group that faithfully acts on an elliptic curve (and so is isomorphic to a subgroup of $\Cr_2(\mathbb{C})$ by Corollary~\ref{cor-elliptic-type}) by a finite abelian subgroup of $\Cr_1(\mathbb{C})$. Hence, $G$ is of product type by Proposition~\ref{prop-abstact-extension}.
\end{proof}

\subsection{Simple elliptic and cusp singularities}
\label{subsec-simple-ell-cusp}
We recall some types of surface singularities.
By the classification (see \textit{e.g.}~\cite[Section~4.4]{KM98}),
strictly lc Gorenstein surface singularities are either simple elliptic or cusp singularities. 
\label{subsec-cusp-elliptic}
Recall that a normal surface $S$ has a \emph{simple elliptic singularity} at a point $P\in S$ if for the minimal resolution $f\colon \widetilde{S}\to S$, we have $(R^1f_*\OOO_{\widetilde{S}})_P=\mathbb{C}$ and the exceptional divisor $\Exc(f)$ is a smooth elliptic curve~$E$.  
By \cite[Proposition~4.54]{KM98}, we have 
\begin{equation}
\label{eq-elliptic-curve-embedding-dimension}
\dim T_P S = \max(3, -E^2). 
\end{equation}
For a simple elliptic singularity, there are the following possibilities (we also define a weight vector $w$ that we will use later): 
\begin{enumerate}
\item
  either $\mult_P(S)\geq 3$,
  and then $P\in S$ embeds into $\mathbb{C}^k$ for $k=-E^2$; in this case,  
\begin{equation}
\label{eq-simple-elliptic-sing-0}
w(x_1,\ldots,x_k)=(1,\ldots,1),  \quad -E^2\geq 3; 
\end{equation}
\item
or $P\in S$ embeds into $\mathbb{C}^3$, and after an analytic change of coordinates, the local equation $\phi=0$ of $P\in S$ is
\begin{equation}
\label{eq-simple-elliptic-sing-1}
x_1^2 + q(x_2, x_3) = 0,  \quad w(x_1,x_2,x_3)=(2,1,1),  \quad -E^2=2, 
\end{equation}
where $q$ is a homogeneous polynomial of degree $4$; 
\item
or $P\in S$ embeds into $\mathbb{C}^3$, and after an analytic change of coordinates,  the local equation $\phi=0$ of $P\in S$ is
\begin{equation}
\label{eq-simple-elliptic-sing-2}
x_1^2+x_2^3+ax_2x_3^4+bx_3^6=0,  \quad w(x_1,x_2,x_3)=(3,2,1),  \quad -E^2=1, 
\end{equation}
where $a,b\in\mathbb{C}$. 
\end{enumerate}

\label{hypersurface-cusp-sing}
A normal surface $S$ has a \emph{cusp singularity} $P\in S$ if $S$ is strictly lc at $P$ and, in the minimal resolution $f\colon \widetilde{S}\to S$, we have $(R^1f_*\OOO_{\widetilde{S}})_P=\mathbb{C}$, and the preimage of $P$ is either a cycle of smooth rational curves or an irreducible rational curve with one node. 

\begin{proposition}
\label{prop-small-orbit-on-elliptic-curve}
Let $P\in S$ be a germ of a simple elliptic singularity, and let $G$ be a finite abelian group such that $G\subset \Aut(P\in S)$.
Let $f\colon \widetilde{S}\to S$ be the minimal resolution that is automatically $G$-equivariant. 
Then there exists a $G$-orbit on an elliptic curve $E=f^{-1}(P)$ of cardinality at most $-E^2$. 
\end{proposition}

\begin{proof}
We see that $P\in S$ is embedded into $\mathbb{C}^k$, where $k$ depends on $-E^2$, as discussed in Section~\ref{subsec-simple-ell-cusp}. 
We show that the weighted blow-up $g\colon \widetilde{\mathbb{C}^k}\to \mathbb{C}^k$ with  weight $w$ is $G$-equivariant, where the weight $w$ is defined in Equations~\eqref{eq-simple-elliptic-sing-0}, \eqref{eq-simple-elliptic-sing-1}, or~\eqref{eq-simple-elliptic-sing-2} depending on $-E^2$. Also,~$g$ induces the minimal resolution $f\colon \widetilde{S}\to S$, where~$\widetilde{S}$ is the strict transform of $S$ on~$\widetilde{\mathbb{C}^k}$ (\textit{cf.} \cite[Theorem~4.57]{KM98}).

Assume that $-E^2\geq 3$. Then $g$ is the blow-up of the closed point $P\in \mathbb{C}^k$, so $g$ is $G$-equivariant. 
By Lemma~\ref{lem-faithful-action}, we have that $G$ faithfully acts on the tangent space $T_P S$. Since $G$ is abelian, there exists a $G$-invariant hyperplane in $T_P S$ that corresponds to a $G$-invariant hyperplane~$H$ in the $g$-exceptional divisor $\mathbb{P}^{k-1}$. 
Then the intersection of the $f$-exceptional elliptic curve $E\subset \mathbb{P}^{k-1}=\mathbb{P}(T_P S)$ with $H$ is $G$-invariant. Using the fact that 
$-E^2=k=H\cdot E$ (\textit{cf.} \cite[Proposition~4.54]{KM98}), we obtain the result.

Now we consider the case $-E^2=2$, the case $-E^2=1$ being analogous. If $-E^2=2$, then for the local equation $\phi=0$ of $P\in S$ as in Equation~\eqref{eq-simple-elliptic-sing-1}, we have $\rk(\phi_2)=1$, where $\phi_2$ is the quadratic term of $\phi$.  Hence we can choose a (unique up to scaling) coordinate $x_1$ on $\mathbb{C}^3$ such that $x_1^2=\phi_2$.  Since $G$ is abelian, the representation of $G$ in $(T_P S)^*$ splits, so there is a $G$-invariant complement of $\langle x_1 \rangle$ in $\mathbb{C}^3$ spanned by some coordinates $x_2, x_3$. Hence the blow-up $g$ with weights $w=(2,1,1)$ with respect to coordinates $x_1,x_2,x_3$ is a $G$-equivariant morphism. Then the intersection of the $f$-exceptional elliptic curve $E\subset \mathbb{P}(2,1,1)$, which is a Cartier divisor with a $G$-invariant Weil divisor $H=\{x_2=0\}\subset \mathbb{P}(2,1,1)$, is $G$-invariant. Since $H\cdot E=2$, the claim follows.
\end{proof}

\section{Warm-up: Finite abelian groups in \texorpdfstring{$\boldsymbol{\Cr_2(\mathbb{C})}$}{Cr\_2(C)}}
\label{sec-warm-up}
In this section, we re-prove Theorem~\ref{thm-surfaces} using the methods that we will apply in dimension $3$ later on. 

\begin{proof}[Proof of Theorem~\ref{thm-surfaces}]
Let $G\subset \Cr_2(\mathbb{C})$ be a finite abelian group. We will obtain the list of possibilities for isomorphism classes of $G$ in three steps. 

\emph{Step~$1$. Reduction to Fano}.~
Passing to a $G$-equivariant resolution of singularities, we may assume that $G\subset \Aut(X)$, where~$X$ is a smooth projective rational surface. Applying the $G$-equivariant minimal model program, we may assume that either $X$ has a $G$-conic bundle structure $\pi\colon X\to Z=\mathbb{P}^1$, or $X$ is a $G$-del Pezzo surface. In the former case, we have an exact sequence
\[
0\lra G_f \lra G\lra G_Z \lra 0, 
\]
where $G_Z$ faithfully acts on the base $Z= \mathbb{P}^1$ and $G_f$ faithfully acts on the schematic generic fiber $\mathbb{P}^1_{\overline{\mathbb{C}(Z)}}$, where $\overline{\mathbb{C}(Z)}$ is the algebraic closure of the function field of $Z$.  By Lemma~\ref{lem-lefschetz}, we see that $G_f$ is isomorphic to a subgroup of $\Cr_1(\mathbb{C})$.  Proposition~\ref{prop-abstact-extension-cr2} implies that $G$ is of type~\eqref{cremona-plane-abelian-1}, ~\eqref{cremona-plane-abelian-2}, or~\eqref{cremona-plane-abelian-5} as in Theorem~\ref{thm-surfaces} in this case. So we may assume that $G\subset \Aut(X)$, where $X$ is a smooth del Pezzo surface with $\rho^G(X)=1$.

Note that $G$ acts on the linear system $\mathrm{H}^0(X, \OOO(-K_X))$. Since $G$ is abelian and we have $h^0(X, \OOO(-K_X))\geq 2$, there exist at least two $G$-invariant anti-canonical elements, say $C$ and $C'$. 
We consider two possibilities: when one of them is singular, and when both are smooth. 

\emph{Step~$2$. Singular anti-canonical curve}.~
Let $(X, C)$ be a $G$-invariant pair, where $X$ is a smooth del Pezzo surface and $C$ is an anti-canonical element; that is, $C\sim -K_X$. We assume that $C$ is singular.  
If the pair $(X, C)$ is not lc, then by a standard argument (\textit{cf.} {\cite[Lemma~2.9]{Pr11}}), 
it is possible to construct a $G$-invariant lc-center of dimension $0$, which gives us a $G$-fixed point. In this case, by Corollary~\ref{cor-faithful-action}, we obtain $\mathfrak{r}(G)\leq 2$, so $G$ is of type~\eqref{cremona-plane-abelian-1} as in Theorem~\ref{thm-surfaces}.

Thus, we may assume that $(X, C)$ is lc. By adjunction, $C$ has arithmetic genus $p_a(C)=~1$. 
We claim that~$C$ is either a nodal curve, or a cycle that consists of (at least two) smooth rational curves. Indeed, if $C$ is irreducible, then since by our assumption $C$ is singular, $(X, C)$ is lc, and $p_a(C)=1$, it follows that $C$ is a rational curve with exactly one node. If $C=\sum C_i$ is reducible (note that $C$ is connected since $C\sim -K_X$ is ample) and $(X, C)$ is an lc CY pair, then by adjunction we have that $(C_i, \sum_{j\neq i} C_j)$ is an lc CY pair. So each $C_i$ is a smooth rational curve and $\sum_{j\neq i} C_j|_{C_i}$ is a union of two points. We conclude that $C=\sum C_i$ is a cycle of smooth rational curves in this case. 

Assume that $C$ is a nodal curve. Then the node is a $G$-fixed point, and we are done by Corollary~\ref{cor-faithful-action}.  Assume that $C$ is a cycle of smooth rational curves $C=\sum_{i=1}^N C_i$, where $N\geq 2$ (if $N=2$, we have two smooth rational curves transversally intersecting in two points).  Without loss of generality, we may assume that $C_i\cap C_{i+1}\neq \emptyset$ for $1\leq i\leq N-1$ and $C_N\cap C_{1}\neq \emptyset$.  Note that the action of $G$ on $C=\sum C_i$ is faithful. Indeed, if a non-trivial element $g\in G$ acts on $C$ trivially, then $g$ trivially acts on the tangent space to a singular point $C_i\cap C_{i+1}$, which contradicts Lemma~\ref{lem-faithful-action}. Thus, we have an exact sequence
\[
0\lra H \lra G \lra G_{C} \lra 0, 
\]
where $G_{C}$ faithfully acts on the set $\{C_i\}$, and $H\subset \prod_i \Aut(C_i)$ with $C_i\simeq \mathbb{P}^1$. 
If the group $G_{C}$ acts on the set of components $\{C_i\}$ transitively, then applying Lemma~\ref{lem-diagonal-in-product}, we conclude that $H\subset \Aut(C_i)=\Cr_1(\mathbb{C})$. 
Since~$C$ is a cycle of curves, by Lemma~\ref{lem-action-on-spheres-1}, we have that $G_C$ is either cyclic or isomorphic to $(\mathbb{Z}/2)^2$. In both cases, $G_C\subset \Cr_1(\mathbb{C})$. By Proposition~\ref{prop-abstact-extension-cr2}, we conclude that $G$ is of type~\eqref{cremona-plane-abelian-1},~\eqref{cremona-plane-abelian-2}, or~\eqref{cremona-plane-abelian-5} as in Theorem~\ref{thm-surfaces} in this case.

Now assume that $C=C'+C''$, where both $C'$ and $C''$ are $G$-invariant curves. Since $\rho^G(X)=1$, both $C'$ and $C''$ are proportional to $-K_X$. In particular, both $C'$ and $C''$ are ample divisors, and $-K_X\sim C'+C''$, so the index of $-K_X$ is greater than~$1$. By the classification of del Pezzo surfaces, we conclude that $X$ is isomorphic either to $\mathbb{P}^1\times \mathbb{P}^1$ or to $\mathbb{P}^2$. If both $C'$ and $C''$ are reducible, then $X\simeq\mathbb P^1\times\mathbb P^1$ and $C'\sim C''\sim(1,1)$. The $G$-invariant pencil $\langle C',C''\rangle$ has an irreducible rational general member, so Proposition~\ref{prop-invariant-pencil} applies. Otherwise, after interchanging $C'$ and $C''$, we may assume that $C'$ is an irreducible $G$-invariant rational curve. We can write an exact sequence
\[
0\lra \mathbb{Z}/m \lra G\lra G_{C'}\lra 0,
\]
where $G_{C'}\subset \Aut(C')=\Cr_1(\mathbb{C})$ and $\mathbb{Z}/m$ faithfully acts on the normal bundle to $C'$. 
By Proposition~\ref{prop-abstact-extension-cr2}, we conclude that $G$ is of type~\eqref{cremona-plane-abelian-1},~\eqref{cremona-plane-abelian-2}, or~\eqref{cremona-plane-abelian-5} as in Theorem~\ref{thm-surfaces} in this case as well.

\emph{Step $3$. Smooth anti-canonical curves}.~
Let $X$ be a smooth $G$-del Pezzo surface, and assume that $C$ and $C'$ are smooth $G$-invariant anti-canonical elliptic curves. 
Consider an exact sequence
\[
0\lra G_N \lra G\lra G_C \lra 0, 
\]
where $G_C$ faithfully acts on the elliptic curve $C$, and where $G_N$ fixes $C$ pointwise and faithfully acts on the normal bundle $N_{C/X}$. In particular, $G_N$ is cyclic, so $G_N=\mathbb{Z}/n$ for some $n\geq 1$. If $\mathfrak{r}(G_C)\leq 1$, we have $\mathfrak{r}(G)\leq 2$, and thus $G$ is of type~\eqref{cremona-plane-abelian-1} as in Theorem~\ref{thm-surfaces}. Hence we may assume that $\mathfrak{r}(G_C)\geq 2$. 
Since $C\cap C'$ is $G$-invariant,~$G_C$ admits an orbit on $C$ of length at most $(-K_X)^2\leq 9$. By Proposition~\ref{prop-elliptic-curve}, this implies that~$G_C$ is a subgroup of one of the following groups: 
\begin{equation}
\label{eq-extra-groups-cremona2}
(\mathbb{Z}/3)^2,  \quad (\mathbb{Z}/2)^3,  \quad \mathbb{Z}/4\times \mathbb{Z}/2.
\end{equation}
Consider an analogous exact sequence for the elliptic curve $C'$:
\[
0\lra G_{N'} \lra G\lra G_{C'} \lra 0,
\]
where $G_{C'}$ faithfully acts on $C'$, and where $G_{N'}$ fixes $C'$ pointwise and faithfully acts on the normal bundle $N_{C'/X}$. 
Note that $G_{N'}$ is cyclic, so $G_{N'}=\mathbb{Z}/n'$ for some $n'\geq 1$. As above, we may assume that $G_{C'}$ is a subgroup of one of the groups \eqref{eq-extra-groups-cremona2}. 

We claim that $G_N\cap G_{N'}=1$. If $g$ belongs to this intersection, choose $P\in C\cap C'$. The fixed locus of the finite-order automorphism $g$ is smooth and contains two distinct curve germs through $P$ on the smooth surface $X$. So $g$ acts trivially on $T_PX$ and is trivial by Lemma~\ref{lem-faithful-action}.
Hence, $G_N\subset G_{C'}$ and $G_{N'}\subset G_{C}$. Consider a diagram
\begin{equation}
\label{big-diagram}
\begin{tikzcd}
& & 0 \ar[d]  & & \ \\
& & \mathbb{Z}/n' \ar[d] \ar[dr, hook] & & \ \\
0 \ar[r] & \mathbb{Z}/n \ar[r] \ar[dr, hook] & G \ar[r] \ar[d] & G_C \ar[r] & 0\\
& & G_{C'} \ar[d] & & \ \\
& & 0\rlap{.} & & \ 
\end{tikzcd}
\end{equation}
We consider the case $G_C=G_{C'}=(\mathbb{Z}/3)^2$, the other cases being similar. It follows that $n\leq 3$. Since $G$ is an abelian extension of $(\mathbb{Z}/3)^2$ by $\mathbb{Z}/n$, either $G=(\mathbb{Z}/3)^3$ (in which case the exact row in Diagram~\eqref{big-diagram} splits), or $\mathfrak{r}(G)\leq 2$. Similar arguments imply that either $\mathfrak{r}(G)\leq 2$, or $G$ is a subgroup of one of the groups 
\begin{equation}
\label{eq-three-exceptional-groups}
(\mathbb{Z}/3)^3, \quad  (\mathbb{Z}/2)^4,  \quad (\mathbb{Z}/4)^2\times \mathbb{Z}/2.
\end{equation}
This gives us cases~\eqref{cremona-plane-abelian-3},~\eqref{cremona-plane-abelian-4},~\eqref{cremona-plane-abelian-5} of Theorem~\ref{thm-surfaces}. 
Moreover, observe that if $G$ is isomorphic to one of the three groups in the list \eqref{eq-three-exceptional-groups}, then in Diagram \eqref{big-diagram},  both the exact row and the exact column split. 
The theorem is proven.
\end{proof}

\begin{remark} 
Our proof of Theorem~\ref{thm-surfaces} gives us more than just a list of groups. The given argument implies the following: 
\begin{itemize}
\item
If $G$ faithfully acts on a Mori fiber space $X\to Z$, where $X$ is a rational surface and the base $Z$ is non-trivial (so that $Z=\mathbb{P}^1$), then $G$ has type~\eqref{cremona-plane-abelian-1},~\eqref{cremona-plane-abelian-2}, or~\eqref{cremona-plane-abelian-5} as in Theorem~\ref{thm-surfaces}. In particular, $G=H\times K$, where $H, K\in \Cr_1(\mathbb{C})$.
\item
  If $G$ faithfully acts on a smooth del Pezzo surface $S$ and there exists a singular $G$-invariant curve $C$ with $C\sim -K_S$, then $G=H\times K$, where $H, K\in \Cr_1(\mathbb{C})$.
\end{itemize}
Our goal is to generalize these observations to dimension $3$. In the next corollary, we explain the relation between groups of types~\eqref{cremona-plane-abelian-3} and~\eqref{cremona-plane-abelian-4} as in Theorem~\ref{thm-surfaces} and certain elliptic curves with complex multiplication.
\end{remark}

\begin{corollary}
A faithful action of a group $G$ of type~\eqref{cremona-plane-abelian-3} or~\eqref{cremona-plane-abelian-4} as in Theorem~\ref{thm-surfaces} on a rational surface is conjugate to an action of\, $G$ on a smooth del Pezzo surface $X$ such that each of its $G$-invariant anti-canonical elements is a smooth elliptic curve. Moreover, such an elliptic curve admits an automorphism of order $N$, and there are the following possibilities $($here $j(C)$ denotes the $j$-invariant of\, $C)$:  
\begin{enumerate}\setlength{\itemsep}{2pt}
\item
$G=(\mathbb{Z}/3)^3$, $(-K_X)^2=3$, $N=6$, $j(C)=0$;  
\item
$G=(\mathbb{Z}/4)^2\times\mathbb{Z}/2$, $(-K_X)^2=2$, $N=4$,  $j(C)=1728$.
\end{enumerate}
\end{corollary}

\begin{proof}
Let $C$ and $C'$ be distinct $G$-invariant anti-canonical elements on $X$. By our assumption on the group~$G$ and the proof of Theorem~\ref{thm-surfaces},  we know that $C$ and $C'$ are smooth elliptic curves. 
In the notation of Diagram \eqref{big-diagram}, the elements of $\mathbb{Z}/n$ (resp.\ $\mathbb{Z}/n'$) fix $C$ (resp.\ $C'$) pointwise. This implies that the generator of $\mathbb{Z}/n$ (resp.\ $\mathbb{Z}/n'$) acts on $C'$ (resp.\ $C$) as an automorphism with a fixed point. In particular, it is not a translation.

Now observe that $C$ and $C'$ have transversal intersection at all points of $C\cap C'$. Indeed, let $P\in C\cap C'$. Since $\mathbb{Z}/n$ (resp.\ $\mathbb{Z}/n'$) acts on $C$ (resp.\ $C'$) trivially, it acts trivially on $T_PC$  (resp.\ $T_PC'$) as well. The tangent lines are distinct: otherwise Lemma~\ref{lem-faithful-action} would give a faithful one-dimensional representation of the non-cyclic group $G_N\times G_{N'}$, a contradiction. Hence $T_PX = T_PC\oplus T_PC'$, so the intersection of $C$ and $C'$ is transversal at $P$. Hence, $m=C\cdot C'$ is equal to the cardinality of the set $C\cap C'$. We also have $(-K_X)^2=m$ since $C\sim C'\sim -K_X$. 

The set $C\cap C'$ is contained in the fixed locus of $\mathbb{Z}/n$ on $C'$. By our assumption on the group $G$, we have $n=3$ or $n=4$; the generator $g$ of $\mathbb{Z}/n$ acts via an automorphism of order $n$. By transversality and Lemma~\ref{lem-faithful-action}, the action of $\mathbb Z/n$ on $C'$ is faithful. If $n=3$, the element $g$ has three fixed points on $C'$, and if $n=4$, then $g$ has two fixed points on $C'$ (\textit{cf.}~Proposition~\ref{prop-elliptic-curve}). Since $C\cap C'$ is $G_{C'}$-invariant and $G_{C'}$ is isomorphic to $(\mathbb{Z}/3)^2$ or $\mathbb{Z}/4\times\mathbb{Z}/2$, respectively, Proposition~\ref{prop-elliptic-curve} gives $m\geq3$ and $m\geq2$, respectively. Thus $m=3$ and $m=2$, respectively. It is well known that an elliptic curve that admits an automorphism of order $3$ (and hence of order $6$ as well) has $j$-invariant equal to $0$, and an elliptic curve that admits an automorphism of order $4$ has $j$-invariant equal to $1728$. This finishes the proof.
\end{proof}

\section{Terminal singularities}
\label{subsec-terminal-sing}

We recall some facts on the classification of $3$-dimensional terminal singularities from \cite{Mo85,Re85}. 
Let $P\in X$ be a germ of a $3$-dimensional terminal singularity. Then the singularity is isolated: $\Sing(X) = \{P\}$. The \emph{index} of $P\in X$ is the minimal positive integer $r$ such that $rK_X$ is Cartier. If $r = 1$, then $P\in X$ is Gorenstein. In this case, $P\in X$ is analytically isomorphic to a hypersurface singularity in $\mathbb{C}^4$ of multiplicity~$2$. Moreover, any Weil $\mathbb{Q}$-Cartier divisor $D$ on $P\in X$ is Cartier. Also, in this case, $P\in X$ is a compound du Val singularity, which means that the general hyperplane section $H$ that contains $P$ is a surface with a du Val singularity at $P$. We say that $P\in X$ has type $cA_n$, $cD_n$, $cE_n$, respectively, if $H$ is a du Val singularity of type $A_n$, $D_n$, $E_n$, respectively. 

 If $r > 1$, then there is a cyclic, \'etale outside~$P$, covering 
\[
\pi\colon \widetilde{P}\in \widetilde{X} \lra X\ni P
\] 
of degree $r$ such that $\widetilde{P}\in \widetilde{X}$ is a Gorenstein terminal singularity (or a smooth point). The map $\pi$ is called the \emph{index $1$ cover} of $P\in X$, and it is defined canonically.

\begin{thm}[\textit{cf.} {\cite{Mo85, Re85}}]
\label{thm-terminal-sing}
In the above notation, 
$\widetilde{P}\in \widetilde{X}$ is $\mathbb{Z}/r$-equivariantly analytically isomorphic to a hypersurface in $\mathbb{C}^4$ with $\mathbb{Z}/r$-semi-invariant coordinates $x_1,x_2,x_3,x_4$, the action of\, $\mathbb{Z}/r$ is given by
\[
(x_1,x_2,x_3,x_4)\longmapsto (\xi^{a_1}x_1,\xi^{a_2}x_2,\xi^{a_3}x_3,\xi^{a_4}x_4)
\]
for some primitive $\supth{r}$ root of unity $\xi$, and one of the following holds:
\begin{enumerate}
\item
$(a_1,\ldots,a_4) = (1,-1,a_3,0)\bmod r$ and $\gcd(a_3, r) = 1$,
\item
$r = 4$ and $(a_1,\ldots,a_4) = (1,-1,1,2)\bmod 4$.
\end{enumerate} 
\end{thm}

Now, we consider the action of a finite abelian group on a $3$-dimensional terminal singularity germ. Recall the following result.

\begin{proposition}[\textit{cf.} {\cite[Lemma 2.4]{Pr11} and \cite[Lemma 2.1]{Pr14}}]
\label{prop-p-action-on-terminal-point}
Let $P\in X$ be a germ of a threefold terminal singularity, and let
$G \subset \Aut(P\in X)$ be a finite abelian group. Then
$\mathfrak{r}(G_2) \leq 4$ and $\mathfrak{r}(G_p)\leq 3$ for $p\geq 3$. 
Moreover, if\, $\mathfrak{r}(G_2)=4$, then $P\in X$ is not a cyclic quotient singularity.
\end{proposition}

The main goal of this section is to improve this result as follows.

\begin{thm}
\label{thm-action-on-terminal-point}
Let $P\in X$ be a germ of a threefold terminal singularity, and let $G \subset \Aut(P\in X)$ be a finite abelian group. Then either $\mathfrak{r}(G)\leq 3$, or
\[
G= (\mathbb{Z}/2)^2\times\mathbb{Z}/2n\times\mathbb{Z}/2m
\]
for $n,m\geq 1$. 
Moreover, if\, $\mathfrak{r}(G)=4$, then $P\in X$ is a Gorenstein singularity of type $cA$. In particular, in both cases, $G$ is of product type.
\end{thm}

Before proving this theorem, we deduce from it the following useful result. 

\begin{corollary}
\label{cor-fixed-point}
Let $X$ be a threefold with terminal singularities, and let $G\subset \Aut(X)$ be a finite abelian group. 
Assume that there exists a $G$-fixed point $P\in X$. Then $G$ is of product type.
\end{corollary}

\begin{proof}
Assume that $P\in X$ is smooth. Then $G$ is of product type by Corollary~\ref{cor-faithful-action}.  If $P\in X$ is singular, then by Theorem~\ref{thm-action-on-terminal-point}, we conclude that $G$ is of product type.
\end{proof}

The latter case in Theorem~\ref{thm-action-on-terminal-point} is realized, as the following example shows.

\begin{example}
\label{ex-maximal-terminal-action}
Consider a singular point $P\in \mathbb{C}^4$ given by the equation
\[
x_1^2 + x_2^2 + x_3^{n} + x_4^{m} = 0 
\]
for $2\leq n\leq m$. This point is an isolated compound du Val singularity of type $cA_{n-1}$, and hence it is terminal; see \cite{KS-B88}. It admits an action of the group $G=(\mathbb{Z}/2)^2\times\mathbb{Z}/n\times\mathbb{Z}/m$, where the generators of this group act on the coordinates in $\mathbb{C}^4$ via multiplication by roots of unity.
\end{example}

\begin{example}
Consider an ordinary double point $P\in X$ given by the equation
\[
x_1 x_2 + x_3 x_4 = 0.
\]
It is a toric singularity. In particular, it admits an action of $(\mathbb{Z}/n)^3$ for any $n\geq 1$.
\end{example}

We will need the following well-known version of the Morse lemma.

\begin{lem}[\textit{cf.} {\cite[11.1]{AGV85}}]
\label{lem-corank-second-differenetial}
Let $f(x_1, \ldots, x_n)$ be a germ of a holomorphic function such that $0\in \mathbb{C}^n$ is its critical point. Assume that the quadratic term of $f$ in the  expansion near $0$ has rank $k\leq n$. Then, up to a holomorphic change of coordinates, we have
\[
f(x_{1},\ldots, x_n) = x_1^2 + \cdots + x_k^2 + g(x_{k+1},\ldots, x_n), 
\]
where $g(x_{k+1},\ldots, x_n)$ is a germ of a holomorphic function whose quadratic term is equal to zero.
\end{lem}

\subsection{Gorenstein singularities}
Let $P\in X$ be a $3$-dimensional terminal Gorenstein singularity germ, and let $G\subset \Aut(P\in X)$ be a finite abelian group.  We may assume that $X\subset \mathbb{C}^4$, the point $P$ coincides with the origin $0\in \mathbb{C}^4$, and $G \subset \GL(4, \mathbb{C})$. It follows that $\mathfrak{r}(G)\leq 4$.  We can write the analytic equation of $X$ in the form $\phi = \sum_{i\geq 1} \phi_i$, where the $\phi_i$ are homogeneous polynomials of degree~$i$. Since $G$ is abelian, we may assume that the action of $G$ is diagonal with respect to the coordinates $x_1, \ldots, x_4$. By \cite{Mo85,Re85}, we can choose $\phi$ to be $G$-semi-invariant. Since~$P$ is a terminal singularity, we have $\phi_2\neq 0$. By the classification of terminal singularities, if $\phi_3 = 0$, then $\rk(\phi_2) \geq 2$.

\begin{defin}
\label{defin-standard-action}
We say that the action of the group $H=(\mathbb{Z}/m)^n\subset \GL(n, \mathbb{C})$ for some $m\geq 1$, $n\geq 1$ is \emph{standard} with respect to the coordinates $x_1, \ldots, x_n$ on $\mathbb{C}^n$ if the $\supth{i}$ generator of~$H$ for $1\leq i\leq n$ acts as follows:
\[
x_1\longmapsto x_1, \quad \ldots \quad  x_i\longmapsto \xi_m x_i,  \quad \ldots  \quad x_n\longmapsto x_n, 
\]  
where $\xi_m$ is a primitive root of unity of degree $m$.
\end{defin}

For a group $H=(\mathbb{Z}/m)^n\subset \GL(n, \mathbb{C})$ with $m\geq 1$, $n\geq 1$, 
after a change of coordinates on $\mathbb{C}^n$ and a change of generators of $H$, the action of~$H$ becomes {standard} on $\mathbb{C}^n$ with respect to the coordinates $x_1, \ldots, x_n$.

\begin{proposition}
\label{prop-action-on-gorenstein-points}
Let $P\in X$ be a $3$-dimensional terminal Gorenstein singularity germ given by the equation
\[
\phi(x_1,x_2,x_3,x_4)=\sum \phi_i(x_1,x_2,x_3,x_4)=0
\]
in $\mathbb{C}^4$ with coordinates $x_1,\ldots, x_4$. 
Let $G$ be a finite abelian group such that $G\subset \Aut(P\in X)$. 
Assume that $\mathfrak{r}(G_2)=4$.
Then
\begin{itemize} 
\item
up to an analytic change of coordinates, $\phi$ depends on $x_i^2$ for $1\leq i\leq 4$; that is, $\phi(x_1,x_2,x_3,x_4)=\psi(x_1^2,x_2^2,x_3^2,x_4^2)$ for some function $\psi$;  
\item
$\rk(\phi_2)\geq 2$, and $P\in X$ is a singularity of type $cA$;  
\item
$\mathfrak{r}(G_p)\leq 2$ for any $p>2$; in particular,  
\[
G= (\mathbb{Z}/2)^2\times\mathbb{Z}/2n\times\mathbb{Z}/2m\quad  \text{for } n,m\geq 1,
\] 
and hence $G$ is of product type.  
\end{itemize}
\end{proposition}

\begin{proof}
Since $\mathfrak{r}(G_2)=4$, there exists a subgroup $H=(\mathbb{Z}/2)^4\subset G$.  We may assume that the action of $H$ is standard with respect to~$x_i$ for $1\leq i\leq 4$.  Since the hypersurface $X=\{\phi=0\}$ is irreducible, the standard action of $H$ forces $\phi$ to be invariant under changing the sign of each $x_i$. Hence $\phi$ depends on~$x_i^2$ for $1\leq i\leq 4$; that is, $\phi(x_1,x_2,x_3,x_4)=\psi(x_1^2,x_2^2,x_3^2,x_4^2)$ for some $\psi$.  In particular, $\phi_i=0$ for odd~$i$.  By the classification of terminal singularities, $\phi_2\neq 0$, and since $\phi_3 = 0$, we have $\rk(\phi_2) \geq 2$.  By Lemma~\ref{lem-corank-second-differenetial}, it follows that $P\in X$ is a singularity of type $cA$.

Without loss of generality and after rescaling coordinates, we may assume that $\phi_2$ contains the terms $x_1^2+x_2^2$. We show that there exist a power of $x_3$ and a power of $x_4$ appearing in $\phi$ with non-zero coefficient. Indeed, if a power of $x_3$ does not belong to $\phi$, then the line $x_1=x_2=x_4=0$ is singular on $X$, and similarly for a power of $x_4$. Thus, up to a rescaling of coordinates, we have
\begin{equation}
\label{eq-gorenstein-rk-2}
\phi(x_1,x_2,x_3,x_4) = x_1^2+x_2^2+x_3^{2k}+x_4^{2l}+\cdots
\end{equation}
for some $l\geq k\geq 1$. Then the polynomial 
\[
\phi'(x_1,x_2,x_3,x_4) = x_1^2+x_2^2+x_3^{2k}+x_4^{2l} 
\]
is also $G$-semi-invariant. 
Let $d=\gcd(k,l)$ and use the auxiliary $G$-surface
$S=\{x_1^2+x_2^2+x_3^{2k}+x_4^{2l}=0\}\subset\mathbb{P}(kl,kl,l,k)$.
The pencil $\mathcal H=\langle x_3^{k/d},x_4^{l/d}\rangle$ is $G$-invariant. Since $k/d$ and $l/d$ are coprime, its general member is irreducible and reduced; its normalization is
$\{x_1^2+x_2^2+c\,u^{2kl/d}=0\}\subset\mathbb{P}(kl/d,kl/d,1)$, hence is rational.
We have the  exact sequence
\begin{equation}
\label{eq-exact-sequence-blow-up-of-point}
0 \lra G_N \lra G \lra G_S \lra 0, 
\end{equation}
where $G_S$ faithfully acts on $S$ and $G_N$ is a finite subgroup of $\mathbb{C}^\times$, and hence is cyclic.
Since~$S$ admits a $G$-invariant pencil of rational curves, by Proposition~\ref{prop-invariant-pencil}, we see that $G_S$ has type~\eqref{cremona-plane-abelian-1},~\eqref{cremona-plane-abelian-2}, or~\eqref{cremona-plane-abelian-5} as in Theorem~\ref{thm-surfaces}.
Now we apply Lemma~\ref{cor-r2-4-rp-2} to $G_S$ and $G_N$ to conclude that $G=(\mathbb{Z}/2)^2\times\mathbb{Z}/2n\times\mathbb{Z}/2m$ for some $n,m\geq 1$. 
This finishes the proof.
\end{proof}

\subsection{Non-Gorenstein singularities}
Assume that a finite abelian group $G$ faithfully acts on a $3$-dimensional terminal singularity germ $P\in X$ of index $r>1$.  By Theorem~\ref{thm-terminal-sing}, there exists a map $\pi\colon \widetilde{X}\to X$, where $\widetilde{P}\in \widetilde{X}\subset\mathbb{C}^4$ is a germ of Gorenstein terminal singularity and $\pi$ is a quotient by the action of the group $\mathbb{Z}/r$ for some $r>1$. Since $\pi|_{\widetilde{X}\setminus\{0\}}\colon\widetilde{X}\setminus\{0\}\to X\setminus\{0\}$ is a topological universal cover, there exists a lifting $\widetilde{G}\subset \Aut(\widetilde{X})$ fitting in the exact sequence
\begin{equation}
\label{eq-lift}
0\lra \mathbb{Z}/r \lra \widetilde{G} \lra G \lra 0, 
\end{equation}
where 
$\widetilde{G}\subset\GL_4(\mathbb{C})$.

\begin{proposition}[\textit{cf.} {\cite[Lemma 2.4]{Pr11} and \cite[Lemma 2.1]{Pr14}}]
\label{prop-action-on-non-gorenstein-points}
Let $P\in X$ be a $3$-dimensional terminal singularity germ, and let $G$ be a finite abelian group such that $G\subset \Aut(P\in X)$. 
Then either the lifting $\widetilde{G}$ of\, $G$ is abelian, or $r=2$.
\end{proposition}

\begin{proof}
Assume that the group $\widetilde{G}$ is not abelian. Clearly, for the Gorenstein index of $P\in X$, we have $r>1$.  We claim that $r=2$ in this case. Assume that $r>2$.  In the notation of Theorem~\ref{thm-terminal-sing}, the subspace of $\mathbb{C}^4$ defined as
\[
T=\mathbb{C}^3=\{x_4=0\}
\] 
is $\widetilde{G}$-invariant. The group $\widetilde{G}$ permutes the eigenspaces of $\mathbb{Z}/r$ on $T$. Since $r>~2$, according to Theorem~\ref{thm-terminal-sing}, we see that $T$ contains a $\widetilde{G}$-invariant $2$-dimensional subspace~$T_2$. Let us denote by $T_1$ the $\widetilde{G}$-invariant subspace of $T$ with the property $T=T_1\oplus T_2$.

By Theorem~\ref{thm-terminal-sing}, the group $\mathbb{Z}/r$ acts on $T_1\subset T$ faithfully.
Observe that the derived subgroup $[\widetilde{G},\widetilde{G}]$ is a non-trivial subgroup of $\widetilde{G}$ contained in $\mathbb{Z}/r$; \textit{cf.}~the exact sequence~\eqref{eq-lift}. Hence $[\widetilde{G},\widetilde{G}]$ is abelian and acts on $T_1$ faithfully. On the other hand, the kernel of the homomorphism $\widetilde{G}\to\GL(T_1)=\mathbb{C}^\times$ contains $[\widetilde{G},\widetilde{G}]$, which gives a contradiction.
This shows that $r=2$.
\end{proof}

\begin{proposition}
\label{prop-terminal-point-abelian}
Let $P\in X$ be a $3$-dimensional terminal singularity germ of Gorenstein index $r>1$, and let~$G$ be a finite abelian group such that $G\subset \Aut(P\in X)$.
Assume that the lifting $\widetilde{G}$ is abelian. Then
$\mathfrak{r}(G)\leq 3$.
\end{proposition}

\begin{proof}
Since $G$ is a quotient of $\widetilde G$, we have $\mathfrak r(G)\leq\mathfrak r(\widetilde G)$. Thus we are done if $\mathfrak r(\widetilde G)\leq3$, and we may assume that $\mathfrak{r}(\widetilde{G})=4$. Then $\widetilde P$ is singular by Corollary~\ref{cor-faithful-action}. Moreover, Proposition~\ref{prop-p-action-on-terminal-point} gives $\mathfrak r(\widetilde G_2)=4$. Hence by Proposition~\ref{prop-action-on-gorenstein-points}, we have that $\widetilde{P}\in\widetilde{X}$ is a cA point, and 
\[
\widetilde{G}=(\mathbb{Z}/2)^2\times\mathbb{Z}/2n\times\mathbb{Z}/2m \quad \text{for } n, m\geq 1.
\] 
We can choose coordinates in such a way that the equation 
\[
\phi(x_1,x_2,x_3,x_4)=\sum \phi_i(x_1,x_2,x_3,x_4)
\]
of $\widetilde{X}\subset\mathbb{C}^4$ is $\widetilde{G}$-semi-invariant, and hence $\mathbb{Z}/r$-semi-invariant, where the covering map $\pi\colon \widetilde{X}\to X$ is the quotient by the action of $\mathbb{Z}/r$. 
We may assume that the action of a subgroup $(\mathbb{Z}/2)^4\subset \widetilde{G}$ is standard (see Definition~\ref{defin-standard-action}) with respect to the coordinates $x_1,\ldots,x_4$. 
Hence, $\phi$ depends on $x_i^2$,  {that is, $\phi(x_1,x_2,x_3,x_4)=\psi(x_1^2,x_2^2,x_3^2,x_4^2)$ for some $\psi$,} and $2\leq \rk(\phi_2)\leq 4$ (\textit{cf.}~Proposition~\ref{prop-action-on-gorenstein-points}). 
We analyze cases according to the classification of terminal singularities as in \cite[Section~6.1]{Re85}.

\begin{enumerate}[wide,label=\textit{Case}~\arabic*:, ref=\arabic*]
\item \label{case1}
 $\mathbb{Z}/r$ acts on $\mathbb{C}^4$ with  weights $(a,-a,1,0)$. We show that $r=2$ in this case. 
By the classification of terminal singularities, the equation $\phi$ is $\mathbb{Z}/r$-invariant  {(and not only $\mathbb{Z}/r$-semi-invariant)}. 
Assume that $r\geq 3$, and let $\xi_r$ be a primitive root of unity of degree $r$. 
Let $T_2\subset \mathbb{C}^4$ be the subspace spanned by the coordinates on which $\mathbb{Z}/r$ acts with weights $a$ and $-a$. By the classification of terminal singularities, we see that $\phi_2|_{T_2}$ is non-degenerate. Since $\phi$ depends on $x_i^2$, it follows that $\phi_2$ contains the terms $x_1^2$ and $x_2^2$. Up to a rescaling of coordinates, we may assume that $\phi_2$ contains $x_1^2+x_2^2$. 

If follows that $\xi_r^{2a} x_1^2 + \xi_r^{-2a} x_2^2$ should be equal to $x_1^2 + x_2^2$, which is possible only if $\xi_r^{2a}=1$. Since $a$ is coprime with $r$, this implies that $r=2$ in this case.
We obtain 
\begin{equation}
\label{eq-terminal-point-abelian-3}
\phi = x_1^2 + x_2^2 + f(x_3,x_4),  \quad r=2,  \quad w=(1,1,1,0;0), 
\end{equation}
where $\mult_0 (f)\geq 2$, $w$ is the weight vector of the action of $\mathbb{Z}/2$ on $\mathbb{C}^4$, and  {its last element is the weight of~$\phi$. It is equal to $0$, indicating that $\phi$ is $\mathbb{Z}/2$-invariant. }

\item\label{case2}
$\mathbb{Z}/r$ acts on $\mathbb{C}^4$ with  weights $(1,3,1,2)$, and $r=4$. The only possibility in this case is
\begin{equation}
\label{eq-terminal-point-abelian-1}
\phi = x_1^2 + x_2^2 + f(x_3,x_4),  \quad r=4,  \quad  w=(1,3,1,2;2), 
\end{equation}
where $\mult_0 (f)\geq 3$ and $w$ is the weight vector of the action of $\mathbb{Z}/4$.

\item\label{case3} $\mathbb{Z}/r$ acts on $\mathbb{C}^4$ with weights $(0,1,1,1)$, and $r=2$. 
We obtain 
\begin{equation}
\label{eq-terminal-point-abelian-2}
\phi = x_1^2 + x_2^2 + f(x_3,x_4),  \quad r=2,  \quad w=(0,1,1,1;0), 
\end{equation}
where $\mult_0 (f)\geq 4$ and $w$ is the weight vector of the action of $\mathbb{Z}/2$.
\end{enumerate}

Let $g$ generate the deck group $\mathbb{Z}/r\subset\widetilde G$. We claim that $\langle g\rangle$ is maximal cyclic in $\widetilde G_2$. Otherwise there exists $h\in\widetilde G_2$ with $h^2=g$. The element $h$ acts on the coordinates $x_1$ and $x_2$ as $x_1\mapsto \xi x_1$, $x_2\mapsto \zeta x_2$, where $\xi$ and $\zeta$ are roots of unity.

In Case~\ref{case2}, we have $\xi^2=\sqrt{-1}$, $\zeta^2=-\sqrt{-1}$, $\xi^2=\zeta^2$. Here the last equality holds since the equation $\phi$ is $\widetilde{G}$-semi-invariant. This gives a contradiction. 
In Case~\ref{case3}, we have $\xi^2=1$, $\zeta^2=-1$, $\xi^2=\zeta^2$. This gives a contradiction. 
In Case~\ref{case1}, $h$ should act via multiplication by $\pm \sqrt{-1}$ on the coordinates $x_1,x_2,x_3$, and by multiplication by $\pm1$ on~$x_4$. 
As in the proof of Proposition~\ref{prop-action-on-gorenstein-points}, one can show that $\phi$ contains a monomial $x_4^{2l}$ with non-zero coefficient. It follows that $\phi$ cannot be $h$-semi-invariant. This gives a contradiction.

Therefore, $\mathbb{Z}/r$ is a maximal cyclic subgroup in $\widetilde{G}_2$ in all three cases considered above. According to Lemma~\ref{lem-decreasing-rank}, we have 
\[
\mathfrak{r}(G_2)=\mathfrak{r}\left(\widetilde{G}_2/\langle g\rangle\right)<\mathfrak{r}\left(\widetilde{G}_2\right)=4.
\] 
Using Proposition~\ref{prop-action-on-gorenstein-points}, we obtain $\mathfrak{r}(G_{\neq 2})\leq \mathfrak{r}(\widetilde{G}_{\neq 2})\leq 2$. We conclude that $\mathfrak{r}(G)\leq 3$, and the claim follows.
\end{proof}

\begin{proposition}
\label{prop-terminal-point-not-abelian}
Let $P\in X$ be a $3$-dimensional terminal singularity germ of Gorenstein index $r>1$, and let~$G$ be a finite abelian group such that $G\subset \Aut(P\in X)$. 
Assume that the lifting $\widetilde{G}$ is not abelian. Then
$\mathfrak{r}(G)\leq 3$.
\end{proposition}

\begin{proof}
By Proposition~\ref{prop-action-on-non-gorenstein-points}, we have $r=2$. 
We may assume that $\mathfrak{r}(G_2)=4$. Indeed, otherwise Proposition~\ref{prop-p-action-on-terminal-point} implies that $\mathfrak{r}(G_p)\leq 3$ for $p\neq 2$, so we have $\mathfrak{r}(G)\leq 3$ and we are done. 

Let $\widetilde{G}_2$ be the Sylow $2$-subgroup of $\widetilde{G}$. 
Note that $\widetilde{G}_2$ surjects to $G_2$ under the natural projection map $\widetilde{G}\to G$ as in the exact sequence~\eqref{eq-lift}. 
We claim that $\widetilde{G}_2$ is abelian. 
In the notation of Theorem~\ref{thm-terminal-sing}, the subspace $T=\mathbb{C}^3=\{x_4=0\}$ is $\widetilde{G}$-invariant. 
We show that there exists a $1$-dimensional subrepresentation of $\widetilde{G}_2$ in $T$. Assume that the representation of $\widetilde{G}_2$ on $T=\mathbb{C}^3=\{x_4=0\}$ is irreducible. Then $3=\dim T$ divides $|\widetilde{G}_2|$, which gives a contradiction. Hence the representation of $\widetilde{G}_2$ on $T$ is reducible. In particular, there exists a $1$-dimensional subrepresentation $T_1\subset T$. 

Observe that the derived subgroup $[\widetilde{G},\widetilde{G}]$ is contained in $\mathbb{Z}/r=\mathbb{Z}/2$. The latter group acts on $T$, and in particular on $T_1\subset T$, faithfully.   
If $\widetilde{G}_2$ is not abelian, then $[\widetilde{G}_2,\widetilde{G}_2]=\mathbb{Z}/2$.  
On the other hand, the kernel of the homomorphism $\widetilde{G}_2\to\GL(T_1)=\mathbb{C}^\times$ contains $[\widetilde{G}_2,\widetilde{G}_2]$, which gives a contradiction. This shows that $\widetilde{G}_2$ is abelian. Proposition~\ref{prop-terminal-point-abelian} applied to $G_2$ and its lifting $\widetilde{G}_2$ shows that $\mathfrak{r}(G_2)\leq 3$. 

For every $p>2$, Proposition~\ref{prop-p-action-on-terminal-point} applied to $G_p$ gives $\mathfrak r(G_p)\leq3$. Together with $\mathfrak r(G_2)\leq3$, this yields $\mathfrak r(G)\leq3$.
\end{proof}

Now, Theorem~\ref{thm-action-on-terminal-point} follows from Propositions~\ref{prop-action-on-gorenstein-points},~\ref{prop-terminal-point-abelian}, and~\ref{prop-terminal-point-not-abelian}. 

\section{Irreducible anti-canonical element}
\label{sec-sing-of-pair}
\label{sec-geom-anti-can-notation}
\label{sect-special-case}

In this section, we start working with anti-canonical elements on a $G\mathbb{Q}$-Fano threefold~$X$, where $G$ is a finite abelian group. 
Assume that $|-K_X|\neq~
\emptyset$. The group $G$ acts on
$\mathrm{H}^0(X, \OOO(-K_X))$. Since $G$ is abelian, the corresponding representation of $G$ splits, so we can pick a $G$-invariant element
$S\in |-K_X|$. 

The main goal in this section is to show that $G$ is either of product type or of K3 type under the assumption that $S$ is irreducible. 
First, we analyze the singularities of the pair $(X, S)$. 
In Proposition~\ref{the-pair-is-lc}, we prove that either the pair $(X, S)$ is lc, or $G$ is of product type. In Proposition~\ref{the-pair-is-plt}, we prove that if $(X, S)$ is plt, then $G$ is of K3 type. 

\begin{proposition}[\textit{cf.} {\cite[2.9]{Pr11}}]
\label{the-pair-is-lc}
In the above notation,
if the pair $(X, S)$ is not lc, then there exists either a
$G$-invariant point on $X$, or a $G$-invariant smooth rational curve on $X$, or a $G$-invariant rational surface on $X$. 
As a consequence, $G$ is of product type. 
\end{proposition}

\begin{proof}
If the pair $(X, S)$ is not lc, then pick a maximal $\lambda<1$ such
that the pair $(X, \lambda S)$ is lc.
Consider a minimal lc-center~$Z$ of $(X, \lambda S)$.  
Take a general $G$-invariant very ample divisor $H$ on $X$ passing through $Z$. 
Then the pair $(X, \lambda S + \epsilon H)$ is not lc at $Z$, 
and it is a log Fano pair for $0<\epsilon\ll 1$. 
Put $D=\mu(\lambda S + \epsilon H)$, where $0<\mu<1$ is the log canonical threshold of $(X, \lambda S + \epsilon H)$. 
Then the pair $(X, D)$ is strictly lc at~$Z$, and~$Z$ is a minimal lc-center of $(X, D)$.
For a sufficiently general choice of $H$, the subvariety $Z$ is the unique lc-center of $(X,D)$ in its neighborhood.
Since the pair $(X, D)$ is $G$-invariant, $g(Z)$ is also a minimal lc-center of $(X, D)$ for any $g\in G$, and  $g(Z)$  is a unique lc-center of $(X, D)$ in its neighborhood.
On the other hand, by Koll\'ar--Shokurov connectedness, the non-klt locus of $(X, D)$ is connected. Hence, $Z$ is $G$-invariant. 

If $\dim Z = 0$, so $Z$ is a point, we apply Corollary~\ref{cor-fixed-point} to conclude that $G$ is of product type. 
If $\dim Z=1$, then $Z$ is a smooth rational curve, and we apply Corollary~\ref{cor-inv-curve-or-surface} to conclude that $G$ is of product type. 
If $\dim Z=2$, then by adjunction $Z$ is a klt log del Pezzo surface. In particular, it is rational. 
Again, we apply Corollary~\ref{cor-inv-curve-or-surface} to conclude that $G$ is of product type. 
\end{proof}

Recall that a K3 surface is a normal projective surface
$S$ with at worst canonical (that is, $1$-lc) singularities such that
$\mathrm{H}^1(S, \OOO_S)=0$ and $K_S\sim 0$. 

\begin{proposition}
\label{the-pair-is-plt}
In the above notation,
if the pair $(X, S)$ is plt, then $G$ is of K3 type.
\end{proposition}

\begin{proof}
Since $S\sim -K_X$ is ample, its support is connected. Since the pair $(X, S)$ is plt, it follows that $S$ is irreducible and reduced.  Moreover, $S$ is normal, and by inversion of adjunction the pair $(S, 0)$ is klt. Since $K_S$ is linearly trivial, from the exact sequence
\[
0 \lra \OOO_X(K_X) \lra \OOO_X \lra \OOO_S \lra 0, 
\]
it follows that $h^1(\OOO_S) = 0$. 
We conclude that $S$ is a K3 surface with at worst canonical singularities. The claim follows.
\end{proof}

Until the end of this section, we work in the following setting. 

\begin{setting}
Let $X$ be a $G\mathbb{Q}$-Fano threefold, where $G$ is a finite abelian group, and let $S\in|-K_X|$ be a $G$-invariant element such that the pair $(X, S)$ is lc.
\end{setting}

Let us explain our plan for the rest of the section. 
We assume that $S$ is irreducible. 
We prove that either~$G$ is of product type, or we are in one of the following situations:  
\begin{itemize}
\item
The surface $S$ is normal (see Proposition~\ref{prop-non-normal-irreducible-surface}).
\item
The surface $S$ has two simple elliptic singularities $P_1$ and $P_2$ interchanged by $G$ (see Proposition~\ref{the-pair-is-strictly-lc-irred}).
\item
Let $H$ be the index $2$ subgroup of $G$ that stabilizes $P_i$. Let $E_i$ be the exceptional elliptic curve over $P_i\in S$ that exists on the minimal resolution of $P_i\in S$. Then~$E_i$ has an $H$-orbit of length at most $7$ (see Propositions~\ref{prop-elliptic-sing-short-orbit} and~\ref{prop-elliptic-sing-short-orbit-2}).
\item
Let $H_E$ be the image of $H$ in the automorphism group of $E_i$. Using the previous results, we obtain that either $\mathfrak{r}(H_E)\leq 1$, or $H_E$ is a subgroup of three special groups (see Corollary~\ref{cor-irreducible-HE}).
\end{itemize}
After that, we show that if $\mathfrak{r}(H_E)\leq 1$, then $G$ is of product type (see Corollary~\ref{prop-rank-HE}). 
Finally, we analyze the three remaining cases. This is done in Propositions~\ref{prop-GE-2},~\ref{prop-GE-1}, and~\ref{prop-E2-2} and Corollary~\ref{cor-GE-3}. So we conclude that $G$ is either of product type or of K3 type under the assumption that $S$ is irreducible. We proceed according to this plan.

\begin{proposition}
\label{prop-non-normal-irreducible-surface}
Assume $S$ is irreducible but not normal. Then $G$ is of product type.  
\end{proposition}

\begin{proof}
Since $S$ is Cohen--Macaulay, the non-normal locus of $S$ is a curve $C$. So we have a $G$-invariant curve $C$ such that $S$ is singular along $C$.  According to Lemma~\ref{lem-fixed-curve-surface}, we have $G\subset \Aut(S)$.  Let $S^\nu$ be the normalization of~$S$. By adjunction, $(S^\nu, \Delta)$ is a $2$-dimensional lc log CY pair with a faithful action of $G$ with boundary $\Delta\neq 0$. Then Corollary~\ref{cor-2-dim-cy-pair-pt} implies that $G$ is of product type.
\end{proof}

\begin{proposition}
\label{the-pair-is-strictly-lc-irred}
Assume that the pair $(X, S)$ is lc but not plt. 
Assume, moreover, that~$S$ is
irreducible and normal. Then either $G$ is of product type, or $S$ has two strictly lc singular points interchanged by $G$, which are simple elliptic singularities
\end{proposition}

\begin{proof}
By adjunction, we have $K_S=(K_X+S)|_S\sim 0$, so the surface $S$ is Gorenstein. By inversion of adjunction,~$S$ has strictly lc singularities. 
Since $-K_S\sim 0$ is nef, the non-klt locus of $S$ consists of either one or two points (\textit{cf.} \cite[Theorem 1.2]{B24}). Clearly, they form a $G$-invariant set. 
If $G$ has a fixed point, we apply Corollary~\ref{cor-fixed-point} and we are done. So we may assume that $G$ interchanges two strictly lc singularities $P_1$ and $P_2$ on $S$. 

Since $P_1$ and $P_2$ are strictly lc singular points on $S$, by Section~\ref{subsec-simple-ell-cusp}, the points $P_i$ are either simple elliptic or cusp singularities. 
We use the notation of Lemma~\ref{lem-2-dim-CY-pairs}. If the $P_i$ are cusp singularities, then all the exceptional curves in the minimal resolution $f\colon (S',\Delta')\to S$ are rational. In particular, all the components of $\Delta'$ are rational curves.
If $\dim Z=1$, then $\Delta''\cdot F=2$, so a rational component of $\Delta''$ dominates $Z$, and hence $Z\simeq\mathbb P^1$; if $\dim Z=0$, then $S''$ is a del Pezzo surface. Thus $S''$ is rational. Therefore, $S$ is rational as well. Hence Corollary~\ref{cor-inv-curve-or-surface} implies that $G$ is of product type. If this is not the case, we conclude that the $P_i$ are simple elliptic singularities.
\end{proof}

\begin{proposition}
\label{prop-elliptic-small-embedding-dimension}
Let $P\in X$ be a threefold terminal singularity germ of Gorenstein index $r\geq 1$. Assume that $P\in S$ is a simple elliptic singularity, where $S\sim -K_X$ is a normal irreducible surface. Let 
\[
\pi\colon \widetilde{P}\in\widetilde{X}\lra X\ni P
\] 
be the canonical cyclic covering, and let $\widetilde{S}=\pi^{-1}(S)$. 
Then $\widetilde{P}\in\widetilde{S}$ also is a simple elliptic singularity.
\end{proposition}

\begin{proof}
  We may assume that $P$ is not Gorenstein, so $r>1$.
By Section~\ref{subsec-terminal-sing}, the canonical covering $\pi$ is \'etale on $\widetilde{X}\setminus\{\widetilde P\}$, and $\widetilde{P}\in\widetilde{X}$ is Gorenstein. The pair $(X, S)$ is lc but not plt at $P$, and hence $(\widetilde{X}, \widetilde{S})$ is lc but not plt at $\widetilde{P}$. Therefore, $\widetilde{P}\in\widetilde{S}$ is either a simple elliptic or a cusp singularity; \textit{cf.}~Section~\ref{subsec-simple-ell-cusp}. 

We show that $\widetilde{P}\in\widetilde{S}$ is a simple elliptic singularity. Indeed, consider a finite morphism $\pi|_{\widetilde{S}}\colon \widetilde{S}\to S$. Let $\widetilde{S}'\to \widetilde{S}$ and $S'\to S$ be the minimal resolutions. Then $\pi|_{\widetilde{S}}$ induces a rational map $\pi'\colon \widetilde{S}'\dashrightarrow S'$. We can blow up some number of (possibly infinitely near) points on $\widetilde{S}'$ to obtain a generically finite morphism $\pi''\colon \widetilde{S}''\to S'$. We summarize these maps in the following commutative diagram: 
\begin{equation}
\label{eq-diagram-elliptic-resolutions}
\begin{tikzcd}
\widetilde{S}'' \arrow{d}{} \arrow{dr}{\pi''} & \\
\widetilde{S}' \arrow[dashed]{r}{\pi'} \arrow[swap]{d}{} &  S' \arrow{d}{} \\
\widetilde{S} \arrow{r}{\pi}  & S\rlap{.}
\end{tikzcd}
\end{equation}
Note that if $\widetilde{P}\in \widetilde{S}$ is a cusp singularity, then a union of rational curves that are exceptional for $\widetilde{S}''\to \widetilde{S}$ would dominate an elliptic curve on $S'$, which cannot happen. So $\widetilde{P}\in \widetilde{S}$ is a simple elliptic singularity, as claimed. 
\end{proof}

\begin{remark}
\label{rem-map-elliptic-morphism}
We keep the same notation as in the proof of Proposition~\ref{prop-elliptic-small-embedding-dimension}. 
Let $q\colon\widetilde S''\to\widetilde S'$ be the birational morphism in Diagram~\eqref{eq-diagram-elliptic-resolutions}. Every $q$-exceptional curve is rational and maps into the elliptic curve $E$, hence is contracted by $\pi''$. Thus $\pi''$ is constant on the fibers of $q$ and descends to a morphism $\pi'\colon\widetilde S'\to S'$. Let $\widetilde E$ be the exceptional elliptic curve over $\widetilde P$. Since the canonical cyclic cover $\pi\colon\widetilde S\to S$ is crepant and $E$ and $\widetilde E$ are the unique exceptional curves of the two simple elliptic germs, we have $\pi'(\widetilde E)=E$. Thus $\pi'|_{\widetilde E}\colon\widetilde E\to E$ is non-constant and hence finite.
\end{remark}

By Proposition~\ref{the-pair-is-strictly-lc-irred}, we may assume that $S$ has two simple elliptic singularities interchanged by $G$.  Let $H$ be a subgroup of $G$ of index $2$ that stabilizes the points $P_1, P_2$ on~$X$. Put $P=P_1$. The local construction in Proposition~\ref{prop-elliptic-small-embedding-dimension} can be made equivariant with respect to a finite group action. Indeed, once~$H$ acts on the germ $P\in X$, there exists a lifting $\widetilde{H}$ of $H$ that acts on $\widetilde{P}\in \widetilde{X}$ as in the exact sequence \eqref{eq-lift}.  Then the surface $\widetilde{S}$ is $\widetilde{H}$-invariant, and so is Diagram \eqref{eq-diagram-elliptic-resolutions}.  In particular, the morphism $\widetilde{E}\to E$ is $\widetilde{H}$-equivariant, where $\widetilde{H}$ acts on~$E$ via $H$.  According to Proposition~\ref{prop-action-on-non-gorenstein-points}, either the lifting $\widetilde{H}$ of $H$ is abelian, or $r=2$, where $r$ is the Gorenstein index of $P\in X$.

\begin{proposition}
\label{prop-elliptic-sing-short-orbit}
Assume that 
\begin{itemize}
\item
$P\in X$ is a threefold terminal singularity germ with a faithful action of a finite abelian group $H$,
\item
an $H$-invariant element $S\in |-K_X|$ has a simple elliptic singularity at $P$,
\item
the lifting $\widetilde{H}$ of\, $H$ to the canonical covering of\, $P\in X$ is abelian.
\end{itemize}
Let $E$ be the preimage of\, $P\in S$ in the minimal resolution $S'$ of\, $S$, which is automatically~$H$-invariant. 
Then $E$ has an $H$-orbit of length at most $4$. 
\end{proposition}

\begin{proof}
By Remark~\ref{rem-map-elliptic-morphism}, the map $\pi'\colon \widetilde{S}'\to S'$ in Diagram \eqref{eq-diagram-elliptic-resolutions} is a morphism. Since $\widetilde{P}\in \widetilde{S}\subset \widetilde{X}$ and the point $\widetilde{P}\in \widetilde{X}$ is terminal Gorenstein, 
we conclude that the embedding dimension of $\widetilde{P}\in \widetilde{S}$ is at most $4$. By Equation~\eqref{eq-elliptic-curve-embedding-dimension}, it follows that $-\widetilde{E}^2\leq 4$, where $\widetilde{E}\subset \widetilde{S}'$ is an elliptic curve over $\widetilde{P}$ in the minimal resolution $\widetilde{S}'\to \widetilde{S}$. Since by assumption the lifting $\widetilde{H}$ of $H$ is abelian, by Proposition~\ref{prop-small-orbit-on-elliptic-curve}, applied to the image of $\widetilde H$ in $\Aut(\widetilde P\in\widetilde S)$, we conclude that there exists an $\widetilde{H}$-orbit of length at most $4$ on $\widetilde{E}$. Since the induced map $\widetilde{E}\to E$ is $\widetilde{H}$-equivariant, this implies that there exists an $H$-orbit of length at most~$4$ on $E$, and we are done.
\end{proof}

\begin{proposition}
\label{prop-elliptic-sing-short-orbit-2}
Assume that  
\begin{itemize}
\item
$P\in X$ is a threefold terminal singularity germ with a faithful action of a finite abelian group $H$,
\item
an $H$-invariant element $S\in |-K_X|$ has a simple elliptic singularity at $P$,
\item
the lifting $\widetilde{H}$ of\, $H$ to the canonical covering of\, $P\in X$ is not abelian. 
\end{itemize}
Let $E$ be the preimage of\, $P\in S$ in the minimal resolution $S'$ of\, $S$, which is automatically~$H$-invariant.  Then $E$ has an $H$-orbit of length at most $7$.
\end{proposition}

\begin{proof}
According to Proposition~\ref{prop-action-on-non-gorenstein-points}, for the Cartier index of $K_X$ at $P\in X$, we have $r=2$. By the classification of terminal singularities, this implies that $\dim T_P X\leq 7$. By Equation~\eqref{eq-elliptic-curve-embedding-dimension}, it follows that $-E^2\leq 7$. 
By Proposition~\ref{prop-small-orbit-on-elliptic-curve}, applied to the image of $H$ in $\Aut(P\in S)$, there exists an $H$-orbit on $E$ of length at most~$7$. 
\end{proof}

We denote by $H_E$ the image of $H$ in the automorphism group of the elliptic curve $E$.

\begin{corollary}
\label{cor-irreducible-HE}
Either $\mathfrak{r}(H_E)\leq 1$, or $H_E$ is isomorphic to a subgroup of one of the following groups:
\begin{equation}
\label{eq-GE-3-cases}
(\mathbb{Z}/3)^2, \quad(\mathbb{Z}/2)^3,\quad \mathbb{Z}/2\times\mathbb{Z}/4. 
\end{equation}
Moreover, if $H_E$ is isomorphic to one of these three groups, then $H_E$ does not consist of translations on $E$, and there exists an $H_E$-orbit on $E$ of length at most $4$. 
\end{corollary}

\begin{proof}
 According to Propositions~\ref{prop-elliptic-sing-short-orbit} and~\ref{prop-elliptic-sing-short-orbit-2}, there exists an orbit of $H_E$ on~$E$ of length at most~$7$. By Proposition~\ref{prop-elliptic-curve}, the claim follows.
\end{proof}

Let $G_S$ be the image of $G$ in $\Aut(S)$. Then we have the  exact sequence
\begin{equation}
\label{eq-exact-seq1}
0\lra G_N \lra G \lra G_S \lra 0, 
\end{equation}
where $G_N$ is a cyclic group that faithfully acts on the tangent space to the generic point of $S$. 
Since $G_N$ acts trivially on $S$, it fixes $P_1$ and $P_2$, and hence $G_N\subset H$. Put $\overline H=H/G_N\subset G_S$. Then $[G_S:\overline H]=2$.
Recall Diagram \eqref{eq-diag-minimal-model-ruled-surface} from the proof of Lemma~\ref{lem-2-dim-CY-pairs}:
\begin{equation}
\label{eq-one-more-diagram}
\begin{tikzcd}
(S',\Delta')  \arrow{r}{g} \arrow[swap]{d}{f} & (S'', \Delta'') \arrow{d}{\pi} \\
S & Z\rlap{,}
\end{tikzcd}
\end{equation}
where $S''$ is a minimal model of the minimal resolution $S'$ of $S$, and the fibers of $S''\to Z$ are smooth rational curves with $Z$ being an elliptic curve. 
Then this diagram is $G_S$-equivariant and hence $\overline H$-equivariant, and $G_S$ faithfully acts on $S''$. 
Note that $E$ belongs to the boundary $\Delta'$, and $E$ maps isomorphically to a component $E''$ of $\Delta''$ on $S''$. Then $E''$ maps isomorphically to $Z$, so we have $Z\simeq E$. 
There exists an exact sequence
\begin{equation}
\label{eq-exact-seq2}
0\lra G_f \lra G_S \lra G_Z \lra 0, 
\end{equation}
where $G_Z$ is the image of $G_S$ in $\Aut(Z)$, and $G_f$ faithfully acts on the generic fiber of $S''\to Z$.
The isomorphism $E\to Z$ is $\overline H$-equivariant, and the image of $\overline H$ in $G_Z$ is $H_E$.
We have that $G_f\subset \Cr_1(\mathbb{C})$.

\begin{proposition}
\label{prop-G-f}
We have $G_f\subset (\mathbb{Z}/2)^2$.
\end{proposition}

\begin{proof} 
Since $(S'', \Delta'')$ is an lc log CY pair, by adjunction to the generic fiber of $S''\to Z$, we see that $\Delta''$ has exactly two components that surject onto $Z$ and intersect each fiber in one point. In fact, these components are exceptional curves over $P_1$ and $P_2$.
Hence $G_f$ admits an invariant set on the generic fiber of the map $S''\to Z$, which consists of two points. 

Let $C\subset G_f$ be the subgroup fixing both points. Then
$C\subset\mathbb C^*$, so $C$ is cyclic. An element of $G_S$
interchanging the two sections acts on $C$ by inversion; since $G_S$
is abelian, $|C|\leq2$. Moreover, $G_f/C\subset\mathbb Z/2$. After
choosing the two points as $0$ and $\infty$, every element of
$G_f\setminus C$ has the form $t\mapsto a/t$ and is therefore an
involution. Hence
$G_f\subset(\mathbb Z/2)^2$.
\color{black}
\end{proof}

Observe that the exact sequence \eqref{eq-exact-seq2} fits into the following commutative diagram:
\begin{equation}
\begin{tikzcd}
& & 0 \ar[d]  & & \ \\
0 \ar[r] & {\overline H}\cap G_f \ar[r] \ar[d, hook] & {\overline H} \ar[d, hook] \ar[r] & H_E \ar[d, hook] \ar[r] & 0 \ \\
0 \ar[r] & G_f \ar[r] & G_S \ar[r] \ar[d] & G_Z \ar[r] & 0\\
& & \mathbb{Z}/2 \ar[d] & & \ \\
& & 0\rlap{.} & & \ 
\end{tikzcd}
\end{equation}
Two cases are possible: 
\begin{enumerate}
\item
either $H_E=G_Z$ and $G_f/({\overline H}\cap G_f)=\mathbb{Z}/2$,
\item
or $G_Z/H_E=\mathbb{Z}/2$ and $G_f/({\overline H}\cap G_f)=0$.
\end{enumerate}

\begin{proposition}
\label{prop-rank-GZ-1}
If\, $\mathfrak{r}(G_Z)\leq 1$, then $G$ is of product type. 
\end{proposition}

\begin{proof}
If $\mathfrak{r}(G_Z)\leq 1$, then applying Propositions~\ref{prop-abstact-extension-cr2} and~\ref{prop-abstact-extension} to the extensions \eqref{eq-exact-seq2} and \eqref{eq-exact-seq1}, respectively, we see that $G$ is of product type. 
\end{proof}

\begin{corollary}
\label{prop-rank-HE}
If\, $\mathfrak{r}(H_E)\leq 1$, then $G$ is of product type. 
\end{corollary}

\begin{proof}
Assume that $\mathfrak{r}(H_E)\leq 1$, so $H_E=\mathbb{Z}/m$ for some $m\geq 1$. If $G_Z$ is cyclic as well, we are done by Proposition~\ref{prop-rank-GZ-1}.  
Otherwise, we have $G_Z=\mathbb{Z}/m \times \mathbb{Z}/2$. In this case, $G_f/(\overline H\cap G_f)=0$, so $G_f\subset\overline H$. Since $H$ stabilizes $P_1$ and $P_2$, the group $\overline H$ preserves the two corresponding sections of $S''\to Z$.
Hence $G_f$ does not contain an element that interchanges two sections of $S''\to Z$ that are components of the boundary divisor $\Delta''$ as in Diagram \eqref{eq-one-more-diagram}. Therefore, $G_f\subset \mathbb{Z}/2$. Thus, $G_S\subset \Cr_2(\mathbb{C})$, and hence $G$ is of product type. 
\end{proof}

\begin{proposition}
\label{prop-reduction-for-GZ}
Either $G$ is of product type, or $G_Z$ is isomorphic to a subgroup of one of the groups \eqref{eq-GE-3-cases}.
\end{proposition}

\begin{proof}
By Corollary~\ref{cor-irreducible-HE}, we know that either $H_E$ is a subgroup of one of the groups \eqref{eq-GE-3-cases}, or $\mathfrak{r}(H_E)\leq 1$.  
By Corollary~\ref{prop-rank-HE}, we may assume that $H_E$ is a subgroup of one of the groups \eqref{eq-GE-3-cases}. Then by Proposition~\ref{prop-elliptic-curve}, there exists an $H_E$-orbit on $E$ of length at most~$4$. 
Since by the above either $G_Z=H_E$ or $G_Z/H_E=\mathbb{Z}/2$, we get that there exists a $G_Z$-orbit on $Z$ of length at most~$8$. Again, applying Proposition~\ref{prop-elliptic-curve}, we see that either $\mathfrak{r}(G_Z)\leq 1$, in which case we are done by Proposition~\ref{prop-rank-GZ-1}, or $G_Z$ is isomorphic to a subgroup of one of the groups \eqref{eq-GE-3-cases}. 
\end{proof}

According to Proposition~\ref{prop-reduction-for-GZ}, we need to consider three possibilities for $G_Z$ as in the list \eqref{eq-GE-3-cases}. 

\begin{proposition}
\label{prop-GE-2}
Assume that $G_Z=(\mathbb{Z}/3)^2$. Then $G$ is of product type. 
\end{proposition}

\begin{proof}
From Proposition~\ref{prop-G-f} and the exact sequence \eqref{eq-exact-seq2}, we have $G_S\subset (\mathbb{Z}/6)^2$. Thus using the exact sequence \eqref{eq-exact-seq1}, we obtain $\mathfrak{r}(G)\leq 3$, and so $G$ is of product type. 
\end{proof}

\begin{proposition}
\label{prop-GE-1}
Assume that $G_Z\subset(\mathbb{Z}/2)^3$. Then either $G$ is of product type, or
\begin{equation}
G=(\mathbb{Z}/2)^3\times \mathbb{Z}/4 \times \mathbb{Z}/2n, \quad H=(\mathbb{Z}/2)^2\times \mathbb{Z}/4\times \mathbb{Z}/2n,  \quad n\geq 1,
\end{equation}
\[
G_Z=(\mathbb{Z}/2)^3,  \quad G_f=(\mathbb{Z}/2)^2,  \quad G_S = (\mathbb{Z}/2)^3\times \mathbb{Z}/4,  \quad G_N = \mathbb{Z}/2n. 
\]
Moreover, the sequence \eqref{eq-exact-seq1} splits.  
The points $P_i$ are terminal Gorenstein of type cA, and $-E^2=2$. 
\end{proposition}

\begin{proof}
We may assume that $G_Z=(\mathbb{Z}/2)^3$. Otherwise, $G_Z\subset \Cr_1(\mathbb{C})$; hence by Proposition~\ref{prop-abstact-extension-cr2} applied to the exact sequence \eqref{eq-exact-seq2}, we have $G_S\subset \Cr_2(\mathbb{C})$, and by Proposition~\ref{prop-abstact-extension} applied to the exact sequence \eqref{eq-exact-seq1}, we conclude that $G$ is of product type. We may assume that $G_f=(\mathbb{Z}/2)^2$. Otherwise, arguing similarly, we see that $G_S\subset \Cr_2(\mathbb{C})$ and $G$ is of product type.

So we have $G_Z=(\mathbb{Z}/2)^3$ and  $G_f=(\mathbb{Z}/2)^2$. Then, as $2$-groups, $G_Z$ has type $[1,1,1]$, $G_f$ has type $[1,1]$, and $(G_N)_2$ is cyclic of order $2^k$ for some $k\geq0$. First, Theorem~\ref{thm-fulton} applied to the exact sequence~\eqref{eq-exact-seq2} gives
\[
[1,1,1]\cdot[1,1]=[2,2,1]+[2,1,1,1]+[1,1,1,1,1].
\]
Thus $G_S$ has one of these three types. If $k=0$, Table~\ref{table-1}, together with the fact that $G_{\neq2}=(G_N)_{\neq2}$ is cyclic, shows that $G$ is of product type. If $k=1$, Theorem~\ref{thm-fulton}, applied to the exact sequence~\eqref{eq-exact-seq1}, gives the following possible types of $G_2$:
\[
[3,2,1],\ [2,2,2],\ [2,2,1,1],\ [3,1,1,1],\ [2,1,1,1,1],\ [1,1,1,1,1,1].
\]
All these types are of product type except possibly the last one. In the last case the $2$-part of the subgroup $H$ has rank~$5$, whereas Theorem~\ref{thm-action-on-terminal-point} gives $\mathfrak r(H)\leq4$, a contradiction. Thus, in the remaining case we may assume $k\geq2$. We compute the possible type of $G_2$ using Theorem~\ref{thm-fulton} applied to the exact sequences \eqref{eq-exact-seq1} and \eqref{eq-exact-seq2}:
\begin{multline*}
[k]\cdot[1,1,1]\cdot[1,1] = [k+2, 1, 1, 1] + 2[k+1, 2, 1, 1] \\+ 2[k+1, 1, 1, 1, 1] + [k, 2, 2, 1] 
+ [k, 2, 1, 1, 1] \\+ [k, 1, 1, 1, 1, 1] + [k+2, 2, 1] + [k+1, 2, 2].
\end{multline*}

Since $\mathfrak{r}(G_{\neq 2})=\mathfrak{r}((G_N)_{\neq 2})\leq 1$, we see that the group $G$ might not be of product type only if $G_2$ has type $[k, 1, 1, 1, 1, 1]$ or $[k, 2, 1, 1, 1]$. 
We exclude the first possibility. In this case, for the index $2$ subgroup $H\subset G$ that stabilizes~$P_i$, we have $\mathfrak{r}(H)\geq 5$. However, by Theorem~\ref{thm-action-on-terminal-point}, we have $\mathfrak{r}(H)\leq 4$. This is a contradiction. Thus, either $G$ is of product type, or $G_2$ has type $[k, 2, 1, 1, 1]$, so
\[
G=\mathbb{Z}/2n\times \mathbb{Z}/4\times (\mathbb{Z}/2)^3, 
\]
where $n\geq 1$. We also justify the asserted splitting. Since $k\geq2$, the cyclic subgroup $(G_N)_2\subset G_2$ has order $2^k$. As $G_2$ has type $[k,2,1,1,1]$, its direct-sum decomposition may be chosen so that $(G_N)_2$ is one of the cyclic factors. Moreover, $G_S$ is a $2$-group, so $G_{\neq2}\subset G_N$ and $G_N=(G_N)_2\times G_{\neq2}$. Therefore $G_N$ is a direct summand of $G$, and the sequence \eqref{eq-exact-seq1} splits. Comparing types in $G_2=(G_N)_2\times G_S$, we also obtain $G_S=\mathbb Z/4\times(\mathbb Z/2)^3$.

Since $\mathfrak{r}(G)=\mathfrak{r}(G_2)= 5$, we have $\mathfrak{r}(H)\geq 4$. 
By Theorem~\ref{thm-action-on-terminal-point}, we see that $\mathfrak{r}(H)=4$ and the points $P_i$ are terminal Gorenstein of type cA. By Equation~\eqref{eq-elliptic-curve-embedding-dimension}, we have $-E^2\leq 4$.  
By Proposition~\ref{prop-elliptic-curve}, since $G_Z=(\mathbb{Z}/2)^3$, there exists a $G_Z$-orbit of length~$4$ on $Z$, and this is the minimal length of a $G_Z$-orbit on $Z$. 
Since $H_E$ has index at most $2$ in $G_Z$ and the isomorphism $E\to Z$ induced by the map $S''\to Z$ as in Diagram \eqref{eq-one-more-diagram} is $\overline H$-equivariant, the minimal length of an $H_E$-orbit on $E$ is either $4$ or $2$. 

Assume that the minimal length of an $H_E$-orbit on $E$ is $4$. 
From Proposition~\ref{prop-small-orbit-on-elliptic-curve}, it follows that $-E^2=4$.
From Equation~\eqref{eq-elliptic-curve-embedding-dimension} and the fact that the points $P_i$ are Gorenstein, we obtain $\dim T_{P_i}S=4$, and hence $T_{P_i}S=T_{P_i}X$. Every element of $G_N$ acts trivially on $S$, and therefore trivially on $T_{P_i}S=T_{P_i}X$. Lemma~\ref{lem-faithful-action} gives $G_N=1$, contradicting the reduction $k\geq2$ above. Thus this case cannot occur unless $G$ is already of product type.
 
Now assume that the minimal length of an $H_E$-orbit on $E$ is $2$.
From Proposition~\ref{prop-small-orbit-on-elliptic-curve}, it follows that $2\leq -E^2$. 
The case $-E^2=4$ is excluded as above. 
We exclude the case $-E^2=3$. Indeed, the proof of Proposition~\ref{prop-small-orbit-on-elliptic-curve} produces an $H_E$-invariant effective divisor of degree~$3$ on $E$. Since $H_E$ is a $2$-group, this divisor contains an $H_E$-fixed point, contradicting Proposition~\ref{prop-elliptic-curve}. Hence we have $-E^2=2$. 
Above we have shown that $\mathfrak{r}(H)=4$. It follows that either $H=(\mathbb{Z}/2)^2\times \mathbb{Z}/4\times \mathbb{Z}/2n$ for $n\geq 1$, or $H=(\mathbb{Z}/2)^3\times \mathbb{Z}/4\times \mathbb{Z}/n$ and $n=1$. In the latter case, $G$ is of product type. 
Since we also know that the points $P_i$ are Gorenstein of type $cA$, the claim follows.  
\end{proof}

\begin{proposition}
\label{prop-E2-2}
Assume that $G_Z\subset\mathbb{Z}/2\times\mathbb{Z}/4$. Then either $G$ is of product type, or
\begin{equation}
G=(\mathbb{Z}/2)^3\times \mathbb{Z}/4 \times \mathbb{Z}/2n, \quad H=(\mathbb{Z}/2)^2\times \mathbb{Z}/4\times \mathbb{Z}/2n,   \quad n\geq 1, 
\end{equation}
\[
G_Z=\mathbb{Z}/2\times \mathbb{Z}/4,  \quad G_f=(\mathbb{Z}/2)^2,  \quad G_S = (\mathbb{Z}/2)^3\times \mathbb{Z}/4, \quad G_N = \mathbb{Z}/2n. 
\]
Moreover, the sequences \eqref{eq-exact-seq1} and  \eqref{eq-exact-seq2} split.  
The points $P_i$ are terminal Gorenstein of type cA, and $-E^2=2$. 
\end{proposition}

\begin{proof}
By Proposition~\ref{prop-rank-GZ-1}, we may assume that $G_Z=(\mathbb{Z}/2)^2$ or $G_Z=\mathbb{Z}/2\times\mathbb{Z}/4$. 
In the former case, we have $G_Z\subset \Cr_1(\mathbb{C})$. By Proposition~\ref{prop-G-f}, we have $G_f\subset \Cr_1(\mathbb{C})$. Applying Proposition~\ref{prop-abstact-extension-cr2} to the exact sequence \eqref{eq-exact-seq2} and Proposition~\ref{prop-abstact-extension} to the exact sequence \eqref{eq-exact-seq1}, we see that $G$ is of product type. 
So we may assume that $G_Z=\mathbb{Z}/2\times\mathbb{Z}/4$. 

 Arguing as at the beginning of the proof of Proposition~\ref{prop-GE-1}, we may assume that $G_f=(\mathbb{Z}/2)^2$. 
Assume that the sequence \eqref{eq-exact-seq2} does not split. Then either 
\[
G_S= \mathbb{Z}/2\times\mathbb{Z}/4\times\mathbb{Z}/4, \quad \text{or} \quad G_S= \mathbb{Z}/2\times\mathbb{Z}/2\times\mathbb{Z}/8,  \quad \text{or} \quad G_S= \mathbb{Z}/4\times\mathbb{Z}/8.
\] 
In all these cases $G_S\subset \Cr_2(\mathbb{C})$; hence applying Proposition~\ref{prop-abstact-extension} to the exact sequence \eqref{eq-exact-seq1}, we see that $G$ is of product type.

Thus we may assume that the exact sequence \eqref{eq-exact-seq2} splits, so $G_S= (\mathbb{Z}/2)^3\times \mathbb{Z}/4$. Using Theorem~\ref{thm-fulton}, we check that if the exact sequence \eqref{eq-exact-seq1} does not split, then $G$ is of product type.  
Now assume that the exact sequence \eqref{eq-exact-seq1} splits, so 
\[
G=\mathbb{Z}/n'\times (\mathbb{Z}/2)^3\times \mathbb{Z}/4
\] 
for some $n'\geq 1$.
Moreover, we may assume that $n'=2n$; otherwise, $G$ is of product type. 

Recall that we denote by $H$ the subgroup of  $G$ that stabilizes the points $P_i$. 
Since $\mathfrak{r}(G)=\mathfrak{r}(G_2)= 5$, we have $\mathfrak{r}(H)\geq 4$. 
By Theorem~\ref{thm-action-on-terminal-point}, we see that $\mathfrak{r}(H)=4$ and the points $P_i$ are terminal Gorenstein of type cA. By Equation~\eqref{eq-elliptic-curve-embedding-dimension}, we have $-E^2\leq 4$.  
By Proposition~\ref{prop-elliptic-curve}, since $G_Z=\mathbb{Z}/2\times \mathbb{Z}/4$, there exists a $G_Z$-orbit of length~$2$ on $Z$, and this is the minimal length of a $G_Z$-orbit on $Z$. 
Since $H_E$ has index at most $2$ in $G_Z$ and the isomorphism $E\to Z$ induced by the map $S''\to Z$ as in Diagram \eqref{eq-one-more-diagram} is $\overline H$-equivariant, the minimal length of an $H_E$-orbit on $E$ is either $2$ or $1$. 
By Corollary~\ref{prop-rank-HE}, we may assume that $H_E$ is not cyclic; hence Proposition~\ref{prop-elliptic-curve} shows that the minimal length of an $H_E$-orbit on $E$ is equal to~$2$.

Above we have shown that $\mathfrak{r}(H)=4$. It follows that either $H=\mathbb{Z}/2n\times(\mathbb{Z}/2)^2\times \mathbb{Z}/4$ for $n\geq 1$, or $H=\mathbb{Z}/n\times (\mathbb{Z}/2)^3\times \mathbb{Z}/4$ and $n=1$, in which case $G$ is of product type. So we have 
\[
G=(\mathbb{Z}/2)^3\times \mathbb{Z}/4 \times \mathbb{Z}/2n, \quad H=(\mathbb{Z}/2)^2\times \mathbb{Z}/4\times \mathbb{Z}/2n
\] 
for some $n\geq 1$. 

If $-E^2= 4$, then the embedding dimension of $P\in S$ is equal to $4$, and we conclude by arguing as in Proposition~\ref{prop-GE-1}. Thus we may assume that $-E^2\leq3$. We exclude the cases $-E^2=1$ and $-E^2=3$.
From Proposition~\ref{prop-small-orbit-on-elliptic-curve}, it follows that $2\leq -E^2$. 
If $-E^2=3$, the $H_E$-invariant effective divisor of degree~$3$ constructed in the proof of Proposition~\ref{prop-small-orbit-on-elliptic-curve} contains an $H_E$-fixed point, since $H_E$ is a $2$-group; this contradicts Proposition~\ref{prop-elliptic-curve}. Hence we have $-E^2=2$. 
Since we also know that the points $P_i$ are Gorenstein of type $cA$, the claim follows.  
\end{proof}

\begin{corollary}
\label{cor-GE-3}
In the case $-E^2=2$, the group $G$ is of product type. 
\end{corollary}

\begin{proof}
  We have to deal with two remaining cases, described in the conclusions of Propositions~\ref{prop-GE-1} and~\ref{prop-E2-2}. In particular, in both cases, $P\in X$ is a Gorenstein point of type cA.  Since in both cases $\mathfrak{r}(H)=4$, we can apply Proposition~\ref{prop-action-on-gorenstein-points}.  Thus, we can choose coordinates in $\mathbb{C}^4$ that contain
  $0=P\in X$ such that the action of $(\mathbb{Z}/2)^4\subset H$ is standard with respect to the coordinates $x_1,\ldots, x_4$ and the local equation $\phi=0$ of $0\in X$ depends on $x_i^2$.  Moreover, for the local equation of $0\in X$ (see Equation~\eqref{eq-gorenstein-rk-2}), we have
\begin{equation}
\label{eq-phi-psi-equation}
\phi = x_1^2 + x_2^2 + x_3^{2k} + x_4^{2l} + \cdots
\end{equation}
for some $k, l\geq 1$. 

Let $\psi=0$ be the equation of the divisor $D$ in~$\mathbb{C}^4$ such that $D|_X = S$. Since $-E^2=2$, by Equation~\eqref{eq-elliptic-curve-embedding-dimension}, we have $\dim T_PS=3$. Therefore, for the linear term $\psi_1$ of $\psi$, we have $\psi_1\neq 0$.  By the classification of simple elliptic singularities, \textit{cf.} Equation~\eqref{eq-simple-elliptic-sing-1}, we have $\rk(\phi_2|_{\psi=0})=1$.  Since $\rk(\phi_2|_{\psi=0})=1$, we see that $\psi_1$ contains either $x_1$ or $x_2$ with a non-zero coefficient.

Thus we may assume that $\psi_1$ contains $x_1$.  Note that $\psi$ (and so $\psi_1$) can be chosen $H$-semi-invariant, and also $(\mathbb{Z}/2)^4$-semi-invariant.  We conclude that~$\psi$ should be divisible by $x_1$. But since $S$ is irreducible, this implies that, up to a rescaling of the coordinates, we have $\psi = x_1$.

By $(\phi|_{\psi=0})_4$ we denote the homogeneous term of degree $4$ of $\phi|_{\psi=0}$. 
By Equation~\eqref{eq-simple-elliptic-sing-1}, we have that $(\phi|_{\psi=0})_4$ has four distinct roots. It follows that $x_3^4$ and $x_4^4$ belong to $(\phi|_{\psi=0})_4$ with non-zero coefficients. So after a rescaling of coordinates, we have
\[
\phi = x_1^2 + x_2^2 + x_3^4 + x_4^4 + \cdots.
\]

By Propositions~\ref{prop-GE-1} and~\ref{prop-E2-2}, it follows that $G_N=H_N$, where $G_N$ is as in the sequence~\eqref{eq-exact-seq1} and $H_N$ is a subgroup of $H$ that acts trivially on~$S$.  
Let $g$ be the generator of $G_N$. Since $g$ acts on $S$ trivially and $S$ is given by the equation $\phi|_{x_1=0}=0$ in $\mathbb{C}^3$, by Lemma~\ref{lem-faithful-action}, we have that $g$ acts on $\mathbb{C}^4$ as follows: $x_1\mapsto \xi_{2n} x_1$ and $x_i\mapsto x_i$ for $2\leq i \leq 4$, where $\xi_{2n}$ is a root of unity of order $2n$. However, since $\phi$ is $H$-semi-invariant, and hence $G_N$-semi-invariant, this implies that $n\leq 1$. 
It follows that $G=\mathbb Z/4\times(\mathbb Z/2)^4$, which is of product type. The claim follows.
\end{proof}

We conclude this section by presenting some examples. 

\begin{example}
Consider an ordinary double point germ $P\in X$ given by the equation
\[
\phi = x_1^2+x_2^2+x_3^2+x_4^2 = 0,
\]
and let $S\subset X$ be given by the equation
\[
\psi = a_1x_1^2+a_2x_2^2+a_3x_3^2+a_4x_4^2 = 0,
\]
where the $a_i\in \mathbb{C}$ are general. Then there exists a diagonal faithful action of the group $H=(\mathbb{Z}/2)^3\times\mathbb{Z}/2n$ on $P\in X$. Note that an extension of $H$ by $\mathbb{Z}/2$ is of product type. 
\end{example}

\begin{example}
Consider a germ of a terminal singularity $P\in X$ given by the equation
\[
\phi = x_1^2+x_2^2+x_3^4+x_4^4 = 0.
\]
Then there exists a diagonal faithful action of the group $H=(\mathbb{Z}/2)^2\times\mathbb{Z}/4\times \mathbb{Z}/4n$ on $P\in X$ for $n\geq 1$. Note that  $S=\{x_1=0\}$ has a simple elliptic singularity at $P$ with $-E^2=2$. However, in this case the image of $\mathbb{Z}/4n$ in $G_S$ is contained in $G_f$, where $G_f$ is as in the exact sequence~\eqref{eq-exact-seq2}, so this case is not realized because by Proposition~\ref{prop-G-f} we have $G_f\subset (\mathbb{Z}/2)^2$.
\end{example}

\section{Reducible anti-canonical element: Dlt case}
\label{sec-dual-complex}
In this section, we treat the case of a $G$-invariant reducible element $S\in |-K_X|$ such that the pair $(X, S)$ is divisorially log terminal. We study the induced action of $G$ on $S$. The key ingredient here is the notion of dual complex. 
We recall some results on dual complexes of Calabi--Yau pairs; see Section~\ref{subsec-dual-complex} for the main definitions. 

\begin{lem}[\textit{cf.} {\cite[Claim~32.1(1)]{KX16}}]
\label{lem-strata-are-rational}
Let $(X, S=\sum S_i)$ be a dlt log CY pair.
Assume that $\dim \mathcal{D}(S)=\dim X-1$. Then all the strata of\, $S$ are rationally connected. In particular, if $\dim X=3$, then all the strata of\, $S$ are rational. 
\end{lem}

\begin{lem}[\textit{cf.} {\cite[Lemma 25]{KX16}}]
\label{lem-dual-cohomology}
Let $(X, S = \sum S_i)$ be an snc pair. Then for any $i\geq0$
there is an inclusion
\[
\mathrm{H}^i (\mathcal{D}(S),\mathbb{C}) \xhookrightarrow{\hphantom{0pt}} \mathrm{H}^i (S, \OOO_S).
\]
\end{lem}

\begin{corollary}
\label{lem-dual-cohomology-dlt}
Let $(X, S = \sum S_i)$ be a dlt pair. Then for any $i\geq 0$ there is an inclusion
\[
\mathrm{H}^i (\mathcal{D}(S),\mathbb{C}) \xhookrightarrow{\hphantom{0pt}} \mathrm{H}^i (S, \OOO_S).
\]
\end{corollary}
\begin{proof}

Let $f\colon (\widetilde{X}, \widetilde{S})\to (X, S)$ be a log resolution of the dlt pair $(X, S)$, so that 
\[
f^*(K_X+S)=K_{\widetilde{X}}+\widetilde{S}.
\] 
Since $(X, S)$ is dlt, we can choose the log resolution $f$ in such a way that the log discrepancies of all the $f$-exceptional divisors are positive. Put $S'=\widetilde{S}^{=1}$. 
The induced map $\pi=f|_{S'}\colon S'\to S$ induces a bijection between the components of $S'$ and the components of $S$ and is birational on each component of $S'$. 
The argument in \cite[proof of Theorem~59(1)]{KX16} gives isomorphisms
$
\mathrm{H}^i(S,\OOO_S)\simeq\mathrm{H}^i(S',\OOO_{S'})
$
for every $i$.
We obtain the following commutative diagram:
\begin{equation}
\begin{tikzcd}
\mathrm{H}^i (\mathcal{D}(S'),\mathbb{C}) \arrow[hookrightarrow]{r}{g} \arrow[swap]{d}{p} & \mathrm{H}^i (S', \OOO_{S'}) \arrow{d}{h}\\
\mathrm{H}^i (\mathcal{D}(S),\mathbb{C})&  \mathrm{H}^i (S, \OOO_S)\rlap{,} 
\end{tikzcd}
\end{equation}
where $h$ is an isomorphism. Note that $p$ is an isomorphism as well since
$f$ induces a homeomorphism between $\mathcal{D}(S)$ and $\mathcal{D}(S')$. 
By Lemma~\ref{lem-dual-cohomology}, the map $g$ is an inclusion. We conclude that the induced map $h\circ g \circ p^{-1}$ is an inclusion. The claim follows. 
\end{proof}

\begin{corollary} 
\label{not-a-circle}
Let $X$ be a Fano threefold with terminal singularities.
Assume that the pair $(X, S)$ is lc, where $S=\sum S_i$ and $S\sim -K_X$. Then the dual complex $\mathcal{D}(X, S)$ is homeomorphic to one of the following topological spaces:
\begin{enumerate}
\item
a point,
\item
an interval $I$,
\item
a sphere $\mathbb{S}^2$.
\end{enumerate}
\end{corollary}

\begin{proof}
Let $f\colon (\widetilde{X}, \widetilde{S})\to (X, S)$ be a $\mathbb{Q}$-factorial 
dlt modification of $(X, S)$, so we have
\[
K_{\widetilde{X}}+\widetilde{S}= f^*(K_X+S)\sim 0.
\]  
Put $S'=\widetilde{S}=\widetilde{S}^{=1}$.  Denote by $\mathcal{D}(S')$ the dual complex of $S'$, which is equal to the dual complex $\mathcal{D}(X, S)$ by definition.  Assume that $\dim \mathcal{D}(S')=0$. Since $S$ is ample, its support is connected, and so is the support of~$S'$; hence $\mathcal{D}(S')$ is a point.

Assume that $\dim \mathcal{D}(S')=1$. By the main result of \cite{KX16}, we have that $\mathcal{D}(S')$ is homeomorphic to the quotient of a circle by a finite group. Thus, $\mathcal{D}(S')$ is homeomorphic either to a circle or to an interval. Assume that~$\mathcal{D}(S')$ is homeomorphic to a circle $\mathbb{S}^1$. 
Consider an exact sequence

\[
0 \lra \OOO_{\widetilde X}(-S') \lra \OOO_{\widetilde X} \lra \OOO_{S'} \lra 0.
\]
Since $\widetilde X$ is rationally connected, we have $h^i (\widetilde X, \OOO_{\widetilde X}) = 0$
for $i>0$. 
Using Serre duality, we obtain 
\[
h^2 \left(\widetilde X, \OOO_{\widetilde X}(-S')\right)=h^2 \left(\widetilde X, \OOO_{\widetilde X}(K_{\widetilde X})\right)=h^1 \left(\widetilde X, \OOO_{\widetilde X}\right)=0.
\] 
Therefore, $h^1 (S', \OOO_{S'}) = 0$. This contradicts Corollary~\ref{lem-dual-cohomology-dlt}. We conclude that $\mathcal{D}(S')$ is homeomorphic to an interval in this case. 
\color{black}

Finally, assume that $\dim \mathcal{D}(S')=2$. Then $(\widetilde X,S')$ has maximal intersection, and \cite[Section~33]{KX16} shows that $\mathcal{D}(S')$ is homeomorphic to $\mathbb S^2$.
\end{proof}

\begin{lem}
\label{lem-transitive-action}
Let $G$ be a finite abelian group. 
Assume that 
\begin{itemize}
\item
$X$ is a $G\mathbb{Q}$-Fano threefold and the action of\, $G$ on $X$ is faithful,
\item 
$S\in |-K_X|$ is a $G$-invariant element,
\item
the pair $(X, S)$ is lc, and 
\item
$S=\sum S_i$ is reducible. 
\end{itemize}
Then either $G$ is of product type, or the action of\, $G$ on the set $\{S_i\}$ of irreducible components of\, $S$ is transitive.
\end{lem}

\begin{proof}
Since $X$ is a $G\mathbb{Q}$-Fano variety, we have $\rho^G(X)=1$. If the action of $G$ on the set of irreducible components of $S=\sum S_i$ is not transitive, then $S=S'+S''$, where $S'$ and $S''$ are $G$-invariant Weil $\mathbb{Q}$-Cartier divisors. In particular, $S'$ and $S''$ are proportional in $\Pic(X)$. It follows that the pair $(X, S')$ is an lc log Fano pair. 
The following argument is similar to the one in the proof of Proposition~\ref{the-pair-is-lc}.

Consider a minimal lc-center $Z$ of $(X, S')$. After perturbing the pair by a $G$-invariant divisor with small coefficients, we may assume that $Z$ is an isolated lc-center of an lc log Fano pair $(X,\Delta)$. Since every translate $g(Z)$ is also an isolated lc-center, Koll\'ar--Shokurov connectedness implies that $Z$ is $G$-invariant. By subadjunction, $Z$ is either a point, a smooth rational curve, or a klt log del Pezzo surface (in particular, it is rational).  In the first case, apply Corollary~\ref{cor-fixed-point} to conclude that $G$ is of product type.  In the last two cases, apply Corollary~\ref{cor-inv-curve-or-surface} to conclude that $G$ is of product type.  The claim is proven.
\end{proof}

\begin{remark}
\label{rem-fixed-point-on-dual-complex}
Let $(X, S=\sum S_i)$ be a dlt log CY pair with $\dim X=3$. 
Assume that a finite abelian group $G$ acts faithfully on $(X, S)$ in such a way that the induced action on a dual complex $\mathcal{D}(S)$ has a fixed point. Then there exists a $G$-invariant stratum of $S$. In particular, if $X$ has terminal singularities and $\dim \mathcal{D}(S)=2$, then by Lemma~\ref{lem-strata-are-rational}, all the strata of $S$ are rational, and hence $G$ is of product type by Corollary~\ref{cor-fixed-point} or~\ref{cor-inv-curve-or-surface}.
\end{remark}

\begin{proposition}
\label{prop-reducible}
Let $G$ be a finite abelian group. 
Assume that 
\begin{itemize}
\item
$X$ is a $G\mathbb{Q}$-Fano threefold and the action of\, $G$ on $X$ is faithful,
\item 
$S\in |-K_X|$ is a reducible $G$-invariant element,
\item
the pair $(X, S)$ is dlt, and 
\item
the dual complex $\mathcal{D}(S)$ is homeomorphic to an interval. 
\end{itemize}
Then $G$ is of product type.
\end{proposition}

\begin{proof}
By assumption, $S=\sum S_i$ is reducible, and hence $S$ is singular along a curve. 
By Lemma~\ref{lem-fixed-curve-surface}, the action of $G$ on $S$ is faithful. 
By Lemma~\ref{lem-transitive-action}, we may assume the action of $G$ on the set of components of $S$ is transitive. 
Thus, $S$ has only two components: $S=S_1+S_2$. We claim that both $S_i$ are rational. 

Indeed, the two components are isomorphic because they are
interchanged by $G$. Suppose that they are not rational. By the proof
of Lemma~\ref{lem-2-dim-CY-pairs}, each $S_i$ is then birationally
ruled over an elliptic curve, and therefore
$
h^1(S_i,\OOO_{S_i})=1.
$
Put $C=S_1\cap S_2$. The curve $C$ is connected, and adjunction gives
$p_a(C)\leq1$. As in the proof of
Proposition~\ref{the-pair-is-plt}, the exact sequence
\[
0\lra\OOO_X(K_X)\lra\OOO_X\lra\OOO_S\lra0
\]
gives $H^1(S,\OOO_S)=0$. On the other hand, from
\begin{equation}
\label{eq-normalization-conductor}
0\lra\OOO_S\lra\OOO_{S_1}\oplus\OOO_{S_2}
\lra\OOO_C\lra0
\end{equation}
we obtain an injection
\[
H^1(S_1,\OOO_{S_1})\oplus H^1(S_2,\OOO_{S_2})
\lhook\joinrel\longrightarrow H^1(C,\OOO_C),
\]
which is impossible because the left-hand side has dimension $2$, whereas the
right-hand side has dimension at most $1$. Thus both $S_i$ are rational.
\color{black}
By Lemma~\ref{lem-diagonal-in-product}, we have that $G$ is an extension of a finite abelian subgroup of $\Cr_2(\mathbb{C})$ by $\mathbb{Z}/2$. By Proposition~\ref{prop-abstact-extension}, we conclude that $G$ is of product type. 
\end{proof}

\begin{proposition}
\label{prop-reducible-2}
Let $G$ be a finite abelian group. 
Assume that 
\begin{itemize}
\item
$X$ is a $G\mathbb{Q}$-Fano threefold and the action of\, $G$ on $X$ is faithful,
\item 
$S\in |-K_X|$ is a $G$-invariant element,
\item 
the pair $(X, S)$ is dlt, and
\item
the dual complex $\mathcal{D}(S)$ is homeomorphic to a sphere $\mathbb{S}^2$. 
\end{itemize}
Then $G$ is of product type.
\end{proposition}

\begin{proof}
{
The action on $S$ is faithful by Lemma~\ref{lem-fixed-curve-surface},
and by Lemma~\ref{lem-transitive-action} we may assume that $G$ acts
transitively on the set $\Omega$ of components of $S$. Since
$\dim\mathcal D(S)=2$, all components and strata are rational by
Lemma~\ref{lem-strata-are-rational}.

If the induced action of $G$ on
$\mathcal D(S)\simeq\mathbb S^2$ has a fixed point, we are done
by Remark~\ref{rem-fixed-point-on-dual-complex}. Otherwise,
Lemma~\ref{lem-action-on-spheres} gives a $G$-orbit $\{x_1,x_2\}$.
Let $\sigma_i$ be the unique open cell containing $x_i$, and let $V_i$
be the set of components of $S$ corresponding to the vertices of
$\overline\sigma_i$. Both
$V_1\cap V_2$ and $V_1\cup V_2$ are $G$-invariant subsets of $\Omega$.
Transitivity implies that $V_1\cup V_2=\Omega$ and that either
$V_1\cap V_2=\Omega$ or $V_1\cap V_2=\varnothing$. In the first case
$|\Omega|\leq3$, since a cell of the two-dimensional dual complex has
at most three vertices. In the second case
\[
|\Omega|=2|V_1|\in\{2,4,6\}.
\]

Put $G_\Omega=\operatorname{im}(G\to\operatorname{Sym}(\Omega))$.
Since $G_\Omega$ is abelian and acts transitively on $\Omega$, the stabilizers of all
components coincide. Their intersection is the kernel of the action,
which is trivial; hence every stabilizer is trivial and
$|G_\Omega|=|\Omega|$. Since $S$ is reducible, $|\Omega|\in\{2,3,4,6\}$.
Therefore $G_\Omega$ is cyclic or isomorphic to $(\mathbb Z/2)^2$, and hence
$G_\Omega\subset\Cr_1(\mathbb C)$. By Lemma~\ref{lem-diagonal-in-product},
the kernel of $G\to G_\Omega$ embeds into the automorphism group of a rational
component of $S$, hence into $\Cr_2(\mathbb C)$. Proposition~\ref{prop-abstact-extension}
now shows that $G$ is of product type.
}
\end{proof}

\section{Reducible anti-canonical element: Lc case}
\label{sec-reducible-lc}
In this section, we treat the case of a $G$-invariant reducible element $S\in |-K_X|$ such that the pair $(X, S)$ is lc but not dlt. We study the induced action of $G$ on~$S$. The key ingredient here is the boundedness of the number of components of~$S$ that pass through a terminal point on a threefold $X$. This is done in Proposition~\ref{prop-lc-bounded-number-of-components}. 

\begin{remark}
\label{rem-canonical-surface-boundedness}
Let $P\in X$ be a $2$-dimensional canonical singularity germ. 
Assume that $(X, S=\sum_{i=1}^n S_i)$ is an lc pair such that $P \in S_i$ for every $i$. 
Then using the classification of canonical singularities on surfaces, one can check that $n\leq 2$. 
\end{remark}

\begin{example}
Let $P\in X$ be the germ of an ordinary double point given by the equation 
\[
x_1 x_2 + x_3 x_4 = 0
\]
where $P$ is the origin. 
Consider a divisor $S$ on $X$ given by the equation $x_1x_2=0$. Then $S=S_1+S_2+S_3+S_4$ has four components passing through $P$. The pair $(X, S)$ is toric; in particular, it is lc.
\end{example}

\begin{proposition}
\label{prop-lc-bounded-number-of-components}
Let $P\in X$ be a $3$-dimensional terminal singularity germ. 
Let $(X, S=\sum_{i=1}^n S_i)$ be an lc pair such that $P \in S_i$ for every $i$. 
Then $n\leq 4$. 
\end{proposition}

\begin{proof}
Let 
\[
\pi\colon \widetilde{P}\in\widetilde{X}\lra X\ni P
\]
be the canonical cover as in Section~\ref{subsec-terminal-sing}, and let $(\widetilde{X}, \widetilde{S})$ be the log pullback of $(X, S)$.
We see that the number of components of $S$ is bounded by the number of components of $\widetilde{S}$.  
So we may assume that $P\in X$ is Gorenstein. 

Clearly, if $X$ is smooth, then $n\leq 3$. 
Hence we may assume that $X$ is singular at $P$. Let us assume that $0=P\in X\subset \mathbb{C}^4$ is given by the equation
\[
\phi(x_1,x_2,x_3,x_4)=\sum_{i\geq 2}\phi_i(x_1,x_2,x_3,x_4)=0, 
\]
where $\deg \phi_i = i$. Since the point $P\in X$ is terminal, we have $1\leq \rk(\phi_2)\leq 4$. Consider the blow-up $g\colon \widetilde{\mathbb{C}^4} \to \mathbb{C}^4$ of the closed point $P\in X\subset \mathbb{C}^4$ with  exceptional divisor $E\simeq \mathbb{P}^3$. 
Let $Y$ be the strict transform of $X$ on $\widetilde{\mathbb{C}^4}$. 
By abuse of notation, we denote by $g\colon Y\to X$ the induced blow-up of $X$. 
Since $K_X+S$ is $\mathbb{Q}$-Cartier and~$K_X$ is Cartier on $X$, we have that~$S$ is $\mathbb{Q}$-Cartier, and hence Cartier; see \cite[Lemma 5.1]{Kaw88}. 
We denote by $D$ a Cartier divisor on $\mathbb{C}^4$ such that~$D|_X=S$. Write
\begin{equation}
K_{\widetilde{\mathbb{C}^4}} = g^* K_{\mathbb{C}^4} + 3 E, \quad 
Y = g^* X - 2 E,\quad 
D' = g^*D - b E, 
\end{equation}
where $D'$ is the strict transform of $D$ on $\widetilde{\mathbb{C}^4}$. Put $E_Y=E|_Y$. 
Since the pair $(X, S)$ is lc, for the log discrepancy of~$E_Y$, we have $a(E_Y, X, S) = 2-b\geq 0$, and so 
$b\leq 2$. 

Note that $E|_E = -H_E$, where $H_E$ is a hyperplane on $E\simeq \mathbb{P}^3$.
Write 
\[
0 = g^*D|_E = D'|_E + bE|_E = D'|_E - b H_E.
\]

We now show that if $b=1$, then $n\leq 2$. 
Assume that $b=1$. Then $D$ is smooth at the point~$P$. We claim that $X\cap D=S$, considered as a Cartier divisor on $D$, is given by an equation with a non-zero quadratic term in the expansion near $P$. Indeed, for the discrepancy of $E_Y$, we have 
\[
1=2-b=a(E_Y, X, S) = a(E_Y, D, X|_D)=3-\mult_P (X|_D),
\]  
and so $\mult_P (X|_D)=2$. 
Hence $P\in S=X|_D$, considered as a surface germ in a smooth $3$-dimensional space~$D$, has at most two components meeting at $P$. This shows that $n\leq 2$ in this case. 

   So we may assume that $b=2$ and the log pullback via $g$ of the pair $(\mathbb{C}^4, X+ D)$ is $(\widetilde{\mathbb{C}^4}, Y+ D'+E)$.
By assumption, the pair $(X, S)$ is lc. By inversion of adjunction, it follows that $(\mathbb{C}^4, X+ D)$ is lc. Hence, $(\widetilde{\mathbb{C}^4}, Y+ D'+E)$ is lc as well. Thus, the intersection $Y\cap D'\cap E$ cannot contain a surface. Therefore, $ D'|_Y$ is equal to the strict transform $S'=\sum S'_i $ of $S$ on $Y$. Also, again by the lc condition, $Y$ is normal along $E\cap Y$, and hence~$Y$ is normal by the Serre criterion.  
By adjunction, the pair $(Y, S'+E_Y)$ is lc.

We can write
\[
S'|_{E_Y} = D'|_{E_Y} \sim b H_E|_{E_Y}=2 H_E|_{E_Y}.
\] 
Observe that the intersection $S'_i\cap E_Y$ is a curve for every $i$. Indeed, this follows from the fact that $E_Y$ is a Cartier divisor and $S'_i\cap E_Y$ is non-empty for every $i$, because $P\in S_i$ by assumption.

Since $S$ is Cartier on $X$, we have that $S' = g^*S-2E_Y$ is Cartier on $Y$.   
Recall that the pair $(Y, S'+E_Y)$ is lc, $Y$ is Gorenstein, and both $S'$ and $E_Y$ are Cartier. It follows that $Y$ is smooth at the generic point of $S'\cap E_Y$. Therefore, two components of $S'$ cannot intersect along a curve that belongs to $E_Y$. Thus, it is enough to show that $S'|_{E_Y}$ has at most four components on $E_Y$.
We consider four cases.

\emph{Case}~$\rk(\phi_2)=4$.~
 Then $P\in X$ is an ordinary double point, $Y$ is smooth, and $E_Y$ is isomorphic to $\mathbb{P}^1\times\mathbb{P}^1$. Since $H_E|_{E_Y}\sim(1,1)$, we obtain
\[
S'|_{E_Y}\sim (2,2).
\]
Since $S'|_{E_Y}$ can be decomposed into at most four components, it follows that $S'$ (and~$S$) has at most~four components, so $n\leq 4$.

\emph{Case}~$\rk(\phi_2)=3$.~
Then $P\in X$ is point of type $cA_1^k$ for some $k\geq 2$. 
Note that $Y$ is either smooth or terminal and has exactly one singularity of type $cA_1$. 
Also, $E_Y$ is isomorphic to $\mathbb{P}(1,1,2)$.
Since $H_E|_{E_Y}\sim 2 L$, where $L$ is the positive generator of the class group of $\mathbb{P}(1,1,2)$, we obtain
\[
S'|_{E_Y} \sim 4 L.
\]
By adjunction, the pair $(E_Y, S'|_{E_Y})$ is lc. It follows that~$S'|_{E_Y}$ has at most three components. Hence $S'$ (and~$S$) has at most three components, so $n\leq 3$.

\emph{Case}~$\rk(\phi_2)=2$.~  Then $P\in X$ is point of type $cA_k$ for some $k\geq 2$, and $E_Y$ is reducible: $E_Y=E_1+E_2$, where $E_i\simeq \mathbb{P}^2$ and $E_1$ intersects $E_2$ in a line. Write 
\[
0 = g^*D|_{E_i} = D'|_{E_i} - 2 H_E|_{E_i}; 
\]
hence $S'|_{E_i}\sim 2 H_E|_{E_i}$.
Thus, the restriction of $S'$ to $E_i$ is a (possibly reducible) conic. It follows that~$S'|_{E_i}$ has at most two components. Hence $S'$ (and~$S$) has at most four components, so $n\leq 4$.

\emph{Case}~$\rk(\phi_2)=1$.~
Then $P\in X$ is point of type $cD$ or $cE$, and $E_Y$ is non-reduced. We have $E_Y = 2 E_1$, where $E_1\simeq \mathbb{P}^2$. 
However, this contradicts the fact that the pair $(Y, S' + E_Y)$ is lc, so this case does not occur. 
This completes the proof.
\end{proof}

\begin{corollary}
\label{prop-lc-bounded-number-of-components-gorenstein}
Let $P\in X$ be a $3$-dimensional terminal singularity germ. 
Let $(X, S=\sum_{i=1}^4 S_i)$ be an lc pair such that $P \in S_i$ for every $i$. 
Assume that for the canonical cover $\pi\colon \widetilde{X}\to X$, we have that $\widetilde{P}\in\widetilde{X}$ is an ordinary double point. Then $P\in X$ is Gorenstein, so $X=\widetilde{X}$ and $P\in X$ is an ordinary double point.
\end{corollary}

\begin{proof}
We use the same notation as in the proof of Proposition~\ref{prop-lc-bounded-number-of-components}.
Assume that $P\in X$ is non-Gorenstein of index $r>1$. Since $n=4$, it follows that $S$ has four components meeting at~$P$, and so $\widetilde{S}$ has four components meeting at $\widetilde{P}$.
By assumption, $\widetilde{P}\in\widetilde{X}$ is an ordinary double point.  
By the classification of terminal singularities, see \cite[Section~6.1]{Re85}, for the Gorenstein index $r$ of $P\in X$, we have $r=2$, and the Galois group $\mathbb{Z}/2$ of the covering $\pi$ interchanges two rulings of the normal cone to $\widetilde{P}\in \widetilde{X}$, which is isomorphic to $\mathbb{P}^1\times\mathbb{P}^1$. Thus the action of $\mathbb{Z}/2$ interchanges some components of~$\widetilde{S}$ (\textit{cf.}~the case $\rk(\phi_2)=4$ in the proof of Proposition~\ref{prop-lc-bounded-number-of-components}), so the number of components of $S$ cannot be equal to $4$. This contradiction shows that $P\in X$ is Gorenstein, and hence an ordinary double point. 
\end{proof}

\begin{lem}
\label{lem-quadric-cyclic-permutation}
Let $C=\sum_{i=1}^4 C_i$ be a toric boundary on $\mathbb{P}^1\times\mathbb{P}^1$.  Let $G\subset \Aut(\mathbb{P}^1\times\mathbb{P}^1, C)$ be a finite abelian group.  Assume that the image of\, $G$ in the symmetric group on the set $\{C_i\}$ contains a cycle of length~$4$.  Then either $\mathfrak{r}(G)\leq 1$ or $G\subset \mathbb{Z}/2\times \mathbb{Z}/4$.
\end{lem}

\begin{proof}
The group $\Aut(\mathbb{P}^1 \times \mathbb{P}^1, C)$ is just $(\mathbb{C}^*)^2 \rtimes D_4$, where $D_4$ is the dihedral group acting on the components of $C$. As $G$ is abelian and its image in $D_4$ contains an element of order $4$, its image in $D_4$ is $\mathbb{Z}/4$. The kernel is a subgroup of $\mathbb{C}^\times \times \mathbb{C}^\times$ fixed by the element of order $4$ given by $(x, y) \mapsto (y^{-1},x)$. Thus, the kernel has order at most $2$. It follows that either $\mathfrak{r}(G)\leq 1$ or $G\subset \mathbb{Z}/2\times \mathbb{Z}/4$. 
\end{proof}

Until the end of this section, we work in the following setting.  

\begin{setting}
Let $X$ be a $G\mathbb{Q}$-Fano threefold, where $G$ is a finite abelian group that faithfully acts on $X$. Assume that $S\in|-K_X|$ is a $G$-invariant element such that the pair $(X, S)$ is lc and $S=\sum S_i$ is reducible. 
\end{setting}

\begin{proposition}
\label{prop-lc-gorenstein-4comp}
Assume that $G$ interchanges two points $P_1$ and $P_2$ on $X$ such that four components of\, $S=\sum S_i$ meet at each $P_i$.  Assume further that $P_i\in X$ is an ordinary double point. Then $G$ is of product type.
\end{proposition}

\begin{proof}
Let $H$ be the subgroup of $G$ of index $2$ that stabilizes $P_i$, so we have an exact sequence
\begin{equation}
\label{eq-last-sequence}
0\lra H \lra G\lra \mathbb{Z}/2 \lra 0. 
\end{equation}
Fix $i\in\{1,2\}$. By assumption, four components of $S$ meet at $P_i$.
As in the proof of Proposition~\ref{prop-lc-bounded-number-of-components}, for the blow-up $g\colon Y\to X$ of $P_i$ with exceptional divisor $E_Y\simeq \mathbb{P}^1\times\mathbb{P}^1$ on $Y$, the log pullback of the pair $(X, S)$ is the pair $(Y, S' + E|_Y)$, where $S'$ is the strict preimage of $S$. 
Put $C=S'|_{E_Y}$. Then ${C}$ is a reducible divisor on $\mathbb{P}^1\times\mathbb{P}^1$ of bidegree $(2,2)$. 
We have ${C}=C_1+C_2+C_3+C_4$, where each $C_i$ is isomorphic to $\mathbb{P}^1$.

Let $\mathcal{D}(S)$ be the dual complex of $(X,S)$. 
Since ${C}$ admits a $0$-dimensional stratum, it follows that $\mathcal{D}(S)$ has dimension $2$. By Corollary~\ref{not-a-circle}, we see that $\mathcal{D}(S)$ is homeomorphic to $\mathbb{S}^2$. 

By Lemma~\ref{lem-transitive-action}, we may assume that $G$ acts transitively on the set $\Omega$ of components of $S$. For $i=1,2$, let $\Omega_i$ be the set of the four components of $\Omega$ passing through $P_i$. Since $\Omega_1\cap\Omega_2$ and $\Omega_1\cup\Omega_2$ are $G$-invariant, either
\[
\Omega_1=\Omega_2=\Omega,\quad |\Omega|=4,
\qquad\text{or}\qquad
\Omega=\Omega_1\sqcup\Omega_2,\quad |\Omega|=8.
\]

Assume that $|\Omega|=4$, and put $G_\Omega=\operatorname{im}(G\to\operatorname{Sym}(\Omega))$. The action of $G_\Omega$ on $\Omega$ is faithful and transitive. Since $G_\Omega$ is abelian, $|G_\Omega|=|\Omega|=4$, and hence $G_\Omega\subset\Cr_1(\mathbb C)$. The components of $S$ are rational by Lemma~\ref{lem-strata-are-rational}, applied on a dlt modification, and the induced action of $G$ on $S$ is faithful by Lemma~\ref{lem-fixed-curve-surface}. Thus Lemma~\ref{lem-diagonal-in-product} and Proposition~\ref{prop-abstact-extension} show that $G$ is of product type.

Assume that $|\Omega|=8$. Then $H$ acts transitively on each $\Omega_i$. Put $H_\Omega=\operatorname{im}(H\to\operatorname{Sym}(\Omega))$. Since $G$ is abelian, the action of $H_\Omega$ on $\Omega_1$ is faithful and transitive. Thus
\[
H_\Omega\simeq\mathbb Z/4
\quad\text{or}\quad
H_\Omega\simeq(\mathbb Z/2)^2.
\]
Since $C$ has four components and bidegree $(2,2)$, it is a toric boundary on $E_Y$. Let $H_1$ be the image of $H$ in $\Aut(E_Y,C)$, and put $K=\ker(H\to H_1)$. At a general point of $E_Y$, Lemma~\ref{lem-faithful-action} shows that $K$ acts faithfully on the normal line to $E_Y$; hence $K$ is cyclic. Moreover, $H_\Omega=\operatorname{im}(H_1\to\operatorname{Sym}(\{C_1,\ldots,C_4\}))$.

If $H_\Omega\simeq\mathbb Z/4$, Lemma~\ref{lem-quadric-cyclic-permutation} shows that $H_1$ is cyclic or $H_1\subset\mathbb Z/2\times\mathbb Z/4$. If $H_\Omega\simeq(\mathbb Z/2)^2$, the same argument as in the proof of Lemma~\ref{lem-quadric-cyclic-permutation} shows that $\ker(H_1\to H_\Omega)$ has order at most~$2$. Thus, in all cases, $H_1$ is either cyclic or a $2$-group of order at most~$8$.

The quotient $G/K$ fits into the exact sequence
\[
0\lra H_1\lra G/K\lra\mathbb Z/2\lra0.
\]
Hence $G/K$ either has rank at most~$2$ or is an abelian $2$-group of order at most~$16$. Such groups have type~(1),~(2), or~(5) in Theorem~\ref{thm-surfaces}, and therefore belong to $\Cr_2(\mathbb C)$. Since $K$ is cyclic, Proposition~\ref{prop-abstact-extension}, applied to
\[
0\lra K\lra G\lra G/K\lra0,
\]
shows that $G$ is of product type.
\color{black}
\end{proof}

\begin{proposition}
\label{prop-lc-gorenstein-4comp-2}
Assume that $G$ interchanges two points $P_1$ and $P_2$ on $X$ such that four components of\, $S=\sum S_i$ meet at each $P_i$.  Assume that $P_i\in X$ is Gorenstein and the local equation of\, $P_i\in X$ has a quadratic part of rank $2$. Then $G$ is of product type.
\end{proposition}

\begin{proof}
We use the same notation as in the proof of Proposition~\ref{prop-lc-bounded-number-of-components}. 
Let $H$ be the subgroup of $G$ of index $2$ that stabilizes $P_i$. 
Fix $i\in\{1,2\}$. 
By assumption, four components of $S$ meet at $P_i$. 
Since the local equation of $P_i\in X$ has a quadratic part of rank $2$, for the blow-up $g\colon Y\to X$ of $P_i$ with exceptional divisor $E_Y$, we have $E_Y=E_1+E_2$, where $E_i\simeq \mathbb{P}^2$. 

Arguing as in the proof of Proposition~\ref{prop-lc-bounded-number-of-components}, we have that the log pullback $(Y, S' + E_1+E_2)$ of $(X, S)$ is lc, where $S'$ is the strict preimage of $S$. By adjunction, the pair $(E_1, S'|_{E_1}+E_2|_{E_1}+D_{E_1})$ is lc as well, where $D_{E_1}$ is the different (in fact, from the proof of Proposition~\ref{prop-lc-bounded-number-of-components}, it follows that the different is equal to $0$, but we do not need this). Note that the divisor $E_1$ might not be $\mathbb{Q}$-Cartier; however, the adjunction formula still makes sense. 

As shown in the proof of Proposition~\ref{prop-lc-bounded-number-of-components}, the divisors $E_1$ and $E_2$ intersect in a line, and $S'|_{E_1}$ is a reducible conic, that is, the union of two lines. 
Since the pair $(E_1, S'|_{E_1}+E_2|_{E_1}+D_{E_1})$ is lc, these lines are distinct and do not intersect at a point that belongs to $E_1\cap E_2$. 
Similarly, $S'|_{E_2}$ is a union of two distinct lines whose point of intersection does not belong to $E_1\cap E_2$.  
Put $S'\cap E_1=L_{1}\cup L_{2}$ and $S'\cap E_2=L_{3}\cup L_{4}$. 
Up to renumbering of $L_i$, we may assume that
\[Q_{1,2}=L_1\cap L_2\in E_1\setminus E_2,  \quad Q_{3,4}=L_3\cap L_4\in E_2\setminus E_1,
\] 
\[
Q_{2,3}=L_2\cap L_3\in E_1\cap E_2,  \quad Q_{1,4}=L_1\cap L_4\in E_1\cap E_2.
\] 
By Lemma~\ref{lem-transitive-action}, either $G$ is of product type, or the action of $G$ on the set of components of $S$ is transitive. We may assume that the latter case is realized. 

Let $\Omega$ be the set of components of $S$, and let $\Omega_i$ be the set of the four components through $P_i$. As in the proof of Proposition~\ref{prop-lc-gorenstein-4comp}, either
\[
\Omega_1=\Omega_2=\Omega,\quad |\Omega|=4,
\qquad\text{or}\qquad
\Omega=\Omega_1\sqcup\Omega_2,\quad |\Omega|=8.
\]
In the first case, $Q_{1,2}$ is a zero-dimensional stratum of the pair $(Y,S'+E_1+E_2)$, so the argument for $|\Omega|=4$ in the proof of Proposition~\ref{prop-lc-gorenstein-4comp} shows that $G$ is of product type.

Assume that $|\Omega|=8$. Put $H_\Omega=\operatorname{im}(H\to\operatorname{Sym}(\Omega))$. Since $G$ is abelian, the action of $H_\Omega$ on $\Omega_1$ is faithful and transitive. Put $N=\ker(H\to H_\Omega)$. Then
\[
H_\Omega\simeq\mathbb Z/4
\quad\text{or}\quad
H_\Omega\simeq(\mathbb Z/2)^2.
\]
The group $H_\Omega$ preserves the sets $\{Q_{1,2},Q_{3,4}\}$ and $\{Q_{2,3},Q_{1,4}\}$. If $H_\Omega\simeq\mathbb Z/4$, its generator would interchange these two sets. Hence $H_\Omega\simeq(\mathbb Z/2)^2$.

The points $Q_{1,2},Q_{3,4},Q_{2,3},Q_{1,4}$ determine four coordinate axes in $T_{P_1}X$; choose the corresponding coordinates $x_1,\ldots,x_4$ in this order. Since $H$ is abelian, $N$ acts on $x_1,x_2$ by the same character $\alpha$ and on $x_3,x_4$ by the same character $\beta$. This action is faithful by Lemma~\ref{lem-faithful-action}. Let $\phi$ be an $N$-semi-invariant local equation of $P_1\in X$. Its quadratic part is $x_1x_2$, so $\phi$ has weight $2\alpha$. Since $P_1$ is isolated, $\phi\notin(x_1,x_2)^2$. Thus $\phi$ contains a monomial of degree at most one in $x_1,x_2$. Comparing its weight with $2\alpha$ gives $2\alpha\in\langle\beta\rangle$ or $\alpha\in\langle\beta\rangle$.
Since the action is faithful, these relations imply that $N$ contains a cyclic subgroup $K$ of index at most~$2$. Then $G/K$ is an abelian $2$-group of order at most~$16$, and hence $G/K\subset\Cr_2(\mathbb C)$ by Theorem~\ref{thm-surfaces}. Since $K$ is cyclic, Proposition~\ref{prop-abstact-extension} shows that $G$ is of product type.
\color{black}
\end{proof}

\begin{corollary}
\label{prop-lc-gorenstein-4comp-3}
Assume that $G$ interchanges two points $P_1$ and $P_2$ on $X$ such that four components of\, $S=\sum S_i$ meet at each $P_i$.  Assume that each $P_i\in X$ is not Gorenstein.  Also assume that the local equation of\, $\widetilde{P}_i\in \widetilde{X}$ in $\mathbb{C}^4$ has a quadratic part of rank $2$.  Then $G$ is of product type.
\end{corollary}

\begin{proof}
By Lemma~\ref{lem-transitive-action}, either $G$ is of product type, or the action of $G$ on the set of components of $S$ is transitive. We may assume that the latter case is realized. 

Let $H$ be the subgroup of $G$ of index~$2$ stabilizing $P_1$ and
$P_2$, and let $\Omega_i$ be the set of components of $S$ passing
through $P_i$. As in Propositions~\ref{prop-lc-gorenstein-4comp}
and~\ref{prop-lc-gorenstein-4comp-2}, either
\[
\Omega_1=\Omega_2=\Omega,\qquad |\Omega|=4,
\qquad\text{or}\qquad
\Omega=\Omega_1\sqcup\Omega_2,\qquad |\Omega|=8.
\]

Assume first that $|\Omega|=4$. The four components have unique lifts
to the index-one cover by
Proposition~\ref{prop-lc-bounded-number-of-components}, and the deck
group fixes these lifts. Let $\widetilde Y\to\widetilde X$ be the
blow-up of $\widetilde P_1$. This morphism is equivariant with respect
to the deck group, and its quotient gives a birational morphism
$f\colon Y\to X$. Let $S_Y$ be determined by
\[
K_Y+S_Y=f^*(K_X+S).
\]
Then $(Y,S_Y)$ has a zero-dimensional stratum. Thus
$\mathcal D(S)\simeq\mathbb S^2$, and the argument for $|\Omega|=4$
in Proposition~\ref{prop-lc-gorenstein-4comp} applies.

Assume that $|\Omega|=8$. Put $H_\Omega=\operatorname{im}(H\to\operatorname{Sym}(\Omega))$. Since $G$ is abelian, the action of $H_\Omega$ on $\Omega_1$ is faithful and transitive. Put $N=\ker(H\to H_\Omega)$. Let $\widetilde\Omega_1$ be the set of the four components through $\widetilde P_1$ on the index-one cover. Each of them is the unique lift of a component in $\Omega_1$ and is fixed by the deck group. In the notation of
Proposition~\ref{prop-lc-gorenstein-4comp-2}, the group $H_\Omega$
preserves both sets $\{Q_{1,2},Q_{3,4}\}$ and
$\{Q_{2,3},Q_{1,4}\}$.
If $H_\Omega\simeq\mathbb Z/4$, its generator would interchange these
two sets. Hence $H_\Omega\simeq(\mathbb Z/2)^2$.
Let $\widetilde H$ be the lift of $H$, and put
$\widetilde N=\ker(\widetilde H\to\operatorname{Sym}(\widetilde\Omega_1))$. The
sequence~\eqref{eq-lift} restricts to
\[
0\lra\mathbb Z/r\lra\widetilde N\lra N\lra0,
\]
where $r$ is the Gorenstein index of $P_1\in X$. If $\widetilde H$ is
abelian, the calculation in the proof of
Proposition~\ref{prop-lc-gorenstein-4comp-2} shows
that $\widetilde N$ contains a cyclic subgroup of index at most~$2$.
Its image in $N$ is cyclic and also has index at most~$2$.

Suppose that $\widetilde H$ is not abelian. Then the index of
$P_1\in X$ is~$2$ by
Proposition~\ref{prop-action-on-non-gorenstein-points}, and a Sylow
$2$-subgroup $\widetilde H_2$ is abelian by the proof of
Proposition~\ref{prop-terminal-point-not-abelian}. Every odd Sylow
subgroup $\widetilde H_p$ commutes with $\widetilde H_2$: otherwise an
odd-order element $x$ would be conjugate to $cx$, where $c$ is the
deck involution, although $cx$ has even order. The calculation in the
proof of Proposition~\ref{prop-lc-gorenstein-4comp-2}, applied to
$\widetilde H_2$ and to
$\langle\widetilde H_2,\widetilde H_p\rangle$ shows that $N_2$ is
cyclic or has type $[a,1]$, while $N_p$ is cyclic for every odd
prime~$p$.

Thus, in either case, $N$ contains a cyclic subgroup $K$ of index at
most~$2$. Consequently, $G/K$ is an abelian $2$-group and
\[
|G/K|=[G:H][H:N][N:K]\leq2\cdot4\cdot2=16.
\]
Hence $G/K\subset\Cr_2(\mathbb C)$ by
Theorem~\ref{thm-surfaces}. Since $K$ is cyclic,
Proposition~\ref{prop-abstact-extension} shows that $G$ is of product
type.
\color{black}
\end{proof}

\begin{proposition}
\label{prop-reducible-4}
Let $G$ be a finite abelian group. 
Assume that 
\begin{itemize}
\item
$X$ is a $G\mathbb{Q}$-Fano threefold and the action of\, $G$ on $X$ is faithful,
\item 
$S\in |-K_X|$ is a reducible $G$-invariant element,
\item 
the pair $(X, S)$ is lc and not dlt, and
\item
the dual complex $\mathcal{D}(S)$ is homeomorphic to a sphere $\mathbb{S}^2$. 
\end{itemize}
Then $G$ is of product type.
\end{proposition}

\begin{proof}
{
Let $f\colon(Y,S_Y)\to(X,S)$ be a $G$-equivariant dlt modification.
The action of $G$ on $S$ is faithful by
Lemma~\ref{lem-fixed-curve-surface}. By
Lemma~\ref{lem-transitive-action}, we may assume that $G$ acts
transitively on the components of $S$. All strata of $S_Y^{=1}$, and
in particular the strict transforms of these components, are rational
by Lemma~\ref{lem-strata-are-rational}. If the action on
$\mathcal D(S_Y^{=1})\simeq\mathbb S^2$ has a fixed point,
Remark~\ref{rem-fixed-point-on-dual-complex} gives a $G$-invariant
rational stratum. Its image on $X$ is rationally connected, so
Corollary~\ref{cor-fixed-point} or
Corollary~\ref{cor-inv-curve-or-surface} applies.

Otherwise, Lemma~\ref{lem-action-on-spheres} gives two strata whose
images $C_1,C_2\subset X$ are interchanged by $G$. If $C_1=C_2$, then
$C_1$ is $G$-invariant, so Corollary~\ref{cor-fixed-point} or
Corollary~\ref{cor-inv-curve-or-surface} applies. Assume that they are distinct. If
$\dim C_i>0$, or if $C_i=P_i$ are points and at most three components
of $S$ pass through each $P_i$, the components of $S$ containing $C_1$ or
$C_2$ form a $G$-invariant set of at most six elements and hence
exhaust $S$. Let $\Omega$ be the set of components of $S$, and put
$G_\Omega=\operatorname{im}(G\to\operatorname{Sym}(\Omega))$. The group $G_\Omega$ acts faithfully and transitively on
$\Omega$. Since $G_\Omega$ is abelian, $|G_\Omega|=|\Omega|$;
hence $G_\Omega$ is cyclic or isomorphic to $(\mathbb Z/2)^2$. Thus
Lemma~\ref{lem-diagonal-in-product} and
Proposition~\ref{prop-abstact-extension} show that $G$ is of product
type in this case.

Thus exactly four components of $S$ pass through each $P_i$.
Proposition~\ref{prop-lc-bounded-number-of-components} and
Corollary~\ref{prop-lc-bounded-number-of-components-gorenstein} reduce
the remaining cases to Propositions~\ref{prop-lc-gorenstein-4comp},
~\ref{prop-lc-gorenstein-4comp-2}, and
Corollary~\ref{prop-lc-gorenstein-4comp-3}. The claim follows.
}
\end{proof}

\begin{proposition}
\label{prop-reducible-3}
Let $G$ be a finite abelian group. 
Assume that 
\begin{itemize}
\item
$X$ is a $G\mathbb{Q}$-Fano threefold and the action of\, $G$ on $X$ is faithful,
\item 
$S\in |-K_X|$ is a reducible $G$-invariant element,
\item 
the pair $(X, S)$ is lc and not dlt, and
\item
the dual complex $\mathcal{D}(S)$ is homeomorphic to an interval. 
\end{itemize}
Then $G$ is of product type.
\end{proposition}
\begin{proof}
{
Let $f\colon (Y,S_Y)\to (X,S)$ be a $G$-equivariant dlt modification,
and write $B=S_Y^{=1}=\sum B_i$. Every double curve of $B$, that is, every irreducible component of $B_i\cap B_j$ with $i\ne j$, is rational by
\cite[Claim~32.1]{KX16}. Since $K_Y+B\sim0$ and $Y$ is rationally
connected, the anti-canonical exact sequence, as in the proof of
Proposition~\ref{the-pair-is-plt}, gives $H^1(B,\OOO_B)=0$. Applying
the exact sequence~\eqref{eq-normalization-conductor} successively,
removing an end component of the interval $\mathcal D(B)$ at each step,
and using the rationality of the double curves, we obtain
$H^1(B_i,\OOO_{B_i})=0$ for every $i$. Since $\mathcal D(B)$ is an
interval, every component $B_i$ meets another component of $B$ along
a curve. Adjunction gives
$K_{B_i}+\Delta_i=(K_Y+B)|_{B_i}\sim0$. Since $\Delta_i$ contains
this curve, we have $\Delta_i\ne0$. Lemma~\ref{lem-2-dim-CY-pairs}
shows that $B_i$ is rational unless it is birationally ruled over an
elliptic curve; the latter case would give $H^1(B_i,\OOO_{B_i})\ne0$.
Thus all the components $B_i$ are rational.

The action on $S$ is faithful by Lemma~\ref{lem-fixed-curve-surface}.
By Lemma~\ref{lem-transitive-action}, we may assume that $G$ acts
transitively on the components of $S$. Hence $G$ acts transitively on
the vertices of $\mathcal D(B)$ corresponding to their strict transforms. A finite
group of homeomorphisms of an interval has order at most~$2$.
Since $S$ is reducible, it follows that $S=S_1+S_2$, and its two
rational components are interchanged. Lemma~\ref{lem-diagonal-in-product}
and Proposition~\ref{prop-abstact-extension} finish the proof.
}
\end{proof}

\section{Groups of K3 type}
\label{sect-action-on-K3}
In this section, we work with finite abelian groups of K3 type. Recall that by a K3 surface we mean a normal projective surface
$S$ with at worst canonical (that is, $1$-lc) singularities such that
$\mathrm{H}^1(S, \OOO_S)=0$ and $K_S\sim 0$. 

\begin{defin}
We say that a finite abelian subgroup $G$ is of \emph{K3 type} if $G$ faithfully acts on a $G\mathbb{Q}$-Fano threefold $X$ preserving a K3 surface $S\in|-K_X|$ with at worst du Val singularities.
\end{defin}

Let $X$ be a $G\mathbb{Q}$-Fano threefold, where $G$ is a finite abelian group and $G$ acts on $X$ faithfully. Assume that there exists a $G$-invariant element $S\in |-K_X|$ that is a K3 surface with at worst canonical singularities. Consider the  exact sequence
\[
0\lra K \lra G\lra H\lra 0, 
\]
where $H\subset \Aut(S)$ and $K$ fixes $S$ pointwise.  Note that $K=\mathbb{Z}/m$ for some $m\geq 1$ (\textit{cf.}~Corollary~\ref{cor-G-preserves-K3-type}).  In particular, the group $G$ is of K3 type.  Let $Y=X/K$ be the quotient, and let $f\colon X\to Y$ be the natural projection.

\begin{lem}
\label{lem-for-thm-2}
Assume that $K$ is non-trivial. 
Then the action of\, $K$ is free in codimension~$1$ outside $S$. 
Hence, the divisorial part of the ramification of $f$ is precisely $S$. 
\end{lem}

\begin{proof}
Let $K'$ be a non-trivial subgroup of $K$. 
Assume there exists a $K'$-fixed divisor $D$ such that $D\neq S$. 
Then $K'$ fixes $S\cup D$ pointwise. 
Since $S$ is ample, we have $S\cap D\neq\emptyset$. Since $S$ is $\mathbb{Q}$-Cartier, there exists a curve $C\subset S\cap D$. At a general point of $C$, the variety $X$ is smooth, so the tangent space to $S\cup D$ at some point $P\in C$ is equal to the tangent space to $X$. However, this contradicts Lemma~\ref{lem-faithful-action} applied to the group~$K'$.
\end{proof}

\begin{proposition}
\label{prop-K3-boundedness}
The threefold $Y=X/K$ is an $H\mathbb{Q}$-factorial Fano variety with at worst canonical singularities. Moreover, $-K_Y\sim mB$ for some $B\in \Cl(Y)$.
\end{proposition}

\begin{proof}
Recall that we have $K=\mathbb{Z}/m$ for some $m\geq 1$. We may assume that $m>1$. 
First we show that $Y$ is $H\mathbb{Q}$-factorial. Indeed, this follows from the fact that
\[
\Cl(Y)\otimes\mathbb{Q} = (\Cl(X)\otimes\mathbb{Q})^{K}
\]
and hence 
\[
(\Cl(Y)\otimes\mathbb{Q})^{H} = (\Cl(X)\otimes\mathbb{Q})^{G} = (\Pic(X)\otimes\mathbb{Q})^{G}.
\]

Put $D=f(S)$. Since $K$ fixes $S$ pointwise, we have $S\simeq D$. The Hurwitz formula yields
\[
0\sim K_X + S = f^*( K_Y + D ),
\]
and hence $m(K_Y + D)\sim 0$. Note that $Y$ is a Fano variety with klt singularities. By inversion of adjunction, the pair $(Y,D)$ is plt. Observe that the divisor $D$ is $\mathbb{Q}$-Cartier.
Put $T=K_Y+D$. We show that $T\sim0$. Assume otherwise. Since $mT\sim0$, the divisor $T$ is non-trivial torsion, and hence $\mathrm H^0(Y,\OOO_Y(T))=0$. By adjunction (\textit{cf.}~\cite[Corollary~5.25]{KM98}), there is an exact sequence
\[
0\lra\OOO_Y(K_Y)\lra\OOO_Y(T)\lra\OOO_D(K_D)\lra0.
\]
Kawamata--Viehweg vanishing and Serre duality give $\mathrm H^1(Y,\OOO_Y(K_Y))=0$. It follows that $\mathrm H^0(D,\OOO_D(K_D))=0$, which contradicts $K_D\sim0$. Thus $K_Y+D\sim0$.

We have shown that $-K_Y\sim D$. Since $K_Y+D$ is Cartier and $(Y,D)$ is plt, the log discrepancy $a(E,Y,D)$ of every exceptional divisor $E$ over $Y$ is a positive integer. Therefore,
\[
a(E,Y)=a(E,Y,D)+\operatorname{ord}_E(D)\geq1,
\]
so $Y$ has canonical singularities.
 
We claim that there exists a Weil divisor $B$ such that $mB\sim D$. Indeed, let $Y^0$ be a smooth locus of $Y$, and let $X^0$ be its preimage via the map $f$. We denote the corresponding map by $f^0\colon X^0\to Y^0$. Put $D^0=D\cap Y^0$. Then applying \cite[Theorem 1.2]{Wa68}, we see that there exists an element $B^0\in\Pic(Y^0)$ such that $mB^0\sim D^0$. Since $\Pic(Y^0)=\Cl(Y^0)=\Cl(Y)$, we see that $mB\sim D$, where $B$ is the closure of $B^0$ in $Y$. 
\end{proof}

\begin{remark}
\label{rem-bounded-index}
In the above notation, if $Y$ is terminal (this holds for example if the action of $K$ is free outside~$S$), then by \cite{Su04} we have $m\leq 19$. In general, by \cite{JL25} we have $m\leq 66$. 
\end{remark}

\begin{example}
\label{exam-k3-not-pt}
We present some examples of actions on Fano varieties in weighted projective spaces of groups of K3 type that are not of product type: 
In each case, $-K_X\sim\OOO_X(1)$, and $S=X\cap\{x_0=0\}\in|-K_X|$ is a $G$-invariant K3 surface with at worst du Val singularities.

\begin{enumerate}
\item  
Let $X$ be a hypersurface of degree $4$ in $\mathbb{P}^4$ given by the equation
\[
x_0^4 + x_1^4 + x_2^4 + x_3^4 + x_4^4 = 0,
\]
with the diagonal action of the group $G=(\mathbb{Z}/4)^4$.

\item 
Let $X$ be a hypersurface of degree $8$ in $\mathbb{P}(1, 1, 1, 2, 4)$ given by the equation
\[
x_0^8 + x_1^8 + x_2^8 + x_3^4 + x_4^2 = 0,
\]
with the diagonal action of the group $G=(\mathbb{Z}/8)^2\times\mathbb{Z}/4\times\mathbb{Z}/2$. 

\item
Let $X$ be a hypersurface of degree $6$ in $\mathbb{P}(1,1,1,2,2)$
given by the equation
\[
x_0^6 + x_1^6 + x_2^6 + x_3^3 + x_4^3 = 0,
\]
with the diagonal action of the group 
$G=(\mathbb{Z}/6)^2\times(\mathbb{Z}/3)^2$. 

\item
Let $X$ be a hypersurface of degree $6$ in $\mathbb{P}(1,1,1,1,3)$
given by the equation
\[
x_0^6 + x_1^6 + x_2^6 + x_3^6 + x_4^2 = 0,
\]
with the diagonal action of the group $G=(\mathbb{Z}/6)^3\times \mathbb{Z}/2$.  
\end{enumerate}
\end{example}

\section{Proofs of the main results}

\begin{proof}[Proof of Theorem~\ref{thm-first-main-thm}]
Let $X$ be a $3$-dimensional rationally connected variety, and let~$G$ be a finite abelian group such that $G\subset \Bir(X)$. After passing to a $G$-equivariant resolution of singularities, we may assume that~$X$ is smooth and projective, and that $G\subset \Aut(X)$. We run a $G$-equivariant minimal model program that terminates on a $G\mathbb{Q}$-Mori fiber space $\pi\colon X'\to Z$. If $\dim Z>0$, then by Corollary~\ref{cor-exclude-positive-dim-base}, we conclude that $G$ is of product type. So we may assume from the start that $X=X'$ is a $G\mathbb{Q}$-Fano threefold and $G\subset \Aut(X)$. 

Consider the linear system $|-K_X|$. If it is empty, we are in the last case of Theorem~\ref{thm-first-main-thm}. So we may assume that $|-K_X|$ is non-empty. Since the group $G$ acts on 
$\mathrm{H}^0(X, \OOO(-K_X))$ and $G$ is abelian, there exists a $G$-invariant element $S\in |-K_X|$. We analyze singularities of the pair $(X, S)$. 
By Proposition~\ref{the-pair-is-lc}, if the pair $(X, S)$ is not lc, then $G$ is of product type. So we may assume that $(X, S)$ is lc. 
By Proposition~\ref{the-pair-is-plt}, if $(X, S)$ is plt, then $G$ is of K3 type, and we are in the second case of Theorem~\ref{thm-first-main-thm}. So we may assume that $(X, S)$ is not plt.\looseness=-1 

Assume that $(X, S)$ is lc but not plt and $S$ is irreducible. By Proposition~\ref{prop-non-normal-irreducible-surface},  either $G$ is of product type, or~$S$ is normal. Hence we may assume that $S$ is normal. Similarly, by Proposition~\ref{the-pair-is-strictly-lc-irred}, we may assume that $S$ has two simple elliptic singularities $P_1$ and $P_2$ interchanged by $G$. 
Let $H$ be an index $2$ subgroup of $G$ that stabilizes $P_i$. 
Propositions~\ref{prop-elliptic-sing-short-orbit} and~\ref{prop-elliptic-sing-short-orbit-2} imply that the exceptional elliptic curve $E_i$ over $P_i\in S$ has an $H$-orbit of length at most $7$.  
We denote by $H_E$ the image of $H$ in the automorphism group of $E_i$.  
Corollary~\ref{cor-irreducible-HE} shows that either $\mathfrak{r}(H_E)\leq 1$, or $H_E$ belongs to a list of three groups. 
If $\mathfrak{r}(H_E)\leq 1$, then $G$ is of product type by Corollary~\ref{prop-rank-HE}. 
Then, we show that $G$ is of product type in the three exceptional cases. This is done in Propositions~\ref{prop-GE-2},~\ref{prop-GE-1} and~\ref{prop-E2-2} and Corollary~\ref{cor-GE-3}. 

Assume that $(X, S)$ is dlt and $S$ is reducible. Then by Corollary~\ref{not-a-circle}, the dual complex of $S$ is either an interval or a sphere $\mathbb{S}^2$ (the case of a point is not possible since we assume that $S$ is reducible). Applying Proposition~\ref{prop-reducible} in the former case and Proposition~\ref{prop-reducible-2} in the latter case, we conclude that $G$ is of product type. 

Finally, assume that $(X,S)$ is lc but not dlt and $S$ is reducible. Again, by Corollary~\ref{not-a-circle}, the dual complex of $S$ is either an interval or a sphere $\mathbb{S}^2$. Applying Proposition~\ref{prop-reducible-3} in the former case and Proposition~\ref{prop-reducible-4} in the latter case, we conclude that $G$ is of product type. Suppose that $G$ is neither of product type nor of K3 type, and let $W\to T$ be any $G\mathbb Q$-Mori fiber space with rationally connected total space and faithful $G$-action. Corollary~\ref{cor-exclude-positive-dim-base} forces $T$ to be a point; if $|-K_W|\ne\emptyset$, the preceding argument applied to $W$ shows that $G$ is of product or K3 type. Hence $|-K_W|=\emptyset$, proving the final assertion of the theorem. The theorem is proven. 
\end{proof}

\begin{proof}[Proof of Proposition~\ref{thm-second-main-thm}]
This follows from Lemma~\ref{lem-for-thm-2}, Proposition~\ref{prop-K3-boundedness}, and Remark~\ref{rem-bounded-index}.
 \end{proof}
 
Finally, Proposition~\ref{prop-abstact-extension-intro} is the same as Proposition~\ref{prop-abstact-extension}. 


\end{document}